\documentclass[11pt,a4paper]{article}
\usepackage[T1]{fontenc}
\usepackage{lmodern}
\usepackage[margin=24mm,headheight=14pt]{geometry}
\usepackage{amsmath,amssymb,amsthm,mathtools,bm}
\usepackage{graphicx,xcolor,enumitem,booktabs,longtable,array}
\usepackage{tikz}
\usetikzlibrary{arrows.meta,calc,patterns,positioning,decorations.pathreplacing}
\usepackage{hyperref}
\hypersetup{hypertexnames=false,colorlinks=true,linkcolor=blue!45!black,urlcolor=blue!45!black,citecolor=blue!45!black,pdftitle={Fourier decoupling and its applications — handwritten notes transcription}}
\usepackage{fancyhdr}
\allowdisplaybreaks[2]
\theoremstyle{definition}

\theoremstyle{remark}

\title{Lecture note on Fourier decoupling theory and its applications}
\author{Yudong Liang}
\date{}
\begin{document}
\maketitle
\begin{abstract}
This paper is based on lecture notes taken and organized during the 17th Graduate Mathematics Intensive Training Program (Spring 2025), held at the Beijing International Center for Mathematical Research, in which Professor Shaoming Guo lectured on "An Introduction to Fourier Decoupling Theory."~\cite{Guo2025} The main content includes an introduction to and proofs of $\ell^2$ decoupling inequality, together with its applications in number theory, harmonic analysis, and geometric measure theory.

The writing follows Professor Guo's teaching style, proving theorems in an intuitive and heuristic manner; consequently, the text contains many figures to explain the ideas at hand, while also lacking some rigor. However, all of these geometric intuitions can be rigorously justified, albeit in a more complicated and technical way. The author hopes that readers will be able to understand, through these lecture notes, the geometric intuition hidden behind the cumbersome proofs, and thereby deepen their understanding of the relevant definitions and theorems. Rigorous proofs can be found in Professor Guo's forthcoming book~\cite{GuoBook}, which interested readers may consult.
\end{abstract}

\tableofcontents
\clearpage

\section{Introduction}
\subsection*{Kakeya sets}
Suppose $E\subset\mathbb R^2$ is a Borel set which satisfies
\begin{enumerate}
\item[(i)] For every $e\in S^1$, there exists $x_e\in\mathbb R^2$ such that
\[\{x_e+te:t\in[0,1]\}\subset E;\]
\item[(ii)] $\mathcal L^2(E)=0$.
\end{enumerate}
Then $E$ is called a \emph{Kakeya set}.
Equivalently, we have the following \emph{Kakeya needle problem}:
\begin{quote}
In how small a plane domain can a unit segment be turned around continuously?
\end{quote}
\begin{center}\includegraphics[width=.52\linewidth]{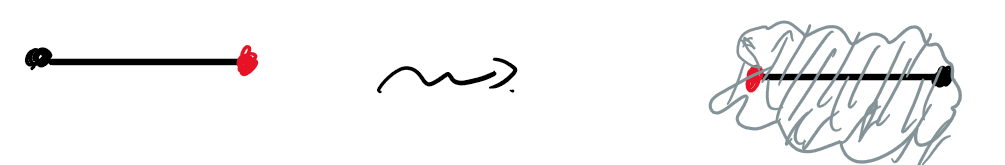}\end{center}
The existence of a Kakeya set gives the answer: arbitrarily small. See~\cite{Besicovitch1928}.

The other motivation is as follows.
\paragraph{Motivation (Besicovitch).}
If $f$ is a Riemann integrable function on $\mathbb R^2$, is it always possible to find a pair of orthogonal coordinate axes with respect to which $\int f(x,y)\,dx$ exists as a Riemann integral for all $y\in\mathbb R$, and with the resulting function of $y$ also Riemann integrable?

The answer is \textbf{NO}. Due to the existence of a Kakeya set, we can construct a counterexample. In fact, suppose $E$ is a compact Kakeya set. After translating $E$, we may assume that a chosen horizontal unit segment in $E$ lies on a line $y=y_0$ with $y_0\notin\mathbb Q$, and a chosen vertical unit segment in $E$ lies on a line $x=x_0$ with $x_0\notin\mathbb Q$. The two translation coordinates can be chosen independently to ensure these conditions. Let
\[E_0=\{(x,y)\in E:x\in\mathbb Q\text{ or }y\in\mathbb Q\},\]
and $f=1_{E_0}$. Since $E$ is compact, the discontinuity points of $f$ are contained in $\overline{E_0}\subset E$. As $\mathcal L^2(E)=0$, $f$ is Riemann integrable by the Lebesgue criterion.

For each non-horizontal, non-vertical direction, choose a unit line segment $L\subset E$ in that direction. Both coordinate functions are nonconstant on $L$, so $E_0\cap L$ is countable and dense in $L$. For the horizontal and vertical directions, use the two segments chosen above: the fixed coordinate is irrational, and the other coordinate takes rational values on a countable dense subset. In every case, both $E_0\cap L$ and $L\setminus E_0$ are dense in $L$. Thus $f|_L$ is discontinuous at every point of $L$, a set of length $1$, so it is \textbf{NOT} Riemann integrable. This gives a non-Riemann-integrable section in every direction.

\subsection*{Nikodym set}
For the original construction, see~\cite{Nikodym1927}.
A set $E\subset[0,1]^2$ is a Nikodym set if $\mathcal L^2(E)=0$ and for any $x\in[0,1]^2\setminus E$, there exists a line $l_x$ through $x$ such that
\[l_x\setminus\{x\}\subset E.\]
\begin{center}
\begin{tabular}{c|c}
$E$ & $[0,1]^2\setminus E$\\\hline
$\mathcal L^2(E)=0$ & $\mathcal L^2([0,1]^2\setminus E)=1$\\\hline
$\exists l_x\text{ such that }l_x\setminus\{x\}\subset E$ & $\forall x\in[0,1]^2\setminus E$\\\hline
$\displaystyle\bigcup_{x\in[0,1]^2\setminus E}l_x\setminus\{x\}$ & $\displaystyle\bigcup_{x\in[0,1]^2\setminus E}\{x\}$
\end{tabular}
\end{center}
\begin{center}\includegraphics[width=.25\linewidth]{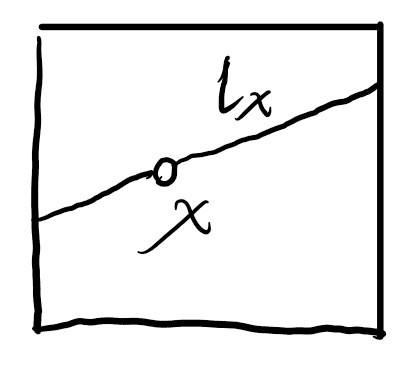}\end{center}
``Larger'' lines (compared to points) are contained in a smaller set $E$, which implies these lines have a significant overlap.

\subsection*{Variants of Kakeya sets and Nikodym sets}
\paragraph{Circular Kakeya set.}
\begin{enumerate}
\item[(1)] $\forall r\in[1,2]$, $\exists x_r\in\mathbb R^2$ such that $C(x_r,r)\subset E$;
\item[(2)] $\mathcal L^2(E)=0$.
\end{enumerate}
T. Wolff, '97: all circular Kakeya sets have dimension $2$. See~\cite{Wolff1997}.

\paragraph{Circular Nikodym set.}
\begin{enumerate}
\item[(1)] $\forall x\in[0,1]^2\setminus E$, $\exists r_x\in[1,2]$, $C(x,r_x)\subset E$;
\item[(2)] $\mathcal L^2(E)=0$.
\end{enumerate}
But circular Nikodym sets do \textbf{NOT} exist! See~\cite{Bourgain1986,Marstrand1987}.
\subsection*{Analytic number theory}
Vinogradov's mean value conjecture. See~\cite{BDG2016,Wooley2016}.

\section{Prologue}
For the Fourier-analytic preliminaries, oscillatory integrals and stationary phase used below, see~\cite{Stein1993}.
\subsection*{Uncertainty principle}
\paragraph{Lemma 1.}
If $\widehat f_j$ is supported in $\Omega_j$ and all $\Omega_j$ are disjoint, then $f_j$ are orthogonal, i.e.,
\[\Big\|\sum_j f_j\Big\|_{L^2(\mathbb R^n)}^2=\sum_j\|f_j\|_{L^2(\mathbb R^n)}^2.\]
\paragraph{Lemma 2.}
If $\widehat f_j$ is supported in $\Omega_j$ and $\operatorname{dist}(\Omega_i,\Omega_j)\ge R^{-1}$ whenever $i\ne j$, then $f_j$ are almost orthogonal, i.e.,
\[\Big\|\sum_jf_j\Big\|_{L^2(B_R)}^2\lesssim\sum_j\|f_j\|_{L^2(w_{B_R})}^2,\]
where $w_{B_R}$ is a weight associated with the ball $B_R=B(c,R)$:
\[w_{B_R}(x)=\frac1{(1+|x-c|/R)^{100n}}.\]
\paragraph{Remark.}
When we deal with the integral of $\sum f_j$ at spatial scale $R$, the minimal scale in frequency space we can distinguish is $R^{-1}$.

\subsection*{Square root cancellation}
\[\Big|\sum_jf_j(x)\Big|\sim\Big(\sum_j|f_j(x)|^2\Big)^{1/2}\quad\text{for a ``generic'' point }x\in\mathbb R^n.\]
\subsection*{Oscillatory integrals}
\paragraph{Lemma 3.}
Suppose $\phi:\mathbb R\to\mathbb R$ satisfies $\phi'(0)=1$. Then
\[\left|\int e^{i\lambda\phi(s)}a(s)\,ds\right|\lesssim_{\phi,a,N}\lambda^{-N}\qquad\text{for all }\lambda\ge1,\ N\in\mathbb N,\]
where $a$ is a smooth cutoff supported near $0$.
\begin{proof}
We will integrate by parts repeatedly:
\begin{align*}
\left|\int e^{i\lambda\phi(s)}a(s)\,ds\right|
&=\left|\int a(s)\frac1{i\lambda\phi'(s)}\,d\big(e^{i\lambda\phi(s)}\big)\right|\\
&=\left|\int\frac d{ds}\left(\frac1{i\lambda\phi'(s)}a(s)\right)e^{i\lambda\phi(s)}\,ds\right|.
\end{align*}
Continuing,
\begin{align*}
\left|\int e^{i\lambda\phi(s)}a(s)\,ds\right|
&=\left|\int\frac d{ds}\left(\frac1{i\lambda\phi'(s)}\frac d{ds}\left(\frac1{i\lambda\phi'(s)}a(s)\right)\right)e^{i\lambda\phi(s)}\,ds\right|\\
&=\cdots=\left|\int D^N(a)(s)e^{i\lambda\phi(s)}\,ds\right|\\
&\le\int_{-1}^1|D^N(a)(s)|\,ds,
\end{align*}
where $D(a)(s):=\frac d{ds}\big(\frac1{i\lambda\phi'(s)}a(s)\big)$.
And by careful calculation,
\begin{align*}
|D^N(a)(s)|
&=\lambda^{-N}\left|\frac d{ds}\left(\frac1{\phi'(s)}\frac d{ds}\left(\frac1{\phi'(s)}\cdots\frac d{ds}\left(\frac1{\phi'(s)}a(s)\right)\cdots\right)\right)\right|\\
&=\lambda^{-N}\Bigg|F\left(\frac1{\phi'},\frac d{ds}\frac1{\phi'},\ldots,\frac{d^N}{ds^N}\frac1{\phi'}\right)(s)a(s)\\
&\qquad+F\left(\frac1{\phi'},\frac d{ds}\frac1{\phi'},\ldots,\frac{d^{N-1}}{ds^{N-1}}\frac1{\phi'}\right)(s)\frac{da}{ds}(s)\\
&\qquad+\cdots+F\left(\frac1{\phi'},\frac d{ds}\frac1{\phi'}\right)(s)\frac{d^{N-1}a}{ds^{N-1}}(s)
+\left(\frac1{\phi'(s)}\right)^N\frac{d^Na}{ds^N}(s)\Bigg|.
\end{align*}
This is at most $C_{N,a,\phi}\lambda^{-N}$, as desired.
\end{proof}
\paragraph{Lemma 4.}
Suppose $a$ is as above. Then for any $\lambda\ge1$,
\[\left|\int e^{i\lambda s^2}a(s)\,ds\right|\lesssim\lambda^{-1/2}.\]
\begin{proof}
Let $\delta$ be a parameter which will be chosen later. Let $\chi$ be a smooth cutoff with $0\le\chi\le1$, equal to $1$ on $[-1/2,1/2]$ and supported in $[-1,1]$. Then
\begin{align*}
\left|\int_{\mathbb R}e^{i\lambda s^2}a(s)\,ds\right|
&\le\underbrace{\left|\int e^{i\lambda s^2}\chi(s/\delta)a(s)\,ds\right|}_{\text{the part without oscillation}}\\
&\quad+\underbrace{\left|\int e^{i\lambda s^2}(1-\chi(s/\delta))a(s)\,ds\right|}_{\text{oscillation, but the phase function is non-degenerate}}\\
&=:I+II.
\end{align*}
For $I$, $e^{i\lambda s^2}$ has no oscillation, meaning $|\lambda s^2|\le1\Rightarrow s\le\lambda^{-1/2}$, so we should choose $\delta=\lambda^{-1/2}$. And now $e^{i\lambda s^2}$ does not provide any cancellation,
so we use the triangle inequality directly:
\[I\le\int_{-\delta}^{\delta}a(s)\,ds\lesssim_a\delta.\]
For $II$, we integrate by parts repeatedly, but note that \emph{the implicit constant depends on the amplitude function}, $(1-\chi(s/\delta))a(s)$, so it actually depends on $\delta$:
\[II\le\left|\int\frac d{ds}\left((1-\chi(s/\delta))a(s)\frac1{2i\lambda s}\right)e^{i\lambda s^2}\,ds\right|.\]
\begin{center}\includegraphics[width=.34\linewidth]{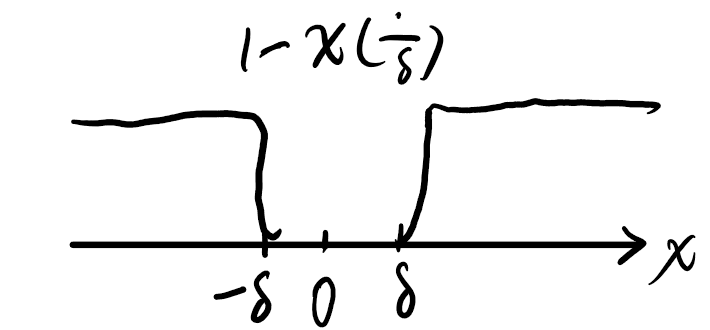}\end{center}
\[\frac d{ds}\big(1-\chi(s/\delta)\big)\approx\frac{1-0}{\delta}=\frac1\delta\quad\text{near }s=\delta,\]
\[\Longrightarrow\quad\left|\frac d{ds}\big(1-\chi(s/\delta)\big)\frac{a(s)}s\right|\lesssim\frac1{\delta^2}\quad\text{on }s\sim\delta.\]
And
\[\left|(1-\chi(s/\delta))\frac d{ds}\left(\frac{a(s)}s\right)\right|\lesssim(1-\chi(s/\delta))\frac1{|s|^2}.\]
Therefore, we have
\begin{align*}
II&\le\frac1\lambda\int_{\delta/2}^{\delta}\left|\frac d{ds}(1-\chi(s/\delta))\frac{a(s)}s\right|\,ds
+\frac1\lambda\int_\delta^\infty\left|(1-\chi(s/\delta))\frac d{ds}\left(\frac{a(s)}s\right)\right|\,ds\\
&\lesssim\frac1\lambda\cdot\delta\cdot\frac1{\delta^2}+\frac1\lambda\int_\delta^\infty\frac{ds}{s^2}
\lesssim\frac1{\lambda\delta}.
\end{align*}
In summary, recalling that we chose $\delta=\lambda^{-1/2}$, we obtain
\[\left|\int e^{i\lambda s^2}a(s)\,ds\right|\lesssim\delta+\frac1{\lambda\delta}=\lambda^{-1/2}.\]
\end{proof}
\paragraph{Remark.}
(1) In the estimate for $II$, we only integrate by parts once. So will there be any gains if we do it a few more times? The answer is \textbf{NO}. In fact,
\begin{align*}
II&\overset{\text{twice}}\le\frac1{\lambda^2}\Bigg(\int_{s\sim\delta}\left|\frac{d^2}{ds^2}(1-\chi(s/\delta))\right|\frac{(\cdot)}{s^2}\,ds\\
&\qquad+\int_{s\sim\delta}\left|\frac d{ds}(1-\chi(s/\delta))\cdot\frac d{ds}\left(\frac{(\cdot)}{s^2}\right)\right|\,ds
+\int_\delta^\infty(1-\chi(s/\delta))\left|\frac{d^2}{ds^2}\left(\frac{(\cdot)}{s^2}\right)\right|\,ds\Bigg).
\end{align*}
\begin{center}\includegraphics[width=.28\linewidth]{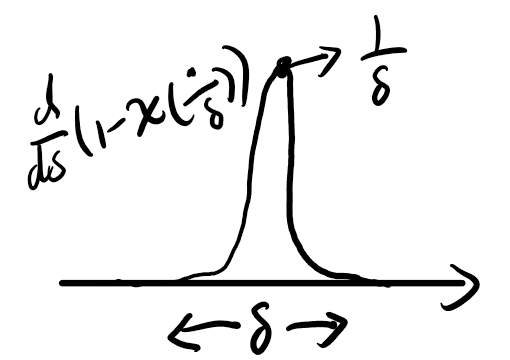}\end{center}
\[\Longrightarrow\quad\frac{d^2}{ds^2}(1-\chi(s/\delta))\sim\frac{1/\delta-0}{\delta}\sim\frac1{\delta^2}\quad\text{on }s\sim\delta.\]
So
\begin{align*}
II&\overset{\text{twice}}\le\frac1{\lambda^2}\left(\delta\cdot\frac1{\delta^2}\cdot\frac1{\delta^2}
+\delta\cdot\frac1\delta\cdot\frac1{\delta^3}+\left|\int_\delta^\infty\frac1{s^4}\,ds\right|\right)\\
&=\frac1{\lambda^2}\delta^{-3}=\lambda^{3/2-2}=\lambda^{-1/2}.
\end{align*}
In general, we have
\begin{align*}
II&\overset{N\text{ times}}\le\frac1{\lambda^N}\Bigg(\int_{s\sim\delta}\left|\frac{d^N}{ds^N}(1-\chi(s/\delta))\right|\frac{(\cdot)}{s^N}\,ds+\cdots\\
&\qquad+\int_\delta^\infty(1-\chi(s/\delta))\left|\frac{d^N}{ds^N}\left(\frac{(\cdot)}{s^N}\right)\right|\,ds\Bigg)\\
&\lesssim\frac1{\lambda^N}\left(\delta\cdot\frac1{\delta^N}\cdot\frac1{\delta^N}
+\delta\cdot\frac1{\delta^{N-1}}\frac1{\delta^{N+1}}+\cdots+\int_\delta^\infty\frac1{s^{2N}}\,ds\right).
\end{align*}
This equals
\[\frac1{\lambda^N}\delta^{1-2N}=\lambda^{-N+\frac12(2N-1)}=\lambda^{-1/2}.\]
(2) The decay of oscillatory integrals is determined by the part without oscillation.

\paragraph{Lemma 5.}
Suppose $\phi:\mathbb R\to\mathbb R$ satisfies $\phi(0)=\phi'(0)=1$ and $\phi''(0)=1$. Then for any $k\ge0$, any $\lambda\ge1$, we have
\[\left|\int e^{i\lambda\phi(s)}s^ka(s)\,ds\right|\lesssim\lambda^{-1/2}\lambda^{-k/2}.\]
\begin{proof}
\begin{align*}
\left|\int_{\mathbb R}e^{i\lambda\phi(s)}s^ka(s)\,ds\right|
&\le\left|\int_{-\lambda^{-1/2}}^{\lambda^{-1/2}}e^{i\lambda\phi(s)}s^ka(s)\,ds\right|
+\left|\int_{|s|>\lambda^{-1/2}}e^{i\lambda\phi(s)}s^ka(s)\,ds\right|\\
&\lesssim\lambda^{-1/2}(\lambda^{-1/2})^k.
\end{align*}
\end{proof}
\paragraph{Corollary 6.}
For any $\alpha\in\mathbb N$, $\phi$ as above, we have
\[\left|\frac{d^\alpha}{d\lambda^\alpha}\int e^{i\lambda\phi(s)}a(s)\,ds\right|\lesssim\lambda^{-1/2-\alpha}.\]
\begin{proof}
When $\alpha=1$, using $\phi(s)\approx s^2$,
\[\left|\frac d{d\lambda}\int e^{i\lambda\phi(s)}a(s)\,ds\right|
\approx\left|\int e^{i\lambda\phi(s)}s^2a(s)\,ds\right|\lesssim\lambda^{-1/2}\lambda^{-1}.\]
In general,
\[\left|\frac{d^\alpha}{d\lambda^\alpha}\int e^{i\lambda\phi(s)}a(s)\,ds\right|
\approx\left|\int e^{i\lambda\phi(s)}(s^2)^\alpha a(s)\,ds\right|\lesssim\lambda^{-1/2}\lambda^{-\alpha}.\]
\end{proof}
Lemma 4, Lemma 5 and Corollary 6 are called \emph{van der Corput's lemma}. Next we give an application of van der Corput's lemma. Recall
\[H^s(\mathbb R^2):=\left\{f\in\mathcal S':\big((1+|\xi|^2)^{s/2}\widehat f\big)^\vee\in L^2(\mathbb R^2)\right\},\]
and
\[\|f\|_{H^s(\mathbb R^2)}:=\left(\int_{\mathbb R^2}(1+|\xi|^2)^s|\widehat f(\xi)|^2\,d\xi\right)^{1/2}.\]
Let
\[Af(x,y)=\int f(x-t,y-t^2)a(t)\,dt.\]
We have the following estimate.
\paragraph{Proposition 7.}
\[\|Af\|_{H^{1/2}(\mathbb R^2)}\lesssim\|f\|_{L^2(\mathbb R^2)}.\]
\begin{proof}
\begin{align*}
\widehat{Af}(\xi,\eta)
&=\int_{\mathbb R^2}\int_{\mathbb R}f(x-t,y-t^2)a(t)e^{-i(x\xi+y\eta)}\,dt\,dx\,dy\\
&=\widehat f(\xi,\eta)\int_{\mathbb R}e^{-i(t\xi+t^2\eta)}a(t)\,dt.
\end{align*}
Using the Littlewood--Paley decomposition,
\[f=\sum_{k=-\infty}^\infty f_k,\quad\text{where }\widehat f_k\text{ is supported on }|\xi|\sim2^k.\]
Then we have
\begin{align*}
\|Af\|_{H^{1/2}}^2
&=\iint(1+|\xi|^2+|\eta|^2)^{1/2}|\widehat{Af}(\xi,\eta)|^2\,d\xi\,d\eta\\
&=\iint(1+|\xi|^2+|\eta|^2)^{1/2}|\widehat f(\xi,\eta)|^2
\left|\int_{\mathbb R}e^{-i(t\xi+t^2\eta)}a(t)\,dt\right|^2\,d\xi\,d\eta\\
&\sim\sum_k\iint_{\sqrt{|\xi|^2+|\eta|^2}\sim2^k}(1+2^{2k})^{1/2}|\widehat f_k(\xi,\eta)|^2
\left|\int_{\mathbb R}\cdots\,dt\right|^2\,d\xi\,d\eta\\
&\lesssim\sum_k\iint_{\sqrt{|\xi|^2+|\eta|^2}\sim2^k}|\widehat f_k(\xi,\eta)|^2\,d\xi\,d\eta
\sim\|f\|_2^2.
\end{align*}
Here the factor $\left|\int_{\mathbb R}\cdots\,dt\right|^2$ is bounded by $2^{-k}$ if $k\ge0$, and by $1$ if $k<0$.
So it is sufficient to prove that when $\sqrt{|\xi|^2+|\eta|^2}\sim2^k$, we have
\[\left|\int e^{-i(t\xi+t^2\eta)}a(t)\,dt\right|^2\lesssim
\begin{cases}2^{-k},&k\ge0,\\1,&k<0.\end{cases}\]
Note that the result when $k<0$ is trivial:
\[\left|\int e^{-i(t\xi+t^2\eta)}a(t)\,dt\right|\le\int|a(t)|\,dt\lesssim_a1.\]
And by rescaling, we will show that for $k\ge0$,
\[\left|\int e^{-i2^k(t\xi+t^2\eta)}a(t)\,dt\right|\lesssim2^{-k/2}
\quad\text{whenever }|\xi|^2+|\eta|^2\sim1.\]
\begin{center}\includegraphics[width=.3\linewidth]{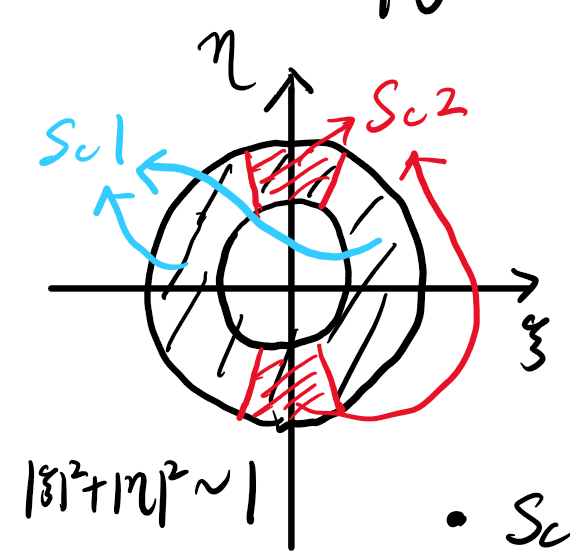}\end{center}
\[\frac\partial{\partial t}(t\xi+t^2\eta)=\xi+2t\eta=0\quad\Longleftrightarrow\quad t=-\frac\xi{2\eta}\to0,\]
\[\frac{\partial^2}{\partial t^2}(t\xi+t^2\eta)=2\eta.\]
\textbf{Sc1:} $|\xi|\ge10|\eta|$.
\[\Longrightarrow\quad\left|\frac\partial{\partial t}(t\xi+t^2\eta)\right|
\ge|\xi|-2|\eta|\ge\frac45|\xi|\ge\frac1{10}.\]
\[\Longrightarrow\quad\left|\int e^{-i2^k(t\xi+t^2\eta)}a(t)\,dt\right|\lesssim_N(2^k)^{-N}
\quad\text{for any }N\in\mathbb N.\]
\textbf{Sc2:} $|\xi|<10|\eta|$.
\[\Longrightarrow\quad\left|\frac{\partial^2}{\partial t^2}(t\xi+t^2\eta)\right|=2|\eta|\gtrsim1.\]
By van der Corput,
\[\left|\int e^{-i2^k(t\xi+t^2\eta)}a(t)\,dt\right|\lesssim(2^k)^{-1/2}.\]
\end{proof}

\subsection*{Stationary phase principle}
Consider the model example
\[I_{\lambda,b}=\int_{\mathbb R}e^{i\lambda(s^2-bs)}a(s)\,ds,\qquad b\in\mathbb R.\]
So far, all the information we know about $I_{\lambda,b}$ is its decay, $\lambda^{-1/2}$, but there is no information about oscillations. However, sometimes we need to know how it oscillates to help us gain additional cancellation.
\paragraph{Lemma 8.}
\begin{align*}
I_{\lambda,b}&\approx\underbrace{e^{i\lambda[(b/2)^2-b\cdot b/2]}}_{\text{oscillation}}a(b/2)\underbrace{\lambda^{-1/2}}_{\text{decay}}\\
&=e^{-i\lambda b^2/4}a(b/2)\lambda^{-1/2}.
\end{align*}
Here, ``$\approx$'' means that there exists a constant $c_0$ such that
\[\left|I_{\lambda,b}-c_0e^{-i\lambda b^2/4}a(b/2)\lambda^{-1/2}\right|\lesssim\lambda^{-1}.\]
\paragraph{Remark.}
The stationary phase principle means that the behavior of $I_{\lambda,b}$ is determined by the critical point of the phase function:
\[\frac\partial{\partial s}(s^2-bs)=0\quad\Longrightarrow\quad\text{critical point }s=b/2.\]
\begin{proof}[Proof of Lemma 8]
\[\int_{\mathbb R}e^{i\lambda(s^2-bs)}a(s)\,ds
=e^{-i\lambda b^2/4}\int_{\mathbb R}e^{i\lambda(s-b/2)^2}a(s)\,ds.\]
(Let $s_0$ be a critical point of $\phi$. Here we are actually doing
\[\phi(s)=\phi(s_0)+\underbrace{\phi'(s_0)(s-s_0)}_{=0}+\tfrac12\phi''(s_0)(s-s_0)^2.\])
\begin{align*}
\int_{\mathbb R}e^{i\lambda(s^2-bs)}a(s)\,ds
&=e^{-i\lambda b^2/4}a(b/2)\int_{\mathbb R}e^{i\lambda(s-b/2)^2}\,ds\\
&\quad+e^{-i\lambda b^2/4}\int_{\mathbb R}e^{i\lambda(s-b/2)^2}[a(s)-a(b/2)]\,ds\\
&=:I+II.
\end{align*}
For $I$, let $s=\lambda^{-1/2}(t+b/2)$. Then
\[\lambda(s-b/2)^2=\bigl(t+(1-\sqrt\lambda)b/2\bigr)^2.\]
Translating once more by $r=t+(1-\sqrt\lambda)b/2$, we have
\[I=e^{-i\lambda b^2/4}a(b/2)\lambda^{-1/2}\underbrace{\int_{\mathbb R}e^{ir^2}\,dr}_{\text{constant}}.\]
For $II$, by van der Corput's lemma,
\[|II|\lesssim\left|\int_{\mathbb R}e^{i\lambda(s-b/2)^2}(s-b/2)\,ds\right|
\lesssim\lambda^{-1/2}\cdot\lambda^{-1/2}=\lambda^{-1}.\]
\end{proof}
In general, we have the following lemma.
\paragraph{Lemma 9.}
Suppose $\phi:\mathbb R\to\mathbb R$ satisfies $\phi'(0)=0$ and $\phi''(0)=1$. Then
\[\int_{\mathbb R}e^{i\lambda\phi(s)}a(s)\,ds=c_0\lambda^{-1/2}e^{i\lambda\phi(0)}a(0)+O(\lambda^{-1}).\]

\paragraph{Theorem 10.}
Let $\phi(x,s):\mathbb R^n\times\mathbb R\to\mathbb R$ satisfy
\[\partial_s\phi(0,0)=0,\qquad\partial_{ss}^2\phi(0,0)\ne0.\]
Therefore, for $x$ small enough, we can find $s(x)$ such that
\[\partial_s\phi(x,s(x))=0.\]
And then
\[\int_{\mathbb R}e^{i\lambda\phi(x,s)}a(x,s)\,ds=e^{i\lambda\phi(x,s(x))}b(\lambda,x),\]
where $b\in C^\infty$ and satisfies
\[|\partial_\lambda^\alpha\partial_x^\beta b(\lambda,x)|\lesssim_{\alpha,\beta,\phi,a}(1+|\lambda|)^{-1/2-\alpha}
\quad\text{for any }\lambda\in\mathbb R.\]
\begin{proof}
Note that
\[b(\lambda,x)=\int_{\mathbb R}e^{i\lambda[\phi(x,s)-\phi(x,s(x))]}a(x,s)\,ds.\]
\begin{align*}
\Longrightarrow\quad\left|\frac\partial{\partial\lambda}b(\lambda,x)\right|
&=\left|\int_{\mathbb R}e^{i\lambda[\phi(x,s)-\phi(x,s(x))]}[\phi(x,s)-\phi(x,s(x))]a(x,s)\,ds\right|.
\end{align*}
Here the bracket is expanded in the notes as
\[\underbrace{\partial_s\phi(x,s(x))(s-s(x))}_{=0}
+\tfrac12\underbrace{\partial_{ss}^2\phi(x,s(x))}_{\text{constant}}(s-s(x))^2.\]
Thus, by van der Corput,
\[\left|\frac\partial{\partial\lambda}b(\lambda,x)\right|\lesssim\lambda^{-1/2}\cdot\lambda^{-1}
\qquad(\lambda\to+\infty).\]
\[\left|\frac\partial{\partial x}b(\lambda,x)\right|
=\left|\int_{\mathbb R}e^{i\lambda[\phi(x,s)-\phi(x,s(x))]}A(x,s)\,ds\right|
\lesssim\lambda^{-1/2}\qquad(\lambda\to+\infty).\]
\end{proof}
\section{Decoupling for parabola}
Given $\delta=2^{-N}$ and an interval $I\subset[-1,1]$, define
\[
\Omega_I(\delta)=\{(\xi_1,\xi_1^2+\delta'):\xi_1\in I\ \&\ |\delta'|\leq\delta^2\}.
\]
\begin{center}\includegraphics[width=.29\linewidth]{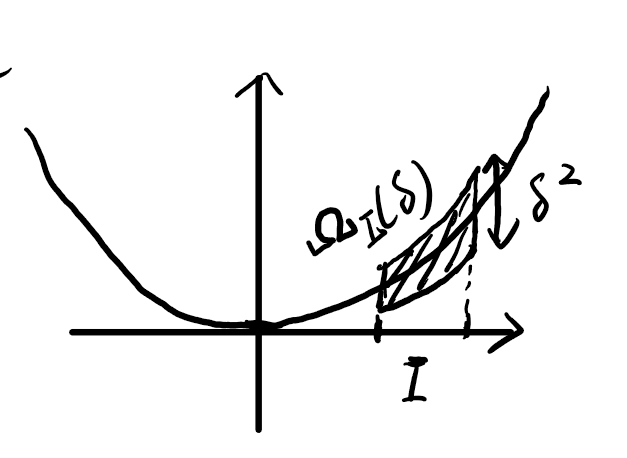}\end{center}
Let $p\in(1,\infty)$, $F:\mathbb R^2\to\mathbb C$ be an $L^p$ function with Fourier transform supported in $\Omega_{[-1,1]}(\delta)$. Define
\[
F_I=(\widehat F\cdot 1_{\Omega_I(\delta)})^\vee.
\]
By the boundedness of the Hilbert transform, it is easy to verify $F_I$ is well-defined. And we have the (uniform) decomposition
\[
F=\sum_{\substack{I\subset[-1,1]\\\ell(I)=\delta}}F_I.
\]
Decoupling inequality aims to compare
\[
\left\|\sum_I F_I\right\|_{L^p(\mathbb R^2)}
\quad\text{v.s.}\quad
\left(\sum_I\|F_I\|_{L^p(\mathbb R^2)}^2\right)^{1/2}.
\]
Let $D_p(\delta)$ be the smallest constant for which the following inequality
\[
\|F\|_{L^p(\mathbb R^2)}\leq D_p(\delta)
\left(\sum_I\|F_I\|_{L^p(\mathbb R^2)}^2\right)^{1/2}
\]
holds for any $F$ with Fourier transform supported in $\Omega_{[-1,1]}(\delta)$.
Trivially, by the triangle inequality and Cauchy--Schwarz inequality,
\begin{align*}
\|F\|_{L^p(\mathbb R^2)}
&\leq\sum_I\|F_I\|_{L^p(\mathbb R^2)}\\
&\leq(\# I)^{1/2}\left(\sum_I\|F_I\|_{L^p(\mathbb R^2)}^2\right)^{1/2}\\
&\lesssim\delta^{-1/2}\left(\sum_I\|F_I\|_{L^p}^2\right)^{1/2}.
\end{align*}
Conversely, fix $I=[0,\delta]$, choose $F$ satisfying $\operatorname{supp}(\widehat F)\subset\Omega_I(\delta)$, then
\[
F_J=\begin{cases}F_I&\text{if }J=I,\\0&\text{if }J\neq I,\end{cases}
\qquad F_I=F.
\]
By the definition of decoupling we get
\[
D_p(\delta)\geq\frac{\|F\|_p}{(\sum_I\|F_I\|_p^2)^{1/2}}
=\frac{\|F\|_p}{\|F\|_p}=1.
\]
So far, we obtain the trivial bound
\[
1\leq D_p(\delta)\lesssim\delta^{-1/2}.
\]
\paragraph{Theorem 1 (Bourgain--Demeter, 2015).} \cite{BD2015}
\[
D_p(\delta)\lesssim_\varepsilon
\begin{cases}
\delta^{-\varepsilon},&2\leq p\leq6,\\
\delta^{-1/2+3/p-\varepsilon},&p\geq6.
\end{cases}
\]
\paragraph{Lemma 2.} $D_2(\delta)=1$, $D_\infty(\delta)\leq\delta^{-1/2}$.
\begin{proof}
For $p=2$, by Plancherel theorem,
\[
\|F\|_2^2=\left\|\sum_I\widehat F_I\right\|_2^2
=\sum_I\|\widehat F_I\|_2^2=\sum_I\|F_I\|_2^2,
\]
since $\{I\}$ are disjoint. ($L^2$ has the best orthogonality, so we can separate completely every frequency $F_I$.)
For $p=\infty$, this is the trivial bound. (For $L^\infty$-norm, we don't expect any cancellation, so we directly separate them using triangle inequality.)
\end{proof}
\paragraph{Example 1.} When $p>6$, we have
\[
D_p(\delta)\gtrsim\delta^{-1/2+3/p},
\]
which implies Theorem 1 is sharp up to an $\varepsilon$-loss.
In fact, let $\chi:\mathbb R\to\mathbb R$ be a characteristic function on $[-1,1]$ and $\chi_0$ be the smooth version of $\chi$. Define $F:\mathbb R^2\to\mathbb C$ by
\[
\widehat F(\xi)=\chi_0(\xi_1)\frac1{\delta^2}
\chi_0\left(\frac{\xi_2-\xi_1^2}{\delta^2}\right)
\qquad\left(\approx\frac1{\delta^2}1_{\Omega_{[-1,1]}(\delta)}\right).
\]
\paragraph{Lemma 3.} $\|F\|_{L^p(\mathbb R^2)}\gtrsim1$ for any $p\in(1,\infty)$. Moreover, $\|F\|_{L^p(\mathbb R^2)}\lesssim1$ for $p\in(4,\infty)$.
\begin{proof}
The maximal rectangle contained in $\Omega_{[-1,1]}(\delta)$ has dimension $\sim\delta\times\delta^2$. So in order to use local constant property, we divide $\Omega_{[-1,1]}(\delta)$ into $\bigsqcup_{\ell(I)=\delta}\Omega_I(\delta)$.
\begin{center}\includegraphics[width=.42\linewidth]{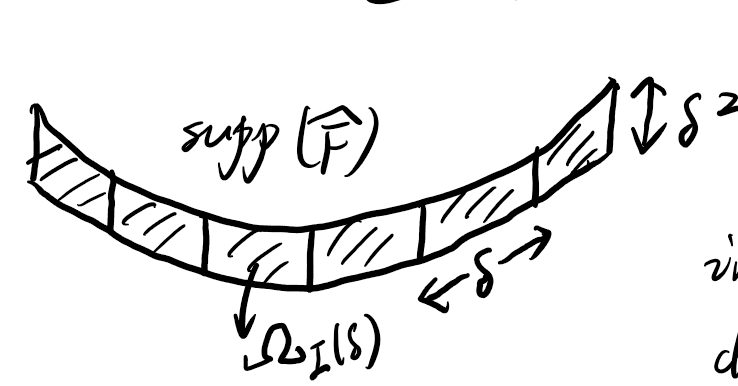}\end{center}
$F=\sum_I F_I$; for each $I$, $\operatorname{supp}(\widehat F_I)\subset\Omega_I(\delta)$.
\begin{center}\includegraphics[width=.75\linewidth]{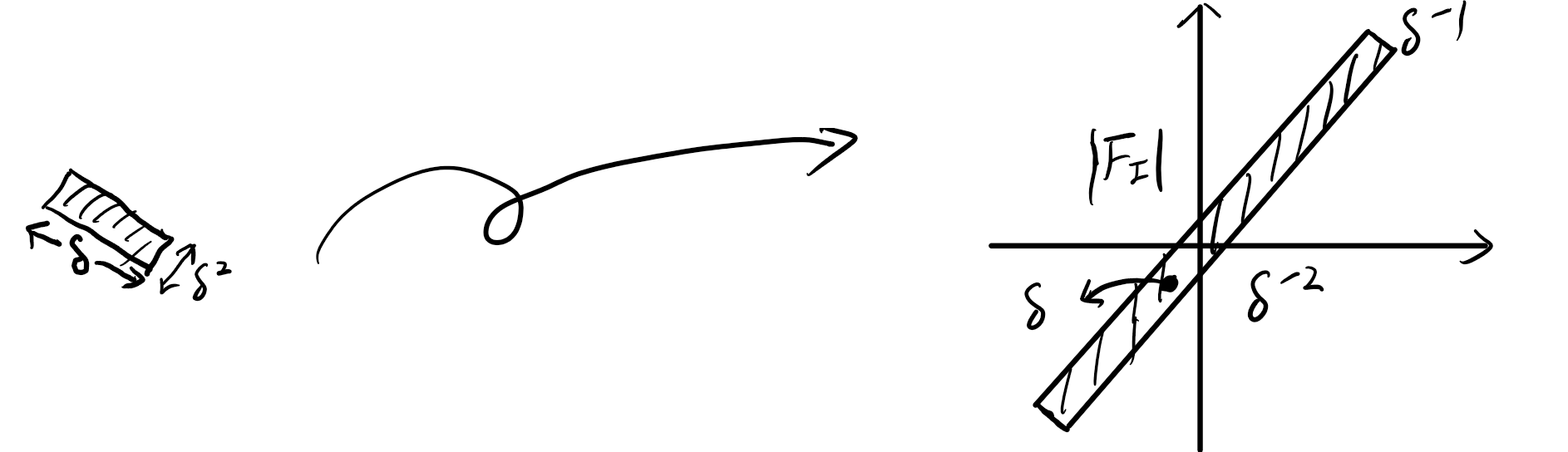}\end{center}
\[
\widehat F_I\sim\frac1{\delta^2}1_{\Omega_I(\delta)}
\]
implies that $|F_I|$ concentrates on a tube dual to $\Omega_I(\delta)$ and centered at the origin.
\[
\|\widehat F_I\|_{L^1}\sim\frac{|\Omega_I(\delta)|}{\delta^2}=\delta
\quad\Longrightarrow\quad
|F_I(x)|\sim\int|\widehat F_I|\sim\delta\quad\text{on this tube}.
\]
Therefore, back to the physical space, we know $F$ concentrates on a brush at the origin. There are interactions between different $F_I$, so we can expect there are some cancellations at some points.
\begin{center}\includegraphics[width=.30\linewidth]{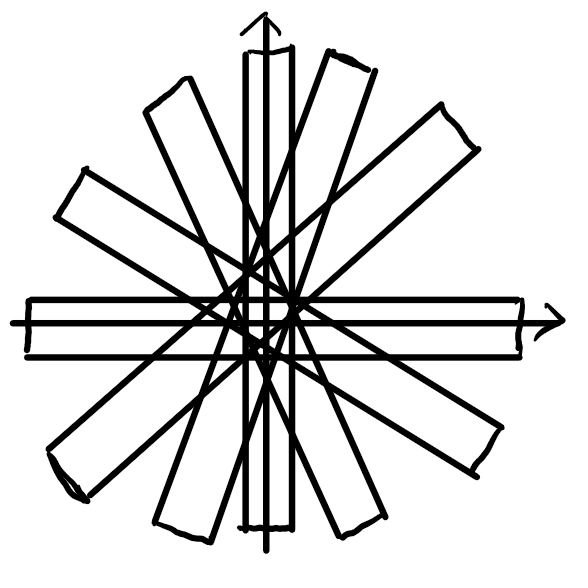}\end{center}
\begin{itemize}
\item On $B(0,\delta^{-2})\setminus B(0,\delta^{-2}/2)$: in this area, all tubes are completely separated. Hence,
\[
|F(x)|=\sum_I|F_I(x)|\sim\delta.
\]
\begin{center}\includegraphics[width=.32\linewidth]{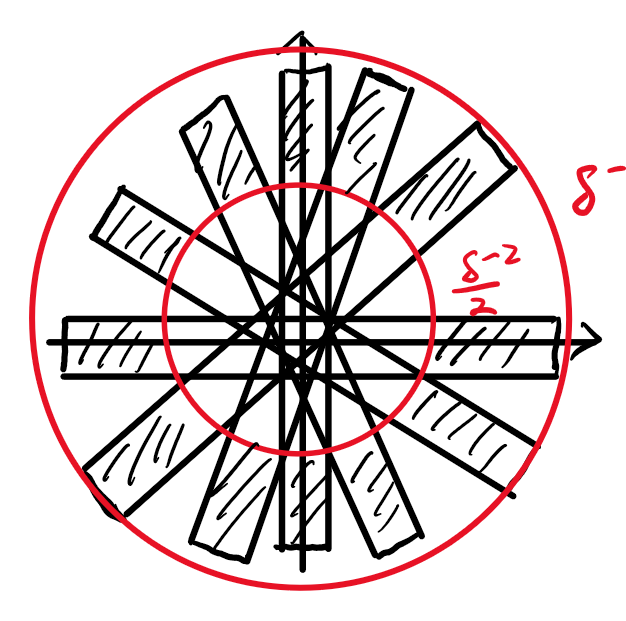}\end{center}
\item On $B(0,\delta^{-1})$: since $|F_I|$ is concentrated on scale $\delta^{-1}$, by the uncertainty principle, on frequency space the smallest scale on which we can distinguish $\widehat F_I$ is $\delta$, so
\[
\widehat F_I\approx1_{\Omega_I(\delta)}\approx1_{B(c_I,\delta)}\approx\delta_{c_I},
\]
the Dirac function at $c_I$.
\begin{center}\includegraphics[width=.75\linewidth]{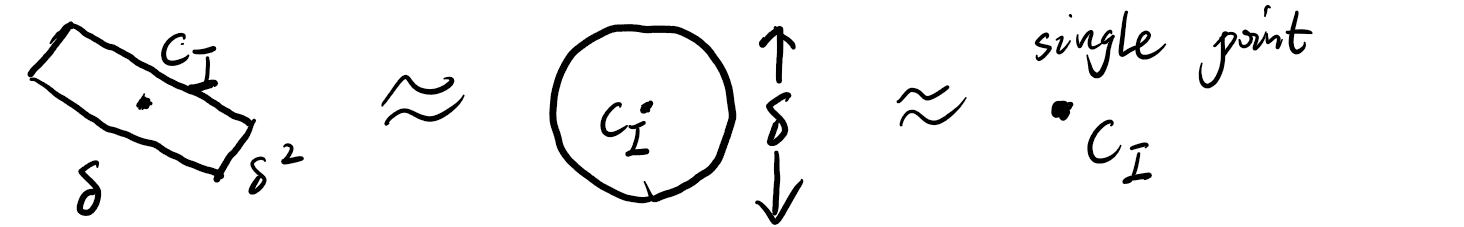}\end{center}
In fact, when $|x|\leq\delta^{-1}$,
\begin{align*}
F_I(x)&=\int\widehat F_I(\xi)e^{ix\cdot\xi}\,d\xi
=e^{ix\cdot c_I}\int\widehat F_I(\xi)e^{ix\cdot(\xi-c_I)}\,d\xi\\
&\sim\|\widehat F_I\|_1e^{ix\cdot c_I}\sim\delta e^{ix\cdot c_I},
\end{align*}
since $|\xi_1-c_I^1|\leq\delta$, $|\xi_2-c_I^2|\leq\delta^2\leq\delta$, implying $|x\cdot(\xi-c_I)|\leq1$ and
$e^{ix\cdot(\xi-c_I)}\approx e^{i0}=1$.
Therefore,
\[
F(x)=\sum_I F_I(x)\approx\sum_I(\pm1)\cdot\delta,
\]
where the signs occur randomly. Hence we obtain the square root cancellation: for ``generic'' $x\in B(0,\delta^{-1})$,
\[
|F(x)|\approx\left(\sum_I|F_I(x)|^2\right)^{1/2}
\quad\Longrightarrow\quad
|F(x)|\sim\delta\cdot(\#I)^{1/2}\sim\delta^{1/2}.
\]
\item On $B(0,1)$: constructive interference / (trivial) major arc.
\[
|F(x)|=\left|\sum_I F_I(x)\right|\approx\sum_I F_I(x)\sim1.
\]
On the one hand, for each $I$,
\[
F_I(x)=\int\widehat F_I(\xi)e^{ix\cdot\xi}\,d\xi
\approx\|\widehat F_I\|_1e^{i0}>0,
\]
since $\xi\in\Omega_I(\delta)\subset B(0,1)$ implies $|x\cdot\xi|\leq1$.
So $|\sum_I F_I(x)|\approx\sum_I F_I(x)\sim1$.
On the other hand, since $\operatorname{supp}(\widehat F)\subset\Omega_{[-1,1]}(\delta)\subset B(0,1)$, we know $|F|$ is essentially constant on $B(0,1)$. Note that
\[
F(0)=\int\widehat F(\xi)\,d\xi
=\frac1{\delta^2}|\Omega_{[-1,1]}(\delta)|\sim1;
\]
we have $|F(x)|\sim1$.
\end{itemize}
\paragraph{Remark.} In fact, we can prove that $|F(x)|\simeq|x|^{-1/2}$.
In summary,
\begin{align*}
\int|F|^p
&\sim\int_{B(0,1)}|F|^p+\int_{B(0,\delta^{-1})}|F|^p
+\int_{B(0,\delta^{-2})\setminus B(0,\delta^{-2}/2)}|F|^p\\
&\sim1+\delta^{p/2}\cdot\delta^{-2}+\delta^p\cdot\delta^{-4}
\sim1+\delta^{-(4-p)}.
\end{align*}
(Since $\displaystyle\frac1{\delta^{-(2-p/2)}}
=\frac{1-\theta}{1}+\frac{\theta}{\delta^{-(4-p)}}$ for some $\theta\in(0,1)$.)
For any $p\in(1,\infty)$, $\|F\|_p\gtrsim1$. If $p>4$, then $\|F\|_p\lesssim1$.
\end{proof}
Come back to the example. When $p>6$,
\[
\|F\|_p\sim1,\qquad
\left(\sum_I\|F_I\|_p^2\right)^{1/2}
\sim\delta\cdot\delta^{-3/p}\cdot\delta^{-1/2},
\]
which implies $D_p(\delta)\gtrsim\delta^{-1/2+3/p}$.
\subsection*{Other versions of decoupling}
For the local and extension-operator formulations, see~\cite{BD2015,Demeter2020}.
Local decoupling: for $R\geq\delta^{-2}$,
\[
\|F\|_{L^p(B_R)}\leq D_p(\delta)
\left(\sum_I\|F_I\|_{L^p(w_{B_R})}^2\right)^{1/2}.
\]
Decoupling of Fourier extension operator:
\[
\|Ef\|_{L^p(B_R)}\leq D_p(\delta)
\left(\sum_I\|E_If\|_{L^p(w_{B_R})}^2\right)^{1/2}.
\]
\paragraph{Remark.} In general, we don't consider the global norm $\|Ef\|_{L^p(\mathbb R^2)}$. Since if we choose $f=1_{[0,1]}$, then $\|Ef\|_{L^p(\mathbb R^2)}=\infty$ for $p\leq4$.
\paragraph{Example.} Discrete distribution $\longleftrightarrow$ exponential sum.
Let $\Lambda\subset\mathbb P^1$ be the maximal $\delta$-separated set, $R=\delta^{-2}$,
\[
f=\sum_{\xi\in\Lambda}\delta_\xi,\qquad
Ef=\sum_{\xi\in\Lambda}E\delta_\xi
=\sum_{\xi\in\Lambda}e^{ix\cdot\xi}
\quad\left(=\sum_{t\in\delta\mathbb Z\cap[0,1]}e^{i(x_1t+x_2t^2)}\right).
\]
When $|x|\ll1$, we expect the constructive interference:
\[
Ef(x)\sim\sum_{\xi\in\Lambda}1=|\Lambda|\quad\text{for }|x|\ll1,
\]
implying $\|Ef\|_{L^p(B_R)}\gtrsim|\Lambda|\sim\delta^{-1}$.
\begin{center}\includegraphics[width=.35\linewidth]{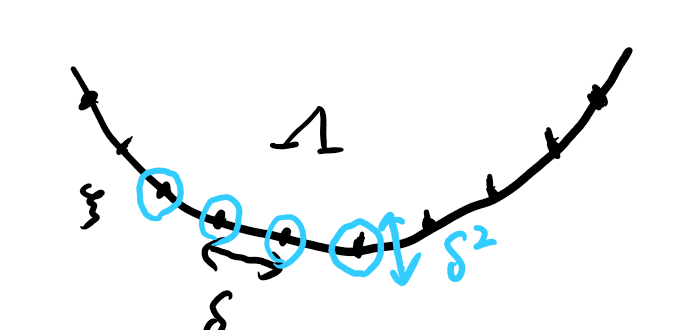}\end{center}
\begin{align*}
\left(\sum_{\ell(I)=\delta}\|E_If\|_{L^p(B_R)}^2\right)^{1/2}
&=\left(\sum_{\xi\in\Lambda}\|E\delta_\xi\|_{L^p(T_\xi)}^2\right)^{1/2}\\
&=(|\Lambda|\cdot|T_\xi|^{2/p})^{1/2}
\sim\delta^{-1/2}\delta^{-3/p}.
\end{align*}

Therefore $D(\delta)\gtrsim\delta^{3/p-1/2}$.

\section{Applications of decoupling in analytic number theory}
\subsection*{Quadratic Vinogradov's system}
Given $s\geq1$ an integer, $X$ a large number. Consider the following equations:
\begin{equation*}
\begin{cases}
x_1+\cdots+x_s=x_{s+1}+\cdots+x_{2s},\\
x_1^2+\cdots+x_s^2=x_{s+1}^2+\cdots+x_{2s}^2,
\end{cases}
\qquad 0\leq x_i\leq X.\tag{$*$}
\end{equation*}
$(x_1,\ldots,x_{2s})\in\mathbb N^{2s}$ satisfying $(*)$ is called a solution of $(*)$.
\paragraph{Question.} $\#\{\text{solutions of }(*)\}=?$ as $X\to\infty$?
When $s=1$,
\[
\begin{cases}x_1=x_2,\\x_1^2=x_2^2,\end{cases}
\qquad\#\{\text{solutions}\}=X+1.
\]
When $s=2$,
\[
\begin{cases}x_1+x_2=x_3+x_4,\\x_1^2+x_2^2=x_3^2+x_4^2,\end{cases}
\Longleftrightarrow
\frac12[(x_1,x_1^2)+(x_2,x_2^2)]
=\frac12[(x_3,x_3^2)+(x_4,x_4^2)].
\]
\begin{center}\includegraphics[width=.55\linewidth]{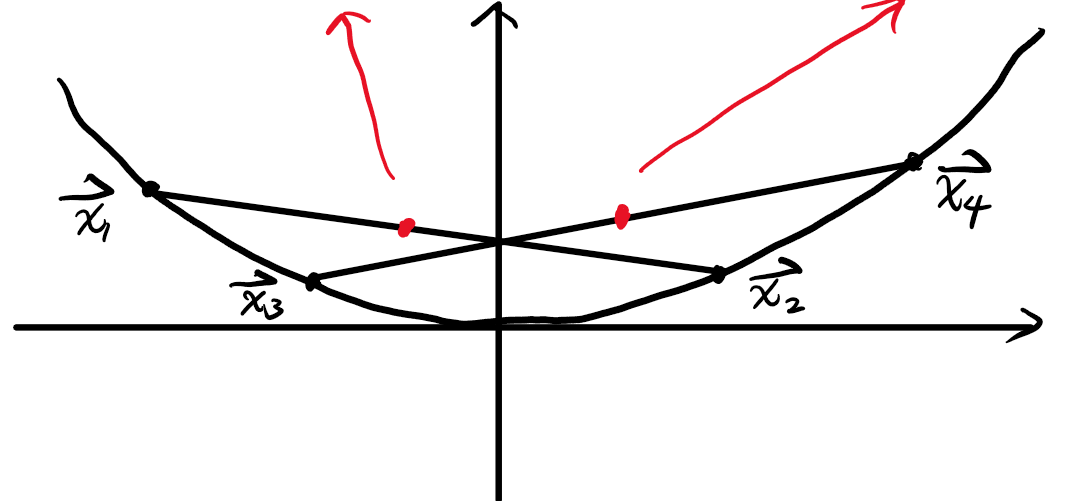}\end{center}
Hence, $(x_1,x_2,x_3,x_4)$ is a solution iff $\vec x_1,\vec x_2,\vec x_3,\vec x_4$ form a parallelogram, iff $x_1=x_3$, $x_2=x_4$ or $x_1=x_4$, $x_2=x_3$. (This solution is called a diagonal solution.) Therefore $\#\{\text{solutions}\}\sim X^2$.
For the quadratic mean-value estimate and its derivation from decoupling, see~\cite{BD2015,Demeter2020}.
\paragraph{Theorem 1.} For any $s\geq1$, denote by $N_s(X)$ the number of solutions of $(*)$, then
\[
N_s(X)\lesssim_{s,\varepsilon}X^{s+\varepsilon}+X^{2s-3+\varepsilon}.
\]
\paragraph{Remarks.}
\begin{enumerate}
\item The term $X^s$ comes from diagonal solutions: $x_1=x_{s+1},\ldots,x_s=x_{2s}$.
The term $X^{2s-3}$ comes from
\[
\#\{(x_1,\ldots,x_{2s})\in\mathbb N^{2s}:0\leq x_i\leq X\}\sim X^{2s},
\]
but $x_1+\cdots+x_s=x_{s+1}+\cdots+x_{2s}$ reduces the number of degrees of freedom by $1$,
$x_1^2+\cdots+x_s^2=x_{s+1}^2+\cdots+x_{2s}^2$ reduces the number of degrees of freedom by $2$.
($x_1^2+x_2^2=x_{s+1}^2+\cdots+x_{2s}^2-x_3^2-\cdots-x_s^2$, RHS takes values $0,\ldots,X^2$. And we know $a^2+b^2=n$, $n=0,\ldots,X^2$, on average, there is only one solution.)
So we obtain $\sim X^{2s-1-2}=X^{2s-3}$ solutions.
\item When $s\geq2s-3$, i.e. $s\leq3$, solutions of $(*)$ are almost always diagonal solutions! (Trivial case.) When $s\gg3$, non-trivial solutions begin to dominate the main term.
\item Using the algebraic method / number theory method, we can prove
\[
N_3(X)=X^3\log X+\text{lower order}.
\]
\end{enumerate}
\paragraph{Theorem 2.} For $p\geq2$, $X\geq0$, we have
\[
\left\|\sum_{n=0}^X e^{2\pi i(nx+n^2y)}\right\|_{L^p([0,1]^2)}
\lesssim_{p,\varepsilon}X^{1/2+\varepsilon}+X^{1-3/p+\varepsilon}.
\]
First, we point out Theorem 2 is equivalent to Theorem 1.
Let $p=2s$, consider
\[
u(x,y)=\sum_{n=0}^X e^{2\pi i(nx+n^2y)},
\]
then we have
\begin{align*}
\|u\|_{L^{2s}([0,1]^2)}^{2s}
&=\iint_{[0,1]^2}\left(\sum_n e^{2\pi i(nx+n^2y)}\right)^s
\overline{\left(\sum_n e^{2\pi i(nx+n^2y)}\right)^s}\,dx\,dy\\
&=\sum_{0\leq n_1,\ldots,n_{2s}\leq X}
\iint_{[0,1]^2}
\exp\Bigl(2\pi i\bigl[(n_1+\cdots+n_s-n_{s+1}-\cdots-n_{2s})x\\
&\hspace{40mm}+(n_1^2+\cdots+n_s^2-n_{s+1}^2-\cdots-n_{2s}^2)y\bigr]\Bigr)\,dx\,dy.
\end{align*}

(Note that $\displaystyle\int_0^1 e^{2\pi int}\,dt=\begin{cases}1&\text{if }n=0,\\0&\text{if }n\neq0.\end{cases}$)
Thus
\[
\|u\|_{L^{2s}([0,1]^2)}^{2s}
=\sum_{\substack{(n_1,\ldots,n_{2s})\\\text{is a solution of }(*)}}1=N_s(X).
\]
Next, Theorem 2 is sharp by the following lemma.
\paragraph{Lemma 3.} Let $u$ be as above. Then
\[
\|u\|_{L^p([0,1]^2)}\gtrsim
\underbrace{X^{1/2}}_{\text{square root cancellation}}
+\underbrace{X^{1-3/p}}_{\text{constructive interference}}.
\]
\begin{proof}
On the one hand, by Parseval's identity,
\[
\|u\|_{L^p([0,1]^2)}\geq\|u\|_{L^2([0,1]^2)}
=\left(\sum_{n=0}^X\|e^{2\pi i(nx+n^2y)}\|_{L^2([0,1]^2)}^2\right)^{1/2}
\sim X^{1/2}.
\]
On the other hand, when $|x|\leq X^{-1}$, $|y|\leq X^{-2}$, $|nx+n^2y|\leq1$, and then $e^{2\pi i(nx+n^2y)}\approx1$ for all $0\leq n\leq X$. Therefore,
\[
\|u\|_{L^p([0,1]^2)}
\geq\left\|\sum_{n=0}^X e^{2\pi i(nx+n^2y)}\right\|_{L^p([0,X^{-1}]\times[0,X^{-2}])}
\sim X\cdot(X^{-3})^{1/p}.
\]
\end{proof}
\begin{proof}[Proof of Theorem 2]
Note that
\[
\left(\sum_{n=0}^X e^{2\pi i(nx+n^2y)}\right)^{\wedge}
=\sum_{n=0}^X\delta_{(n,n^2)}.
\]
Take
\[
\widehat F=\sum_{n=0}^X\delta_{(n/X,(n/X)^2)},\qquad
F=\sum_{n=0}^X e^{2\pi i(x\cdot n/X+y\cdot(n/X)^2)}.
\]
\begin{center}\includegraphics[width=.40\linewidth]{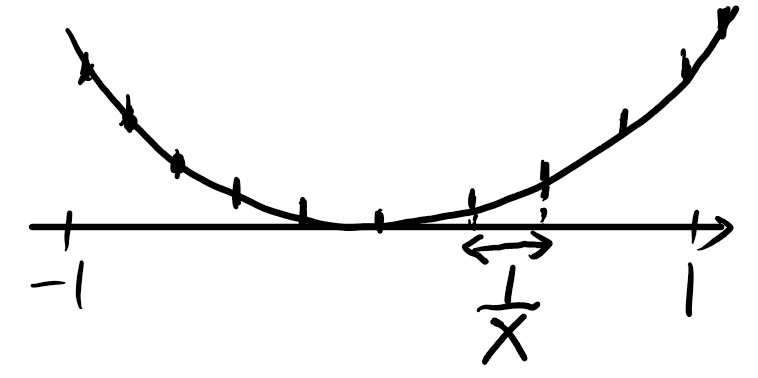}\end{center}
Apply local decoupling inequality for $F$, choose $R=X^2$, then we have
\[
\|F\|_{L^p(B_{X^2})}\lesssim_\varepsilon
\left(\underbrace{X^\varepsilon}_{2\leq p\leq6}
+\underbrace{X^{1/2-3/p+\varepsilon}}_{6\leq p\leq\infty}\right)
\left(\sum_{\ell(I)=1/X}\|F_I\|_{L^p(B_{X^2})}^2\right)^{1/2}.
\]
\[
\left(F_I=\left(\sum_n\delta_{(n/X,(n/X)^2)}1_{\Omega_I(X^{-1})}\right)^\vee
=e^{2\pi i(x\cdot n/X+y\cdot(n/X)^2)}\quad\text{where }n/X\in I\right).
\]
Thus the right hand side is
\[
\sim(X^\varepsilon+X^{1/2-3/p+\varepsilon})X^{4/p}\cdot X^{1/2}.
\]
By scaling, $\widetilde x=x/X$, $\widetilde y=y/X^2$,
\begin{align*}
\|F\|_{L^p(B_{X^2})}^p
&=\int_{B_{X^2}}\left|\sum_{n=0}^X e^{2\pi i(x\cdot n/X+y(n/X)^2)}\right|^p\,dx\,dy\\
&=X^3\int_{\substack{|\widetilde x|\leq X\\|\widetilde y|\leq1}}
\left|\sum_{n=0}^X e^{2\pi i(n\widetilde x+n^2\widetilde y)}\right|^p\,d\widetilde x\,d\widetilde y\\
&=X^4\int_{\substack{|\widetilde x|\leq1\\|\widetilde y|\leq1}}
\left|\sum_{n=0}^X e^{2\pi i(n\widetilde x+n^2\widetilde y)}\right|^p\,d\widetilde x\,d\widetilde y\\
&=X^4\|u\|_{L^p([0,1]^2)}^p,
\end{align*}
where we used periodicity. Finally,
\[
\|u\|_{L^p([0,1]^2)}\lesssim_{\varepsilon,p}X^{1/2+\varepsilon}+X^{1-3/p+\varepsilon}.
\]
\end{proof}
\paragraph{Remark.} Compared to Waring's problem
\[
\begin{cases}
x_1+\cdots+x_s=x_{s+1}+\cdots+x_{2s},\\
x_1^3+\cdots+x_s^3=x_{s+1}^3+\cdots+x_{2s}^3,
\end{cases}
\]
Vinogradov's system has better symmetry / translation invariance: if $(x_1,\ldots,x_{2s})$ is a solution, then $(x_1+a,\ldots,x_{2s}+a)$ is also a solution.
\subsection*{Discrete Strichartz inequality}
For periodic Strichartz estimates, see~\cite{Bourgain1993,BD2015}.
Consider the Schr\"odinger equation on the torus:
\[
\begin{cases}i\partial_tu-\Delta u=0&\text{in }\mathbb T^2,\\
u(\cdot,0)=u_0&\text{on }\mathbb T.
\end{cases}
\]
There is a solution
\[
u(x,t)=\sum_{n\in\mathbb Z}\widehat u_0(n)e^{i(n^2t+nx)},
\qquad(x,t)\in\mathbb T^2.
\]
\paragraph{Theorem 4.} For any $X\geq0$, it holds that
\[
\left\|\sum_{|n|\leq X}\widehat u_0(n)e^{i(n^2t+nx)}\right\|_{L^6(\mathbb T^2)}
\lesssim_\varepsilon X^\varepsilon\left(\sum_{|n|\leq X}|\widehat u_0(n)|^2\right)^{1/2}.
\]
``Reverse H\"older's inequality'':
\[
\left(\sum_{|n|\leq X}|\widehat u_0(n)|^2\right)^{1/2}
=\left\|\sum_{|n|\leq X}\widehat u_0(n)e^{i(n^2t+nx)}\right\|_{L^2(\mathbb T^2)}.
\]
\subsection*{Vinogradov's mean value theorem}
The mean-value bound below follows from~\cite{BD2015,Wooley2016,BDG2016}.
Consider the following system of equations:
\[
\begin{cases}
x_1+\cdots+x_s=x_{s+1}+\cdots+x_{2s},\\
x_1^2+\cdots+x_s^2=x_{s+1}^2+\cdots+x_{2s}^2,\\
\hspace{25mm}\vdots\\
x_1^k+\cdots+x_s^k=x_{s+1}^k+\cdots+x_{2s}^k.
\end{cases}
\]
\[
\#\{\text{solutions}\}\lesssim_\varepsilon X^{s+\varepsilon}
+X^{2s-k(k+1)/2+\varepsilon}.
\]

\section{Construction of Kakeya sets}
For the Besicovitch--Perron construction and its refinements, see~\cite{Besicovitch1928,Perron1928,Schoenberg1962}.
Recall $E\subset\mathbb R^2$ is a Kakeya set if for any $e\in S^1$, there exists $x_e\in\mathbb R^2$ such that the line segment $\{x_e+te:t\in[0,1]\}\subset E$ and $\mathcal L^2(E)=0$.
\subsection*{Construct triangles with heavy overlapping}
\begin{center}\includegraphics[width=.54\linewidth]{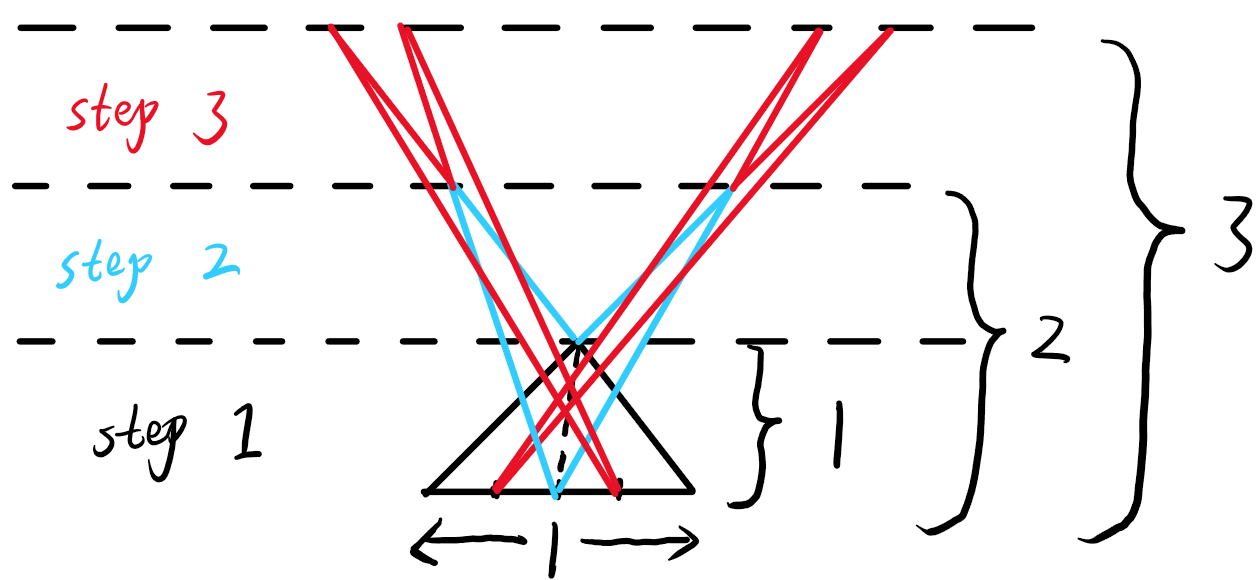}\end{center}
Step 1: initial triangle $\Delta_1$, $|\Delta_1|=1/2$.
Step 2: we obtain two triangles: $\Delta_{2,1},\Delta_{2,2}$; observe that
\[
\Delta_{2,1}\cup\Delta_{2,2}=\Delta_1\cup\text{``two arms''}.
\]
\begin{center}\includegraphics[width=.26\linewidth]{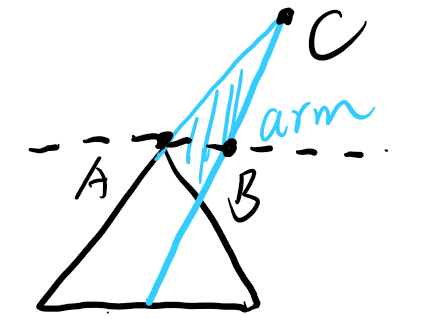}\end{center}
\[
\text{Area of arm}\leq2\cdot S_{\triangle ABC}
=2\cdot\frac12\cdot\underbrace{\left(\frac12\cdot\frac12\right)}_{\text{base}}
\cdot\underbrace{1}_{\text{height}}=\frac14.
\]
Therefore
\[
\text{Area of ``triangles''}:=|\Delta_{2,1}\cup\Delta_{2,2}|
\leq|\Delta_1|+2\cdot\text{Area of arm}\leq|\Delta_1|+\frac12.
\]
Step 3: we obtain four triangles: $\Delta_{3,1},\Delta_{3,2},\Delta_{3,3},\Delta_{3,4}$. Observe that
\[
\Delta_{3,1}\cup\Delta_{3,2}\cup\Delta_{3,3}\cup\Delta_{3,4}
=(\Delta_{2,1}\cup\Delta_{2,2})\cup\text{``four arms''}.
\]
\begin{center}\includegraphics[width=.26\linewidth]{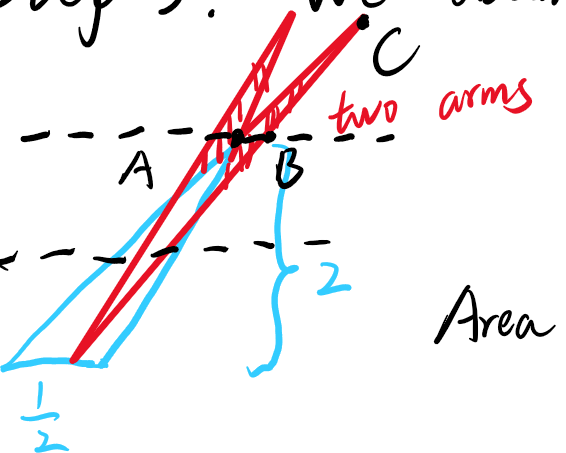}\end{center}
\[
\text{Area of each arm}\leq2\cdot S_{\triangle ABC}
=2\cdot\frac12\underbrace{\left(\frac13\cdot\frac1{2^2}\right)}_{\text{base}}
\cdot\underbrace{1}_{\text{height}}=\frac1{12}.
\]
\[
\text{Area of ``triangles''}=|\Delta_{2,1}\cup\Delta_{2,2}|+4\cdot\text{Area of each arm}
\leq|\Delta_1|+\frac12+\frac13.
\]
Iterate this process, we obtain:
Step $k$: $\#\{\text{small triangles}\}=2^{k-1}$, denoted by $\Delta_{k,i}$, $i=1,\ldots,2^{k-1}$.
\[
\text{Area of each arm}\leq2\cdot\frac12
\underbrace{\left(\frac1k\cdot\frac1{2^{k-1}}\right)}_{\text{base}}
\cdot\underbrace{1}_{\text{height}}=\frac1{k\cdot2^{k-1}}.
\]
\begin{align*}
\text{Area of ``triangles''}:=\left|\bigcup_{i=1}^{2^{k-1}}\Delta_{k,i}\right|
&\leq\left|\bigcup_{i=1}^{2^{k-2}}\Delta_{k-1,i}\right|
+2^{k-1}\cdot\frac1{k\cdot2^{k-1}}\\
&\leq|\Delta_1|+\frac12+\frac13+\cdots+\frac1k\lesssim\log k.
\end{align*}
We know for each $i$, $|\Delta_{k,i}|=\frac12\cdot\frac1{2^{k-1}}\cdot k=k\cdot2^{-k}$, so
\[
\sum_i|\Delta_{k,i}|=2^{k-1}\cdot k\cdot2^{-k}\sim k\gg\log k.
\]
It follows that the overlapping of these triangles is very large.
By dilation $k^{-1}$, we can identify each small triangle $\Delta_{k,i}$ with a tube with length $1$ of radius $2^{-k}/k$. And the volume of the union of these tubes $\sim\log k/k^2\to0$ as $k\to\infty$.
\begin{center}\includegraphics[width=.25\linewidth]{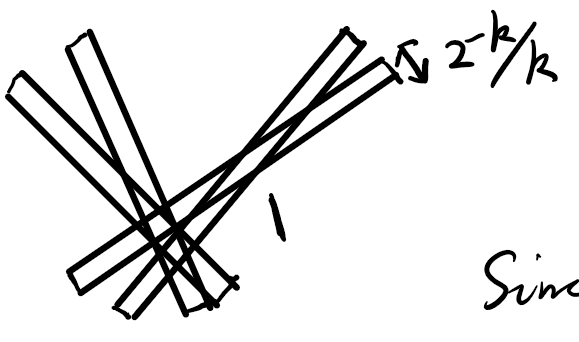}\end{center}
But, this construction doesn't give a Kakeya set. Since $\#\{\text{directions of ``triangles''}\}\sim2^k$, however, for a Kakeya set at scale $2^{-k}/k$, the number of directions should be $1/(2^{-k}/k)=k\cdot2^k\gg2^k$.
\subsection*{Perron tree}
Given a triangle $ABC$ with base length $2b$ and height $h$. Fix a parameter $\alpha\in(1/2,1)$ which will be chosen later.
Step 1:
\begin{center}\includegraphics[width=.92\linewidth]{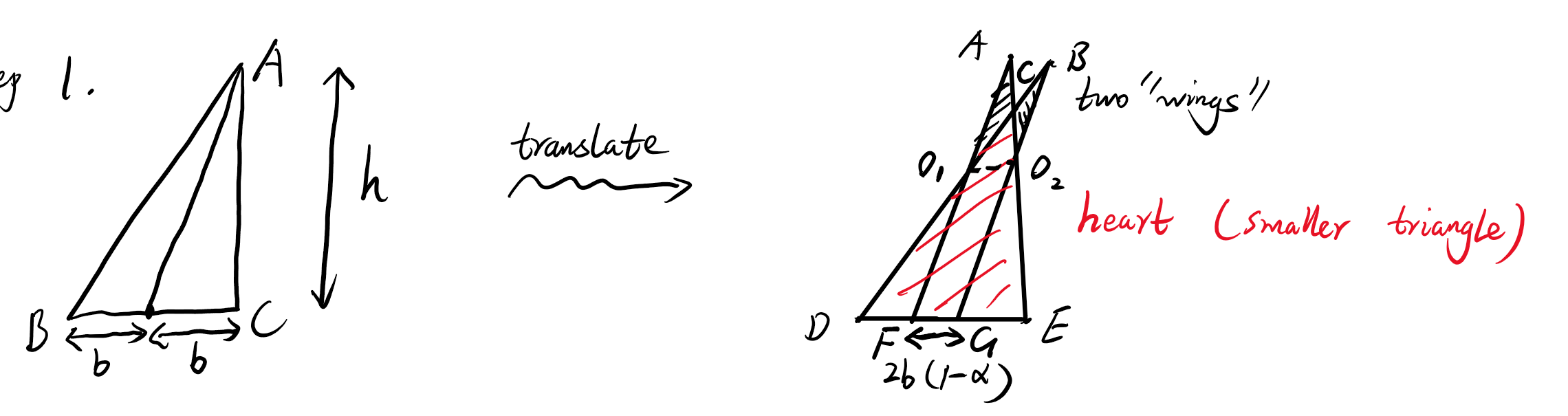}\end{center}
The heart is similar to $\triangle ABC$ and the similarity ratio is
\[
\frac{2b-2b(1-\alpha)}{2b}=\alpha.
\]
Therefore
\[
\boxed{\text{Area of heart}=\alpha^2\mathcal L^2(\triangle ABC).}
\]
Next we compute the area of two wings.
\begin{center}\includegraphics[width=.23\linewidth]{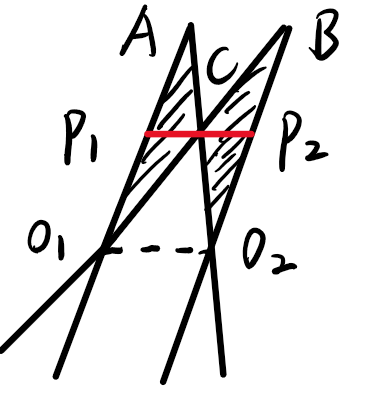}\end{center}
By translation, we know $O_1O_2\parallel DE$ and $|AO_1|=|BO_2|$.
Draw an auxiliary line $P_1P_2\parallel O_1O_2$,
\[
\triangle AO_1C\cong\triangle O_2BC
\Longrightarrow\triangle AP_1C\cong\triangle O_2P_2C
\text{ and }\triangle O_1P_1C\cong\triangle BP_2C.
\]
$\triangle AP_1C\sim\triangle AFE$ and the similarity ratio is $(h-\alpha h)/h=1-\alpha$. Therefore
\[
S_{\triangle AP_1C}=(1-\alpha)^2S_{\triangle AFE}
=\frac12(1-\alpha)^2\mathcal L^2(\triangle ABC).
\]
Similarly, $S_{\triangle BP_2C}=\frac12(1-\alpha)^2\mathcal L^2(\triangle ABC)$.
Thus
\[
\boxed{\text{Area of two ``wings''}=2(1-\alpha)^2\mathcal L^2(\triangle ABC).}
\]
Step 2:
\begin{center}\includegraphics[width=.94\linewidth]{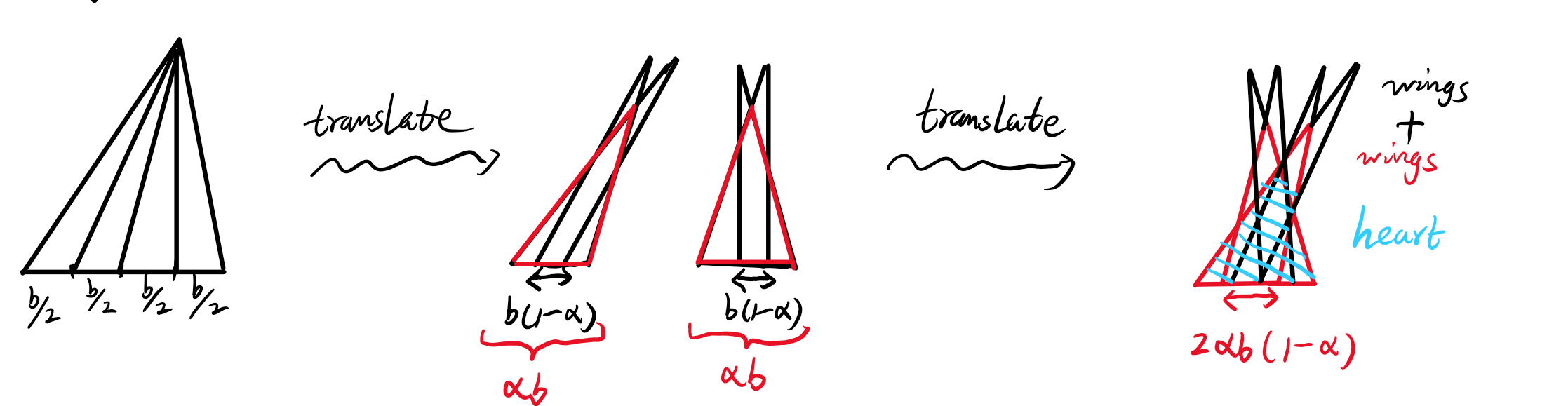}\end{center}
\[
\text{Area of heart}=\alpha^2\cdot2\mathcal L^2(\text{1-step heart})
=2\alpha^2\cdot\alpha^2\cdot\frac12\mathcal L^2(\triangle ABC).
\]
\[
\text{Area of wings}=\text{2-step wings}+(\text{1-step wings})\times2.
\]
\begin{align*}
\text{Area of wings}
&=2(1-\alpha)^2\cdot2\mathcal L^2(\text{1-step heart})
+2(1-\alpha)^2\cdot\frac12\mathcal L^2(\triangle ABC)\cdot2\\
&=2(1-\alpha)^2(\alpha^2+1)\mathcal L^2(\triangle ABC).
\end{align*}
Iterate this process, we obtain:
Step $k$:
\begin{center}\includegraphics[width=.94\linewidth]{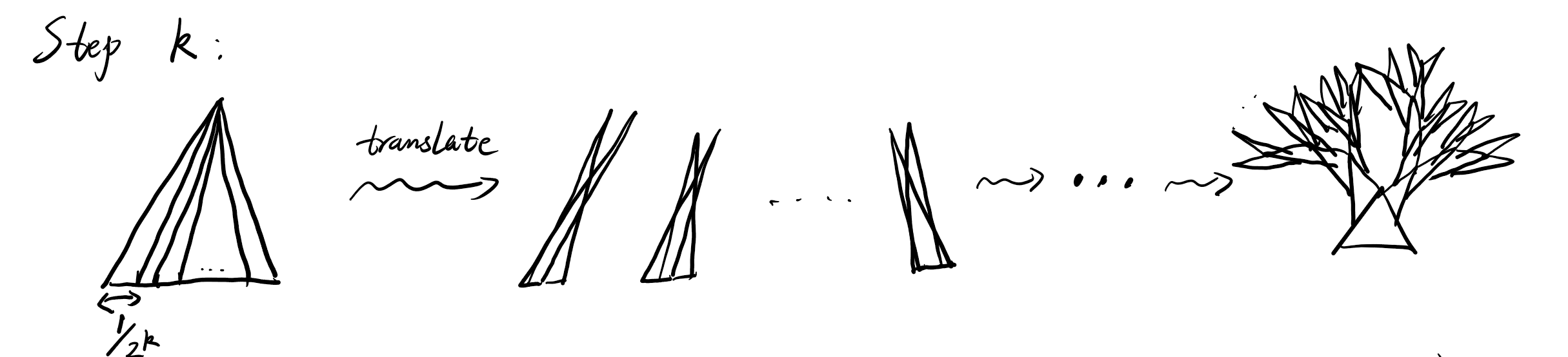}\end{center}
\begin{align*}
\text{Area of heart}
&=\alpha^2\cdot2\mathcal L^2(\text{$k-1$ heart})
=\alpha^2\cdot2\cdot\alpha^2\cdot2\mathcal L^2(\text{$k-2$ heart})\\
&=\cdots=(2\alpha^2)^{k-1}\mathcal L^2(\text{1 heart})\\
&=(2\alpha^2)^{k-1}\cdot\alpha^2\cdot\frac1{2^{k-1}}\mathcal L^2(\triangle ABC)
=\alpha^{2k}\mathcal L^2(\triangle ABC).
\end{align*}
(It follows that $2^i\mathcal L^2(\text{$k-i$ heart})=\alpha^{2(k-i)}\mathcal L^2(\triangle ABC)$.)
\begin{align*}
\text{Area of wings}
&=\text{$k$ wings}+(\text{$k-1$ wings})\times2+\cdots+(\text{1 wings})\times2^{k-1}\\
&=2(1-\alpha)^2\cdot2\mathcal L^2(\text{$k-1$ heart})
\end{align*}
\begin{align*}
&\hspace{8mm}+2(1-\alpha)^2\cdot2\mathcal L^2(\text{$k-2$ heart})\cdot2
+\cdots+2(1-\alpha)^2\cdot\frac1{2^{k-1}}\mathcal L^2(\triangle ABC)\cdot2^{k-1}\\
&=2(1-\alpha)^2\alpha^{2(k-1)}\mathcal L^2(\triangle ABC)\\
&\hspace{8mm}+2(1-\alpha)^2\alpha^{2(k-2)}\mathcal L^2(\triangle ABC)
+\cdots+2(1-\alpha)^2\alpha^2\mathcal L^2(\triangle ABC)
+2(1-\alpha)^2\mathcal L^2(\triangle ABC)\\
&=2(1-\alpha)^2\frac{1-\alpha^{2k}}{1-\alpha^2}\mathcal L^2(\triangle ABC)
=\frac{2(1-\alpha)(1-\alpha^{2k})}{1+\alpha}\mathcal L^2(\triangle ABC).
\end{align*}
In summary,
\begin{align*}
\text{Area of the resulting set}
&=\left(\alpha^{2k}+\frac{2(1-\alpha)(1-\alpha^{2k})}{1+\alpha}\right)\mathcal L^2(\triangle ABC)\\
&=\left(\frac{2(1-\alpha)}{1+\alpha}+
\left(1-\frac{2(1-\alpha)}{1+\alpha}\right)\alpha^{2k}\right)\mathcal L^2(\triangle ABC)\\
&=\left(\frac{2(1-\alpha)}{1+\alpha}+\frac{3\alpha-1}{1+\alpha}\alpha^{2k}\right)\mathcal L^2(\triangle ABC).
\end{align*}
Hence, for any $\varepsilon>0$, first choose $\alpha\to1^-$ such that $2(1-\alpha)/(1+\alpha)<\varepsilon/2$, and then choose $k\to\infty$ such that $\frac{3\alpha-1}{1+\alpha}\alpha^{2k}<\varepsilon/2$, we obtain
\[
\text{Area of the resulting set}<\varepsilon\mathcal L^2(\triangle ABC).
\]
Finally, note that such a figure certainly contains line segments of unit length in all directions at angles at least $60^\circ$ to the horizontal line. We only need to copy two such figures and rotate them by $60^\circ$, $120^\circ$, respectively, to obtain the following theorem.
\paragraph{Theorem.} Given any $\varepsilon>0$, there exists a compact set $E\subset\mathbb R^2$ which contains a line segment of unit length in every direction such that
\[
\mathcal L^2(E)<\varepsilon.
\]
\section{Decoupling for cone}
The cone decoupling theorem is due to Bourgain and Demeter~\cite{BD2015}; see also~\cite{Demeter2020}.

Consider the truncated cone in $\mathbb R^3$:
\[
 \operatorname{Cone}
 =\bigl\{(\xi_1,\xi_2,\sqrt{|\xi_1|^2+|\xi_2|^2}):
 1<\sqrt{|\xi_1|^2+|\xi_2|^2}<2\bigr\}.
\]
We divide the $\delta^2$-neighborhood of the cone,
$\mathcal N_{\operatorname{Cone}}(\delta^2)$, into disjoint thin planks
with dimensions $\sim 1\times\delta\times\delta^2$:
\[
 \mathcal N_{\operatorname{Cone}}(\delta^2)
 =\bigsqcup_{\theta\sim1\times\delta\times\delta^2}\theta.
\]
\begin{center}
\includegraphics[width=.58\linewidth]{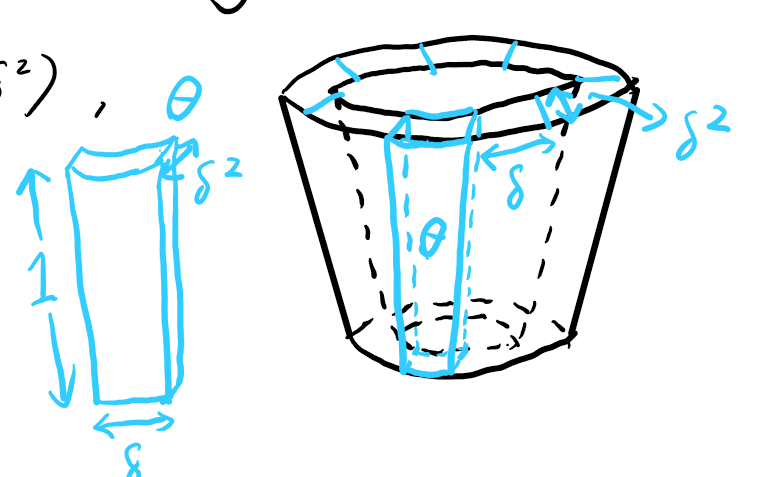}
\end{center}

\paragraph{Theorem.}
If $\operatorname{supp}(\widehat F)\subset
\mathcal N_{\operatorname{Cone}}(\delta^2)$, then
\[
 \|F\|_{L^p(\mathbb R^3)}
 \le D_p(\delta)
 \left(\sum_\theta\|F_\theta\|_{L^p(\mathbb R^3)}^2\right)^{1/2},
\]
where
\[
 D_p(\delta)\lesssim_\varepsilon
 \begin{cases}
  \delta^{-\varepsilon},&2\le p\le6,\\
  \delta^{-1/2+3/p-\varepsilon},&6\le p\le\infty.
 \end{cases}
\]

Without loss of generality, in the following argument, suppose
\[
 \operatorname{Cone}
 =\left\{\left(s,t,\frac{t^2}{s}\right):
 s\in\left[\frac12,1\right],\ t\in[-1,1]\right\}.
\]
For an interval $I\subset[-1,1]$, let
\[
 \Gamma_I(\delta)
 =\left\{\left(s,t,\frac{t^2}{s}+\delta'\right):
 s\in\left[\frac12,1\right],\
 \frac ts\in I,\ 0\le\delta'\le\delta^2\right\},
 \qquad
 \Gamma(\delta):=\Gamma_{[-1,1]}(\delta),
\]
and
\[
 F_I:=\bigl(\widehat F\,1_{\Gamma_I(\delta)}\bigr)^\vee.
\]
\begin{center}
\includegraphics[width=.35\linewidth]{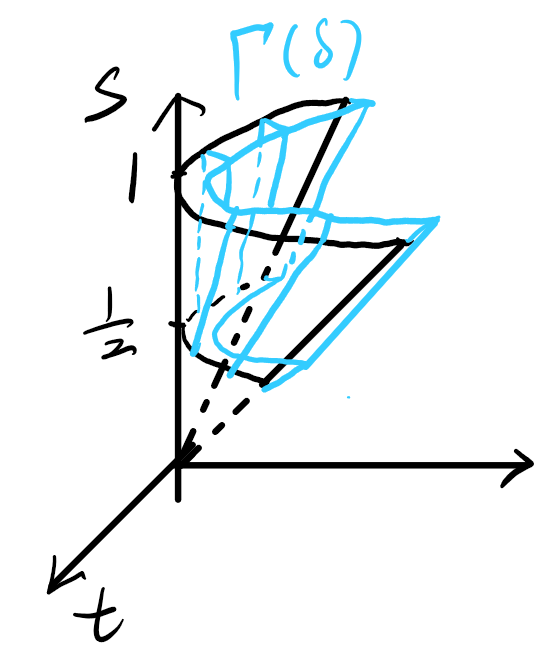}
\end{center}

\subsection*{Parabolic rescaling}
Given an interval $J\subset[-1,1]$ of length $\sigma>\delta$, recall
\[
 \Omega_J(\delta)
 =\{(\xi,|\xi|^2+\delta'):\xi\in J,\ |\delta'|\le\delta^2\}.
\]
Suppose $\operatorname{supp}(\widehat F)\subset\Omega_J(\delta)$.
Since $\Omega_J(\delta)\subset\Omega_{[-1,1]}(\delta)$, we can use
decoupling for the parabola directly:
\[
 \|F\|_p\le D_p(\delta)
 \left(\sum_{I\subset J}\|F_I\|_p^2\right)^{1/2}.
\]
\begin{center}
\includegraphics[width=.49\linewidth]{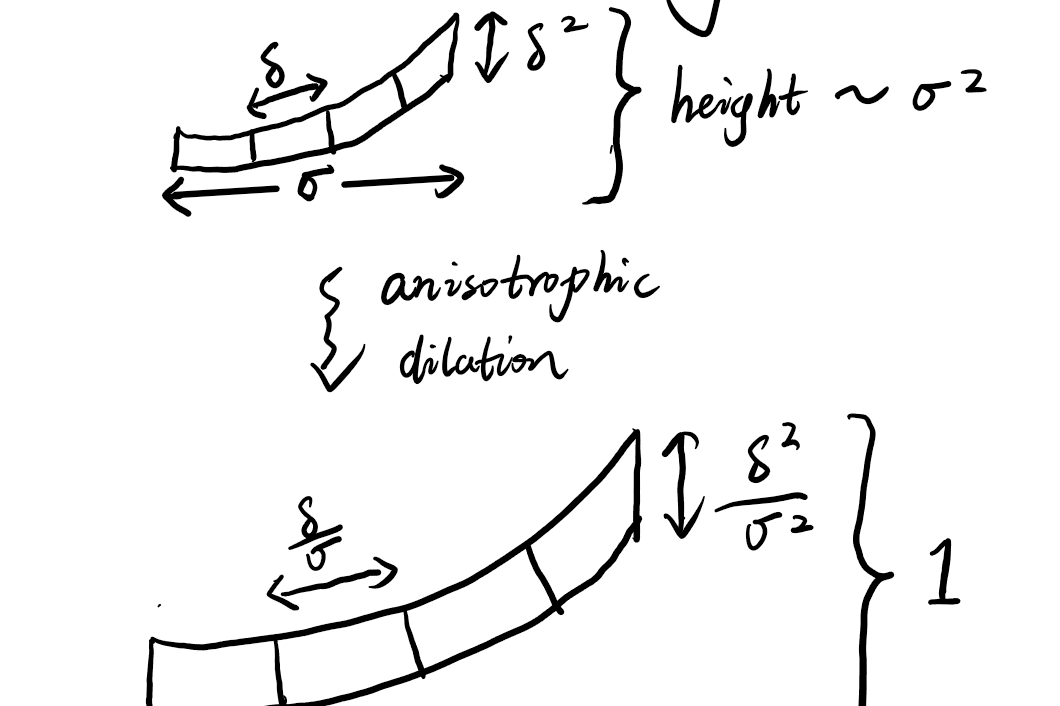}
\end{center}

But now $\widehat F$ is actually supported on a smaller part of the
parabola. In other words, we may decouple from $\sigma$ to $\delta$ rather
than from $1$ to $\delta$. By parabolic rescaling,
\[
 (\xi,\eta)\longmapsto
 \left(\frac{\xi-c_J}{\sigma},
 \frac{\eta-2c_J\xi+|c_J|^2}{\sigma^2}\right),
\]
this is equivalent to decoupling from $1$ to $\delta/\sigma$:
\[
 \|F\|_p\le D_p\left(\frac{\delta}{\sigma}\right)
 \left(\sum_{\substack{I\subset J\\\ell(I)=\delta}}
 \|F_I\|_p^2\right)^{1/2}.
\]

\subsection*{Decoupling for parabolic cylinder}
Consider a parabolic cylinder
\[
 \operatorname{Cyl}
 =\{(s,t,t^2):s\in[0,1],\ t\in[-1,1]\}
\]
and its $\delta^2$-neighborhood
\[
 \mathcal N_{\operatorname{Cyl}}(\delta^2)
 =\{(s,t,t^2+\delta'):s\in[0,1],\
 t\in[-1,1],\ 0\le\delta'\le\delta^2\}.
\]
Suppose
\[
 \operatorname{supp}(\widehat F)
 \subset\mathcal N_{\operatorname{Cyl}}(\delta^2)
 =\bigsqcup_{\substack{\theta:\text{planks}\\
 \theta\sim1\times\delta\times\delta^2}}\theta.
\]
\begin{center}
\includegraphics[width=.38\linewidth]{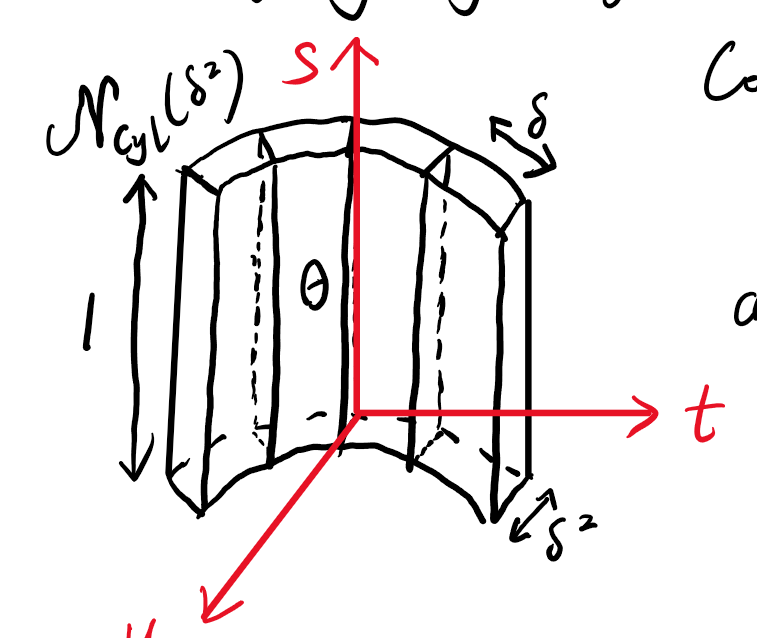}
\end{center}
By Fubini's theorem, we have

\[
 \|F\|_p^p=\int_{\mathbb R}
 \left(\int_{\mathbb R^2}|F(x,y,z)|^p\,dy\,dz\right)\,dx.
\]
For $(x,\xi,\eta)\in\mathbb R^3$, we have
\[
 (F(x,\cdot))^{\wedge_{2,3}}(\xi,\eta)
 =\bigl(\widehat F(\cdot,\xi,\eta)\bigr)^{\vee_1}(x).
\]
In fact, for any $\phi\in\mathcal S(\mathbb R^3)$,
\begin{align*}
 &\int_{\mathbb R^3}(F(x,\cdot))^{\wedge_{2,3}}(\xi,\eta)
 \phi(x,\xi,\eta)\,dx\,d\xi\,d\eta\\
 &\quad=\int_{\mathbb R^3}
 \left(\int_{\mathbb R^2}
 F(x,y,z)e^{-2\pi i(y\xi+z\eta)}\,dy\,dz\right)
 \phi(x,\xi,\eta)\,dx\,d\xi\,d\eta\\
 &\quad=\int_{\mathbb R^3}\int_{\mathbb R^2}
 F(x,y,z)e^{-2\pi i(y\xi+z\eta)}
 \left(\int_{\mathbb R^3}\check\phi(u,v,w)
 e^{-2\pi i(ux+v\xi+w\eta)}\,du\,dv\,dw\right)
 \,dy\,dz\,dx\,d\xi\,d\eta\\
 &\quad=\int_{\mathbb R^5}
 \widehat F(u,\xi,\eta)\check\phi(u,v,w)
 e^{-2\pi i(v\xi+w\eta)}\,du\,dv\,dw\,d\xi\,d\eta.
\end{align*}

Continuing the same calculation,
\begin{align*}
 &\int_{\mathbb R^5}
 \left(\int_{\mathbb R}
 \bigl(\widehat F(\cdot,\xi,\eta)\bigr)^{\vee_1}(x)
 e^{-2\pi iux}\,dx\right)
 \check\phi(u,v,w)e^{-2\pi i(v\xi+w\eta)}
 \,du\,dv\,dw\,d\xi\,d\eta\\
 &\qquad=
 \int_{\mathbb R^3}
 \bigl(\widehat F(\cdot,\xi,\eta)\bigr)^{\vee_1}(x)
 \phi(x,\xi,\eta)\,dx\,d\xi\,d\eta.
\end{align*}
Therefore, for fixed $x$, if
$(\xi,\eta)\notin
\pi_{tu\text{-plane}}(\mathcal N_{\operatorname{Cyl}}(\delta^2))$,
then
\[
 (F(x,\cdot))^{\wedge_{2,3}}(\xi,\eta)
 =\int\widehat F(s,\xi,\eta)e^{2\pi isx}\,ds=0.
\]
It follows that
\[
 \operatorname{supp}\bigl((F(x,\cdot))^{\wedge_{2,3}}\bigr)
 \subset\pi_{tu\text{-plane}}
 \bigl(\mathcal N_{\operatorname{Cyl}}(\delta^2)\bigr)
 =\mathcal N_{\mathbb P^1}(\delta^2).
\]
Thus, continuing the estimate of $\|F\|_p^p$, by two-dimensional
decoupling and then Minkowski,
\begin{align*}
 \|F\|_p^p
 &\le\int_{\mathbb R}
 \left[
 D_p^{(2)}(\delta)
 \left(\sum_{\substack{I\subset\mathbb R_t\\\ell(I)=\delta}}
 \|F(x,\cdot)_I\|_{L^p}^2\right)^{1/2}
 \right]^p\,dx\\
 &\le \bigl(D_p^{(2)}(\delta)\bigr)^p
 \left(\sum_I
 \left(\int_{\mathbb R}\|F(x,\cdot)_I\|_p^p\,dx\right)^{2/p}
 \right)^{p/2}.
\end{align*}

Here $(\mathbb R\times\Omega_I(\delta))\cap
\operatorname{supp}(\widehat F)\simeq\theta$.

The integral in the last line is
\begin{align*}
 \int_{\mathbb R}\|F(x,\cdot)_I\|_p^p\,dx
 &=\int_{\mathbb R}\int_{\mathbb R^2}
 \left|
 \left((F(x,\cdot))^{\wedge_{2,3}}
 1_{\Omega_I(\delta)}\right)^{\vee_{2,3}}(y,z)
 \right|^p\,dy\,dz\,dx\\
 &=\int_{\mathbb R^3}|F_\theta|^p.
\end{align*}

Hence the preceding estimate becomes
\[
 \bigl(D_p^{(2)}(\delta)\bigr)^p
 \left(\sum_\theta\|F_\theta\|_p^2\right)^{p/2}.
\]

In summary, decoupling for a three-dimensional parabolic cylinder is
essentially equivalent to decoupling for a two-dimensional parabola:
\[
 \|F\|_{L^p(\mathbb R^3)}
 \le D_p^{(2)}(\delta)
 \left(\sum_\theta\|F_\theta\|_{L^p(\mathbb R^3)}^2\right)^{1/2}.
\]

\subsection*{Proof of the theorem}
We first prove the estimate for $2\le p\le6$; the range $p\ge6$ then follows by interpolation with the trivial $L^\infty$ estimate.
First, we cover the $st$-plane by a family of disjoint cubes of length
$\sim\delta^\varepsilon$. Since the number of cubes is
$\sim\delta^{-2\varepsilon}$ and we allow a $\delta^{-2\varepsilon}$
loss, it is sufficient to prove the theorem on each cube. So now consider
\[
 \left\{\left(s,t,\frac{t^2}{s}\right):
 1-\delta^\varepsilon\le s\le1,\
 0\le\frac ts\le\delta^\varepsilon\right\}.
\]
\begin{center}
\includegraphics[width=.25\linewidth]{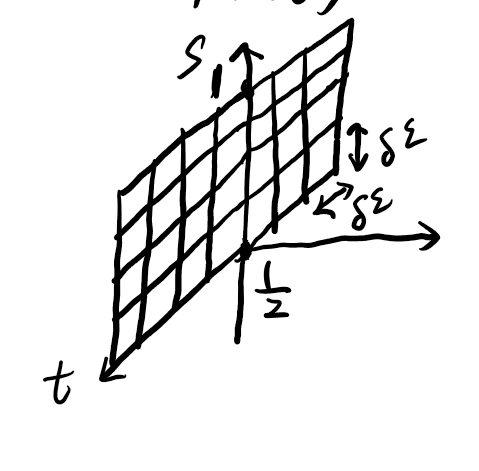}
\end{center}

Second, observe that on the scale $\delta^\varepsilon$, the parabolic
cylinder is very close to the cone. Recall the parabolic cylinder
\[
 \left\{(s,t,t^2):
 1-\delta^\varepsilon\le s\le1,\
 0\le\frac ts\le\delta^\varepsilon\right\}.
\]
We have
\[
 \left|\left(s,t,\frac{t^2}{s}\right)-(s,t,t^2)\right|
 =|t|^2\frac{|1-s|}{|s|}
 \lesssim\delta^{3\varepsilon}.
\]
Therefore, the cone is contained in an $O(\delta^{3\varepsilon})$-neighborhood
of a parabolic cylinder.
\begin{center}
\includegraphics[width=.95\linewidth]{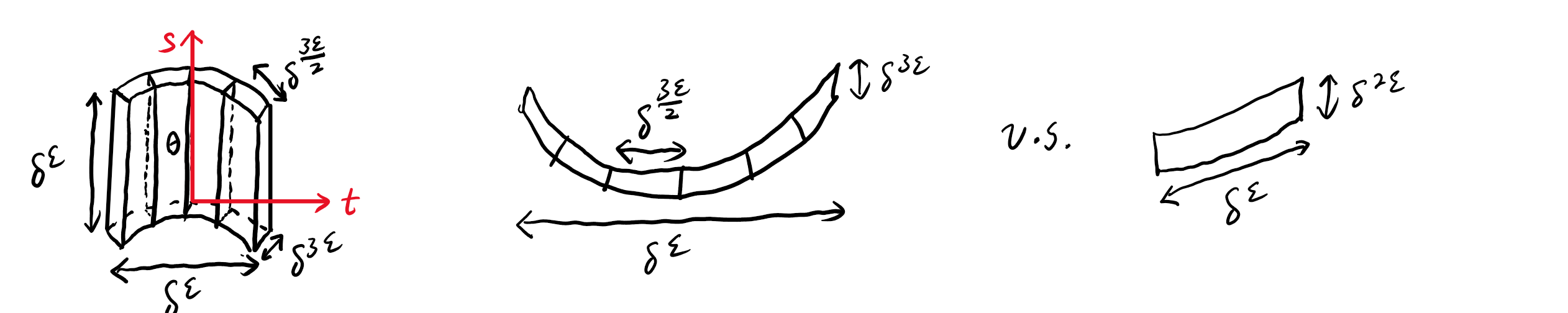}
\end{center}
Since $\delta^{3\varepsilon}\ll(\delta^\varepsilon)^2$, the cylinder is
still \emph{curved}. So we can continue to divide it into smaller planks
with dimensions $\sim\delta^\varepsilon\times
\delta^{3\varepsilon/2}\times\delta^{3\varepsilon}$.
Apply decoupling for the parabolic cylinder and parabolic rescaling;
we obtain
\[
 \|F\|_p\le
 D_p^{(2)}\left(\frac{\delta^{3\varepsilon/2}}{\delta^\varepsilon}\right)
 \left(\sum_{\ell(I)=\delta^{3\varepsilon/2}}
 \|F_I\|_p^2\right)^{1/2}.
\]

Thus
\[
 \|F\|_p\lesssim_\varepsilon
 \delta^{-\varepsilon^2/2}
 \left(\sum_{\ell(I)=\delta^{3\varepsilon/2}}
 \|F_I\|_p^2\right)^{1/2}.
\]

Third, for each $I$ of length $\delta^{3\varepsilon/2}$, we see that
$\operatorname{supp}(\widehat F_I)$ is contained in a cone:
\[
 \left\{\left(s,t,\frac{t^2}{s}\right):
 1-\delta^\varepsilon\le s\le1,\
 0\le\frac ts\le\delta^{3\varepsilon/2}\right\}.
\]
Similarly, we may approximate this cone by a parabolic cylinder
$\{(s,t,t^2)\}$, and
\[
 \left|\left(s,t,\frac{t^2}{s}\right)-(s,t,t^2)\right|
 =|t|^2\frac{|1-s|}{|s|}
 \lesssim\delta^{4\varepsilon}
 \ll(\delta^{3\varepsilon/2})^2.
\]
Again, we can decouple from $\delta^{3\varepsilon/2}$ to
$\delta^{2\varepsilon}$:
\begin{align*}
 \|F_I\|_p
 &\le D_p^{(2)}
 \left(\frac{\delta^{2\varepsilon}}{\delta^{3\varepsilon/2}}\right)
 \left(\sum_{\ell(I_1)=\delta^{2\varepsilon}}
 \|F_{I_1}\|_p^2\right)^{1/2}\\
 &\lesssim_\varepsilon
 \delta^{-\varepsilon^2/2}
 \left(\sum\|F_{I_1}\|_p^2\right)^{1/2}.
\end{align*}
Iterate this process. To sum up, we decouple at these scales:
\begin{align*}
 t^2\text{-direction}:&\quad
 \delta^{3\varepsilon}\longrightarrow
 \delta^{4\varepsilon}\longrightarrow
 \delta^{5\varepsilon}\longrightarrow\cdots
 \longrightarrow\delta^{(3+k)\varepsilon}
 \longrightarrow\cdots\delta^{2+\varepsilon},\\
 (t/s)\text{-direction}:&\quad
 \delta^\varepsilon\longrightarrow
 \delta^{3\varepsilon/2}\longrightarrow
 \delta^{2\varepsilon}\longrightarrow\cdots
 \longrightarrow\delta^{(1+k/2)\varepsilon}
 \longrightarrow\cdots\delta.
\end{align*}
The passage from the upper row to the lower row takes a square root,
and the passage from the lower row to the next upper term squares and
multiplies by $\delta^\varepsilon$.

In the $k$-th step, we have a loss
\[
 D_p^{(2)}
 \left(
 \frac{\delta^{(3/2+k/2)\varepsilon}}
 {\delta^{(1+k/2)\varepsilon}}
 \right)
 \lesssim_\varepsilon\delta^{-\varepsilon^2/2},
\]
and we will decouple $\sim O(1/\varepsilon)$ times. So we lose a total of
\[
 \left(C_\varepsilon\delta^{-\varepsilon^2/2}\right)^{1/\varepsilon}
 =C_\varepsilon\delta^{-\varepsilon/2}
 \le C_\varepsilon\delta^{-\varepsilon},
\]
as desired.

\section{Bourgain's circular maximal function}

Let $\chi_0:\mathbb R\to\mathbb R_{\ge0}$ be a smooth bump function
supported on $[-1,1]$ and satisfying $\chi_0(s)=1$ for any
$s\in[-1/2,1/2]$. For $u\in[1,2]$, define the averaging operator
\[
 A_uf(x,y)=\int_{\mathbb R}f(x-us,y-us^2)\chi_0(s)\,ds
\]
and the maximal function
\[
 Mf(x,y)=\sup_{u\in[1,2]}|A_uf(x,y)|.
\]
Here $A_uf$ is the integral of $f$ along the parabola
$\mathbb P^1_{u,(x,y)}$.
\begin{center}
\includegraphics[width=.4\linewidth]{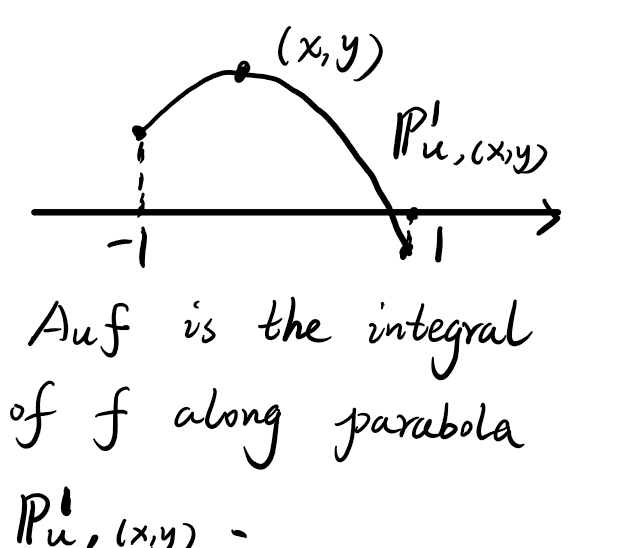}
\end{center}

\paragraph{Theorem 1 (Bourgain, '86).} \cite{Bourgain1986}
\[
 \|Mf\|_{L^p(\mathbb R^2)}
 \lesssim_p\|f\|_{L^p(\mathbb R^2)}
 \qquad\text{if }p>2.
\]
\paragraph{Remark.}
In higher dimensions $n\ge3$, the same result for the spherical maximal
operator with $p>n/(n-1)$ is due to Stein, '76. See~\cite{Stein1976}.

\subsection*{Necessity of the condition \texorpdfstring{$p>2$}{p>2}}

\paragraph{Example 1: $p<2$.}
Consider a $\delta$-ball $B$ centered at the origin. Take $f=1_B$.
Choose a point $(1,1)$, and observe that the parabola
$(1,1)-(s,s^2)$ passes through the origin when $s=1$ and the length of
$B\cap\text{parabola}$ is $\sim\delta$. So
\[
 Mf(1,1)\ge\int 1_B(1-s,1-s^2)\,ds\sim\delta.
\]
For $(x,y)\in B((1,1),1/100)$, there is $u\in[0.9,1.1]$ such that
the parabola $(x,y)-u(s,s^2)$ passes through the origin. In fact,
\[
 \exists s\text{ such that }(x,y)-u(s,s^2)=(0,0)
 \quad\Longleftrightarrow\quad
 \begin{cases}x=us,\\y=us^2,\end{cases}
 \quad\Longrightarrow\quad u=\frac{x^2}{y}.
\]
Note that $0.99\le x,y\le1.01$; we have $0.9\le u\le1.1$.

Then the length of $B\cap\text{parabola}$ is $\sim\delta$. So
\[
 Mf(x,y)\ge\int1_B(x-us,y-us^2)\,ds\sim\delta,
\]
and consequently
\[
 \|Mf\|_{L^p(\mathbb R^2)}
 \ge\|Mf\|_{L^p(B((1,1),1/100))}\gtrsim\delta.
\]
\begin{center}
\includegraphics[width=.38\linewidth]{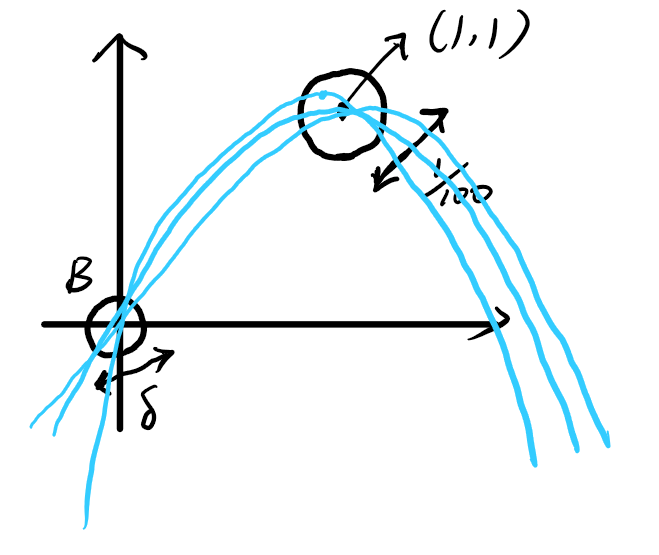}
\end{center}

Since $\|f\|_{L^p(\mathbb R^2)}\sim\delta^{2/p}$, we get
\[
 \frac{\|Mf\|_{L^p(\mathbb R^2)}}{\|f\|_{L^p(\mathbb R^2)}}
 \gtrsim\delta^{1-2/p}
 \quad\Longrightarrow\quad p\ge2.
\]

\paragraph{Example 2: $p=2$ (Stein's example).}
Take
\[
 f(x)=\frac1{|x|\log(1/|x|)}
 1_{\{|x|\le1/2\}}(x).
\]
Then
\begin{align*}
 \int|f|^2
 &=2\pi\int_0^{1/2}\frac1{r^2\log^2(1/r)}\,r\,dr\\
 &=2\pi\int_0^{1/2}\frac{-d\log(1/r)}{\log^2(1/r)}
 =2\pi\int_{\log2}^{\infty}\frac{dt}{t^2}<\infty.
\end{align*}
We claim $Mf(x)=\infty$ for \emph{many} $x$. Observe that, when
$s\to1$, the parabola $(1-s,1-s^2)$ is approximately the line
$(1-s,2(1-s))$, and
\[
 |(1-s,2(1-s))|\sim|1-s|=:r.
\]
\begin{center}
\includegraphics[width=.35\linewidth]{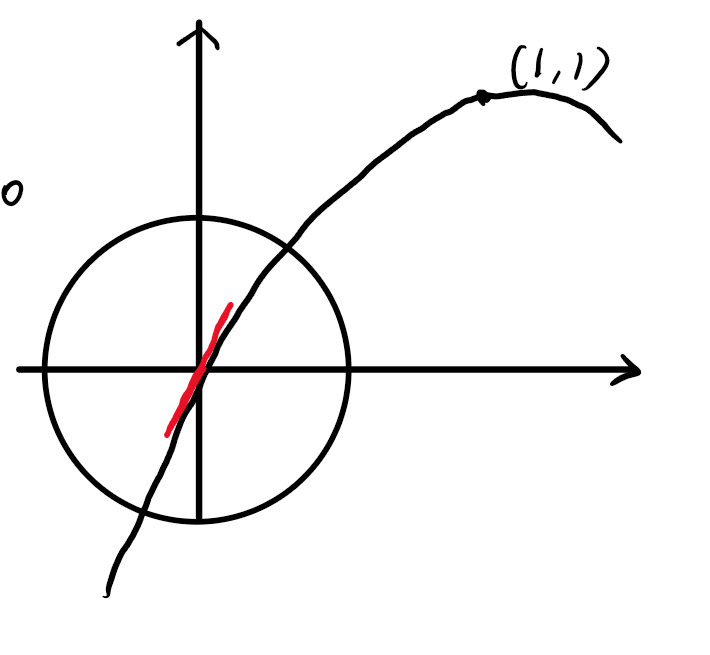}
\end{center}
So, directly compute
\[
 Mf(1,1)\ge\int f(1-s,1-s^2)\,ds.
\]

Continuing,
\[
 Mf(1,1)\ge\int_0^{1/100}\frac1{r\log(1/r)}\,dr=\infty.
\]
Similarly, on $B((1,1),1/100)$, we still have $Mf(x,y)=\infty$.
It follows that $\|Mf\|_{L^2(\mathbb R^2)}=\infty$.

\subsection*{Connection to geometric measure theory}

\paragraph{Theorem 2 (Marstrand, '86).} \cite{Marstrand1987}
If $E\subset\mathbb R^2$ is a Borel set and satisfies that for any
$x\in\mathbb R^2$ there exists $r_x\in[1,2]$ such that the circle
$C(x,r_x)\subset E$, then $\mathcal L^2(E)>0$.

\paragraph{Remark.}
This theorem shows that circular Nikodym sets \emph{do not exist}!

\paragraph{Proof of Theorem 2 by Theorem 1.}
Suppose $f=1_E$. For any $x\in\mathbb R^2$, by the definition of $E$,
\[
 Mf(x)\ge\int_{C(x,r_x)\subset E}1_E(s)\,d\ell(s)
 =2\pi r_x\sim1.
\]

Here we modify the definition of $M$ to integrate $f$ along the circle
$C(x,r)$. Hence,
\[
 1\lesssim\|Mf\|_p\lesssim_p\|f\|_p
 =|E|^{1/p}
 \quad\Longrightarrow\quad |E|>0.
 \qquad\square
\]

\subsection*{Proof of Theorem 1}

Littlewood--Paley decomposition:
\[
 f=P_{\le0}f+\sum_{k\ge1}P_kf,
\]
where each $P_kf$ has Fourier transform supported on the annulus
$|(\xi,\eta)|\sim2^k$.

\paragraph{Low frequency $P_{\le0}f$.}
Since $\operatorname{supp}((P_{\le0}f)^\wedge)\subset B(0,1)$, on
scale $1$ we know $|P_{\le0}f|$ is essentially constant. Hence
$|M(P_{\le0}f)|\approx|P_{\le0}f|$.
More precisely, choose a Schwartz function $\varphi$ satisfying
$\operatorname{supp}(\widehat\varphi)\subset B(0,1)$ and
$\varphi\approx1$ on $B(0,1)$. 
Then, for each $u\in[1,2]$,
\[
 A_u(P_{\le0}f)
 =\bigl((P_{\le0}f)^\wedge\cdot\text{multiplier}
 \cdot\widehat\varphi\bigr)^\vee
 =P_{\le0}f*A_u\varphi
\]
and
\[
 |A_u\varphi(x)|\lesssim_N(1+|x|)^{-N}
 \quad\text{uniformly in }u.
\]

Consequently,
\[
 |A_u(P_{\le0}f)(x)|\lesssim M_{\mathrm{HL}}(P_{\le0}f)(x),
\]
where $M_{\mathrm{HL}}$ is the Hardy--Littlewood maximal function.
Therefore,
\[
 \|M(P_{\le0}f)\|_p
 \lesssim\|M_{\mathrm{HL}}(P_{\le0}f)\|_p
 \lesssim\|P_{\le0}f\|_p
 \lesssim\|f\|_p,
\]
as desired.

\paragraph{High frequency.}
We aim to prove
\[
 \|M(P_kf)\|_p\lesssim2^{-\gamma_pk}\|f\|_p
 \quad\text{for some }\gamma_p>0.
\]
First of all, note that
\begin{align*}
 \widehat{A_u(P_kf)}(\xi,\eta)
 &=\int_{\mathbb R^2}\int_{\mathbb R}
 (P_kf)(x-us,y-us^2)\chi_0(s)
 e^{-i(x\xi+y\eta)}\,ds\,dx\,dy\\
 &=\widehat{P_kf}(\xi,\eta)
 \int_{\mathbb R}e^{-iu(s\xi+s^2\eta)}\chi_0(s)\,ds.
\end{align*}
Because $\widehat{P_kf}$ is now localized on $|(\xi,\eta)|\sim2^k$
and $|s|\le1$, the smallest scale on which we can distinguish $u$
is $2^{-k}$. Therefore we can choose $u$ as some dyadic number:
\[
 u_j\in[1+(j-1)2^{-k},1+j2^{-k}],\qquad1\le j\le2^k.
\]

Then
\[
 \sup_{u\in[1,2]}|A_uP_kf|
 \sim\sup_{1\le j\le2^k}|A_{u_j}P_kf|
 \lesssim\left(\sum_j|A_{u_j}P_kf|^p\right)^{1/p}.
\]
Now, if we use van der Corput's lemma, we get (as we did in the first
section)
\[
 \left|\int_{\mathbb R}
 e^{-iu_j(s\xi+s^2\eta)}\chi_0(s)\,ds\right|
 \lesssim(2^{-k})^{1/2}.
\]
By Plancherel's theorem,
\begin{align*}
 \|A_{u_j}P_kf\|_{L^2}
 &=\left\|\widehat{P_kf}\cdot
 \int_{\mathbb R}e^{-iu_j(s\xi+s^2\eta)}\,ds\right\|_{L^2}\\
 &\lesssim2^{-k/2}\|\widehat{P_kf}\|_2
 =2^{-k/2}\|P_kf\|_2.
\end{align*}
Hence
\begin{align*}
 \|M(P_kf)\|_{L^2}
 &\lesssim\left(\sum_{j=1}^{2^k}
 \|A_{u_j}P_kf\|_{L^2}^2\right)^{1/2}\\
 &\lesssim2^{k/2}\cdot2^{-k/2}\|P_kf\|_2
 =\|P_kf\|_2.
\end{align*}
It follows that $\gamma_2=0$, but we need $\gamma_p>0$!

\subsection*{Sobolev embedding}
Suppose $\operatorname{supp}(\widehat f)
\subset\{x\in\mathbb R:|x|\sim2^k\}$.
Then, for any $p\ge1$,
\[
 \sup|f(x)|\lesssim2^{k/p}\|f\|_p.
\]
The picture illustrates that, on an interval of length $2^{-k}$,
$|f|\ge\frac12\|f\|_\infty$, so
$\|f\|_p\gtrsim\|f\|_\infty\,2^{-k/p}$.
\begin{center}
\includegraphics[width=.52\linewidth]{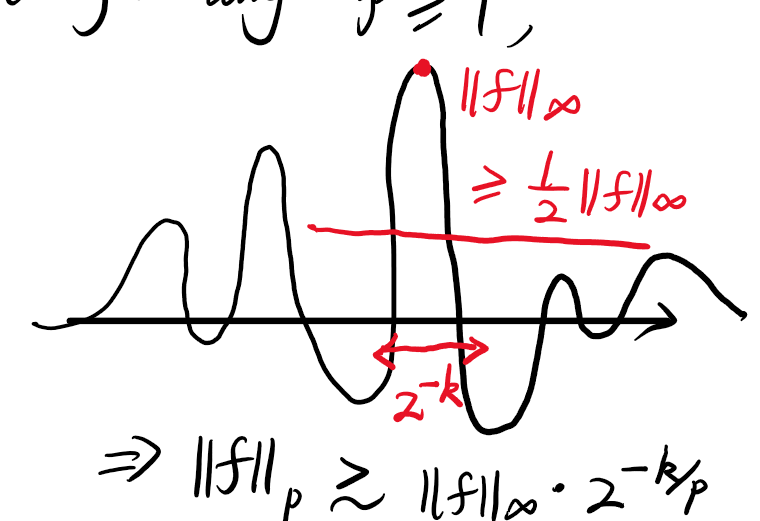}
\end{center}
Consequently,
\[
 |A_{u_j}(P_kf)|
 \lesssim2^{k/p}\|A_uP_kf\|_{L^p_u([1,2])},
\]
and
\[
 \sup_{u\in[1,2]}|A_u(P_kf)|
 \lesssim2^{k/p}\|A_uP_kf\|_{L^p_u([1,2])}.
\]

\subsection*{Local smoothing estimate}
For the connection between local smoothing and circular maximal estimates, see~\cite{MSS1992}; the argument below uses~\cite{BD2015}.
\[
 \|A_uP_kf\|_{L^p(\mathbb R^2\times[1,2])}
 \lesssim2^{-k/p}2^{-\gamma_pk}
 \|P_kf\|_{L^p(\mathbb R^2)}
\]
for every $p>2$ and some $\gamma_p>0$.
The factor $2^{-\gamma_pk}$ means that we obtain more gain than before.
Rather than just using the decay of
$\int_{\mathbb R}e^{-iu(t\xi+t^2\eta)}\chi_0(t)\,dt$,
we should pay more attention to its oscillation. By the stationary phase
principle,
\begin{align*}
 \int_{\mathbb R}e^{-iu(t\xi+t^2\eta)}\chi_0(t)\,dt
 &\approx|\eta|^{-1/2}
 e^{-iu(-\xi^2/(2\eta)+(\xi^2/(4\eta^2))\eta)}\\
 &=2^{-k/2}e^{iu\xi^2/(4\eta)}.
\end{align*}

If $|\xi|\ge10|\eta|$, we have better decay:
\[
 \left|\int_{\mathbb R}e^{-iu(t\xi+t^2\eta)}
 \chi_0(t)\,dt\right|\lesssim_N2^{-kN}.
\]
Hence, in the following argument, we assume $|\xi|\le10|\eta|$,
which implies $|\eta|\sim2^k$. Consequently,
\begin{align*}
 A_uP_kf(x,y)
 &\approx2^{-k/2}\int_{\mathbb R^2}
 \widehat{P_kf}(\xi,\eta)e^{iu\xi^2/(4\eta)}
 e^{i(x\xi+y\eta)}\,d\xi\,d\eta\\
 &=2^{-k/2}E_{\operatorname{Cone}}(\widehat{P_kf})(x,y,u).
\end{align*}
Here the phase is
$\langle(x,y,u),(\xi,\eta,\xi^2/(4\eta))\rangle$,
and $(\xi,\eta,\xi^2/(4\eta))$ lies on a cone.
This reminds us that we can use decoupling for a truncated cone.
\begin{center}
\includegraphics[width=.93\linewidth]{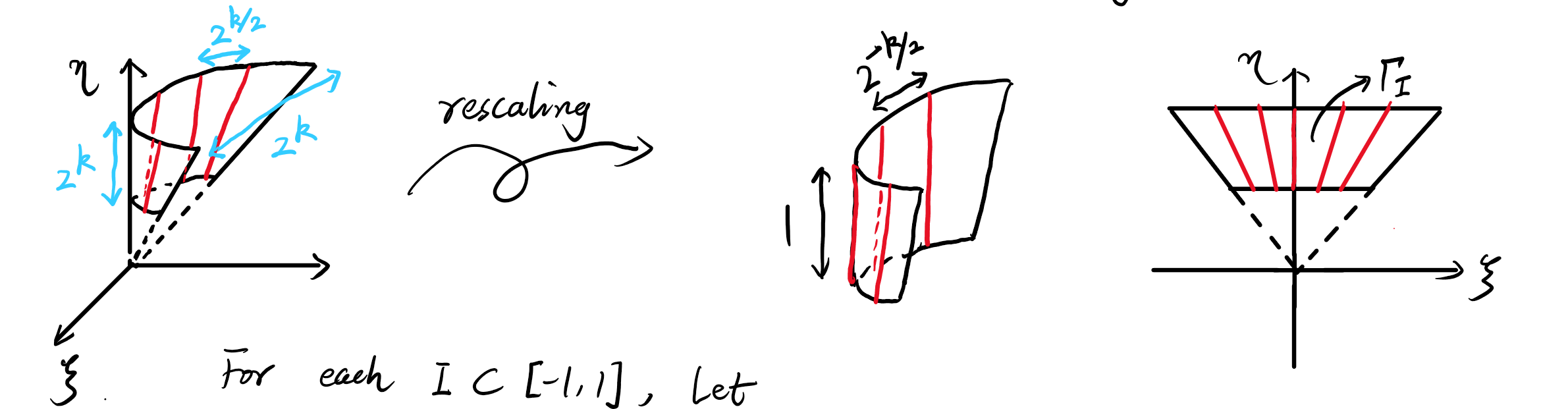}
\end{center}
For each $I\subset[-1,1]$, let

\[
 \Gamma_I:=\{(\xi,\eta):\xi/\eta\in I,\ \eta\sim2^k\},
 \qquad
 P_{k,I}f:=\bigl(\widehat{P_kf}\,1_{\Gamma_I}\bigr)^\vee.
\]
For $p\ge6$, apply decoupling for the cone; we obtain
\begin{align*}
 \|A_uP_kf\|_p
 &\approx2^{-k/2}\|E_{\operatorname{Cone}}(\widehat{P_kf})\|_p\\
 &\lesssim2^{-k/2}
 (2^{-k/2})^{-1/2+3/p-\varepsilon}
 \left(\sum_{\substack{I\subset[-1,1]\\\ell(I)=2^{-k/2}}}
 \|E_{\operatorname{Cone}}(\widehat{P_{k,I}f})\|_p^2
 \right)^{1/2}.
\end{align*}
For each $I$,
\[
 E_{\operatorname{Cone}}(\widehat{P_{k,I}f})(x,y,u)
 =\int_{\Gamma_I}\widehat{P_kf}(\xi,\eta)
 e^{i(x\xi+y\eta+u\xi^2/(4\eta))}\,d\xi\,d\eta.
\]
On $\Gamma_I$, we know $|\xi/\eta-c_I|\le2^{-k/2}$ and
$|\eta|\sim2^k$. Thus
\begin{align*}
 \frac{\xi^2}{4\eta}
 &=\frac14\left[
 \left(\frac{\xi}{\eta}-c_I\right)^2
 +2c_I\frac{\xi}{\eta}-c_I^2
 \right]\eta\\
 &=\frac14\left(\frac{\xi}{\eta}-c_I\right)^2\eta
 +\frac12c_I\xi-\frac14c_I^2\eta.
\end{align*}

Consequently,
\[
 e^{i(x\xi+y\eta+u\xi^2/(4\eta))}
 =
 e^{\frac14iu(\xi/\eta-c_I)^2\eta}
 e^{i[(x+\frac12c_Iu)\xi+(y-\frac14c_I^2u)\eta]},
\]
and
\[
 \left|u\left(\frac{\xi}{\eta}-c_I\right)^2\eta\right|
 \lesssim(2^{-k/2})^2\,2^k=1.
\]
On $\Gamma_I$, $E_{\operatorname{Cone}}$ is essentially flat, and we
remove the curvature term by using decoupling!

Without loss of generality, assume all $c_I=0$. Then
\[
 E_{\operatorname{Cone}}(\widehat{P_{k,I}f})\approx P_{k,I}f,
\]
and therefore
\[
 \|A_u(P_kf)\|_p
 \lesssim2^{-k/2}2^{(k/2)(1/2-3/p+\varepsilon)}
 \left(\sum_I\|P_{k,I}f\|_p^2\right)^{1/2}.
\]
For $p=2$, by Plancherel's theorem,
\[
 \left(\sum_I\|P_{k,I}f\|_2^2\right)^{1/2}
 =\left\|\sum_I P_{k,I}f\right\|_2
 =\|P_kf\|_2.
\]
For $p=\infty$, by Young's inequality,
\[
 \left(\sum_I\|P_{k,I}f\|_\infty^2\right)^{1/2}
 \sim
 \left(\sum_I\|P_kf*\check1_{\Gamma_I}\|_\infty^2\right)^{1/2}.
\]

Continuing,
\begin{align*}
 \left(\sum_I\|P_{k,I}f\|_\infty^2\right)^{1/2}
 &\lesssim
 \left(\sum_I\|\check1_{\Gamma_I}\|_1^2\right)^{1/2}
 \|P_kf\|_\infty\\
 &\lesssim2^{k/4}\|P_kf\|_\infty.
\end{align*}
Here
\[
 \check1_{\Gamma_I}\approx|\Gamma_I|1_{\Gamma_I^*}
 \quad\Longrightarrow\quad
 \|\check1_{\Gamma_I}\|_1\lesssim1,
 \qquad \#I\sim2^{k/2}.
\]
Interpolation between $2$ and $\infty$ gives
\[
 \left(\sum_I\|P_{k,I}f\|_p^2\right)^{1/2}
 \lesssim2^{(k/4)(1-2/p)}\|P_kf\|_p.
\]
Consequently,
\begin{align*}
 \|A_uP_kf\|_p
 &\lesssim_\varepsilon2^{-2k/p+k\varepsilon/2}\|P_kf\|_p\\
 &=2^{-k/p}2^{-(1/p-\varepsilon/2)k}\|P_kf\|_p.
\end{align*}
Choose $\varepsilon=1/p$ and $\gamma_p=1/(2p)>0$, as desired.
Finally, if $2<p<6$, there are some little modifications:

\begin{align*}
 \|A_uP_kf\|_p
 &\approx2^{-k/2}
 \|E_{\operatorname{Cone}}(\widehat{P_kf})\|_p\\
 &\lesssim2^{-k/2}(2^{-k/2})^{-\varepsilon}
 \left(\sum_{\substack{I\subset[-1,1]\\\ell(I)=2^{-k/2}}}
 \|E_{\operatorname{Cone}}(\widehat{P_{k,I}f})\|_p^2
 \right)^{1/2}.
\end{align*}
Also,
\[
 \left(\sum_I\|P_{k,I}f\|_p^2\right)^{1/2}
 \lesssim2^{(k/4)(1-2/p)}\|P_kf\|_p.
\]
It follows that
\[
 \|A_uP_kf\|_p
 \lesssim
 2^{-k/p}2^{-(1/4-1/(2p)-\varepsilon/2)k}
 \|P_kf\|_p.
\]
Choose
\[
 \varepsilon=\frac14-\frac1{2p}>0
 \quad(\Longleftrightarrow p>2),
 \qquad
 \gamma_p=\frac12\left(\frac14-\frac1{2p}\right)>0.
\]
In summary, by Sobolev embedding and the local smoothing estimate,
we have
\begin{align*}
 \|M(P_kf)\|_p
 &=\left(\int_{\mathbb R^2}
 \left(\sup_{u\in[1,2]}|A_uP_kf|\right)^p\right)^{1/p}\\
 &\lesssim2^{k/p}
 \left(\int_{\mathbb R^2}
 \left(\int_1^2|A_uP_kf(x,y)|^p\,du\right)\,dx\,dy
 \right)^{1/p}.
\end{align*}

Continuing,
\[
 \|M(P_kf)\|_p
 \lesssim2^{k/p}2^{-k/p}2^{-\gamma_pk}\|P_kf\|_p
 \lesssim2^{-\gamma_pk}\|f\|_p.
\]
This completes the proof.

\paragraph{Remark.}
Most people believe that, for theorems related to positive operators,
for example Bourgain's circular maximal function $M$, a \emph{negative sign}
should not appear during the proof. In other words, the idea of using
interference cancellation is unnatural, and people tend to look for
a positive proof. Then, in 2001, Schlag gives a proof only using
incidence geometry and combinatorics. See~\cite{Schlag1998}.
\section{Wolff's circular maximal function}

Recall the circular Kakeya set: $E\subset\mathbb R^2$ is a Borel set satisfying: for every $r\in[1,2]$, there exists $x_r\in\mathbb R^2$ such that the circle $C(x_r,r)\subset E$, and $\mathcal L^2(E)=0$.

\paragraph{Minkowski dimension.} Let $N(\delta)$ be the smallest number of $\delta$-squares required to cover $E$. Then
\[
\dim_M(E)\approx\frac{\log N(\delta)}{\log(1/\delta)}.
\]
More precisely,
\[
\overline{\dim}_M(E):=\limsup_{\delta\to0}\frac{\log N(\delta)}{\log(1/\delta)},\qquad
\underline{\dim}_M(E):=\liminf_{\delta\to0}\frac{\log N(\delta)}{\log(1/\delta)}.
\]
The existence of circular Kakeya sets has been resolved by Besicovitch and Rado, '68. See~\cite{BesicovitchRado1968}. Its construction is very similar to Kakeya sets.

\noindent\textbf{Question.} Is the Minkowski dimension of any circular Kakeya set
always $2$?

\subsection{Maximal function}
Naturally, let
\[
 A_uf(x,y)=\int_{\mathbb R}f(x+ut,y+ut^2)\chi_0(t)\,dt,
 \qquad Mf(u)=\sup_{(x,y)\in\mathbb R^2}|A_uf(x,y)|,\quad u\in[1,2],
\]
and we ask whether
\begin{equation}\|Mf\|_{L_u^p([1,2])}\lesssim\|f\|_p?\tag{$*$}\end{equation}
Choose $f=1_{E_\delta}$, where $E_\delta=\mathcal N_\delta(E)$ and $E$ is a circular Kakeya set. Then by the definition of $E$,
\[
 Mf(u)\gtrsim1\quad\Longrightarrow\quad\|Mf\|_{L_u^p([1,2])}\gtrsim1.
\]
Therefore,
\[
 1\lesssim\|f\|_p=|E_\delta|^{1/p}\sim(N(\delta)\delta^2)^{1/p}
 \quad\Longrightarrow\quad N(\delta)\gtrsim\delta^{-2}
 \quad\Longrightarrow\quad\dim_M(E)\ge2.
\]
But in fact, we do not expect $N(\delta)\gtrsim\delta^{-2}$ to hold, since it essentially implies $|E|>0$, which contradicts our assumption. So we should prove $(*)$ with an $\varepsilon$-loss. To match the parameter $\delta$, we need the additional condition
\[
 \operatorname{supp}\widehat f\subset B_{\delta^{-1}}
 \quad\Longleftrightarrow\quad\text{modify the average operator by}
\]
\[
 A_{u,\delta}f(x,y)=\frac1\delta\int_{\mathbb R}\int_{-\delta}^{\delta}
 f(x+ut,y+ut^2+s)\chi_0(t)\,ds\,dt,
\]
and the corresponding circular maximal function
\[
 M^\delta f(u):=\sup_{(x,y)\in\mathbb R^2}|A_{u,\delta}f(x,y)|.
\]
\begin{center}\includegraphics[width=.35\linewidth]{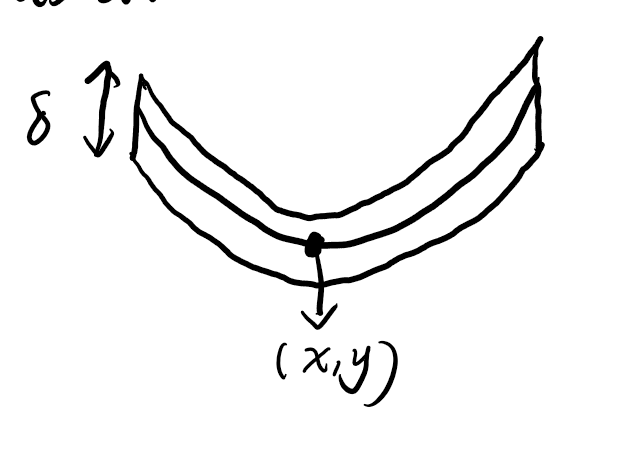}\end{center}

\noindent\textbf{Theorem 1 (Wolff, '97).} \cite{Wolff1997} For $p\ge3$,
\begin{equation}\|M^\delta f\|_{L_u^p([1,2])}\lesssim_\varepsilon\delta^{-\varepsilon}\|f\|_p.\tag{$**$}\end{equation}
\paragraph{Remarks.}
\begin{enumerate}
\item Using decoupling, we can only prove $(**)$ for $p\ge6$. Using the method from Guth--Wang--Zhang, we can prove $p\ge4$. See~\cite{GWZ2020}. But in fact $(**)$ holds whenever $p\ge3$ and is sharp for $p=3$.
\item In higher dimensions $n\ge3$, we can also ask:
\[
 \forall r\in[1,2],\quad\exists\ (n-1)\text{-dimensional sphere }S_r\subset E
 \quad\Longrightarrow\quad\dim_M(E)=n?
\]
Let
\[
 M^\delta f(u)=\sup_{x\in\mathbb R^n}\frac1{|S^\delta(x,u)|}\left|\int_{S^\delta(x,u)}f\right|.
\]
Do we have
\begin{equation}\|M^\delta f\|_{L_u^p([1,2])}\lesssim_\varepsilon\delta^{-\varepsilon}\|f\|_p?\tag{1}\end{equation}
This is proved by Kolasa and Wolff, '99, for $p\ge2$. See~\cite{KolasaWolff1999}. A reason why the higher dimensional case is simpler is that
\begin{center}\includegraphics[width=.95\linewidth]{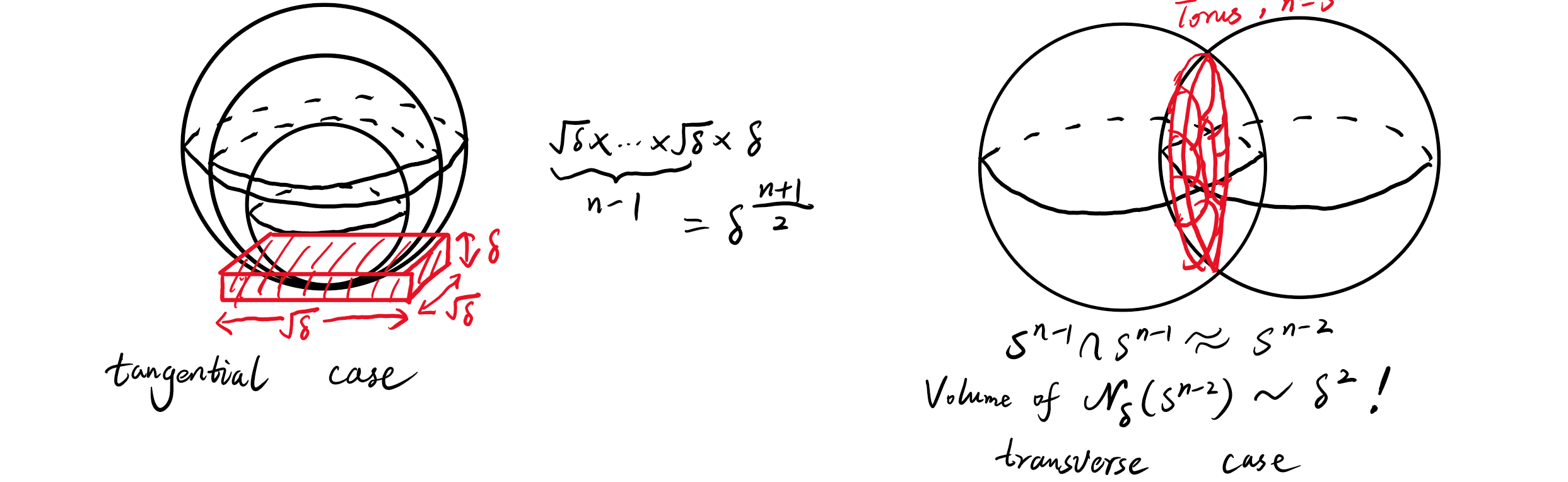}\end{center}
In the tangential case the volume is
\[
 \underbrace{\sqrt\delta\times\cdots\times\sqrt\delta}_{n-1}\times\delta
 =\delta^{(n+1)/2}.
\]
In the transverse case, $S^{n-1}\cap S^{n-1}\approx S^{n-2}$; the volume of $\mathcal N_\delta(S^{n-2})\sim\delta^2$ (a torus when $n=3$). Since $\delta^2\gg\delta^{(n+1)/2}$ as $n\to\infty$, the transverse case dominates the norm.
\item A ``correct'' question is, for example when $n=3$,
\begin{equation}
 \left\|\sup_{(x,y,z)\in\mathbb R^3}\left|\frac1\delta\int_{-\delta}^{\delta}\int_{\mathbb R}
 f(x+t,y+ut,z+ut^2+s)\chi_0(t)\,dt\,ds\right|\right\|_{L_u^p([1,2])}
 \lesssim_\varepsilon\delta^{-\varepsilon}\|f\|_p?\tag{2}
\end{equation}
Since in higher dimensions a sphere has codimension $1$, but a line/curve still has codimension $n-1$, while (1) is simple, (2) is
very hard.
\item By $(**)$, we get for any circular Kakeya set $E$,
\[
 N(\delta)\gtrsim\delta^{-2+\varepsilon}
 \quad\Longrightarrow\quad\dim_M(E)\ge2-\varepsilon\quad\text{for any }\varepsilon>0
 \quad\Longrightarrow\quad\dim_M(E)=2.
\]
\end{enumerate}
\subsection{Proof of Theorem 1 for \texorpdfstring{$p\ge6$}{p >= 6}}
This proof is very similar to Bourgain's circular maximal function. Take the Fourier transform of $A_{u,\delta}f$:
\begin{align*}
 \widehat{A_{u,\delta}f}(\xi,\eta)
 &=\frac1\delta\int_{\mathbb R}\int_{-\delta}^{\delta}\int_{\mathbb R^2}
 f(x+ut,y+ut^2+s)\chi_0(t)e^{-i(x\xi+y\eta)}\,dx\,dy\,ds\,dt\\
 &=\widehat f(\xi,\eta)\frac1\delta\int_{\mathbb R}\int_{-\delta}^{\delta}
 e^{iu(t\xi+t^2\eta)}e^{is\eta}\chi_0(t)\,ds\,dt\\
 &=\widehat f(\xi,\eta)\left(\int_{\mathbb R}e^{iu(t\xi+t^2\eta)}\chi_0(t)\,dt\right)
 \underbrace{\left(\frac1\delta\int_{-\delta}^{\delta}e^{is\eta}\,ds\right)}_{\displaystyle\approx1_{B(0,\delta^{-1})}(\eta)}.
\end{align*}
Without loss of generality assume $\operatorname{supp}\widehat f\subset B(0,\delta^{-1})$,
and it suffices to consider $A_uf$ rather than $A_{u,\delta}f$.

Littlewood--Paley decomposition:
\[
 f=P_{\le0}f+\sum_{k=1}^\infty P_kf.
\]
For the low frequency $P_{\le0}f$, the Fourier support lies in a fixed ball. By Bernstein's inequality,
\[
 \|P_{\le0}f\|_{\infty}\lesssim\|P_{\le0}f\|_p\lesssim\|f\|_p.
\]
Consequently, uniformly in $u\in[1,2]$,
\[
 M(P_{\le0}f)(u)=\sup_{(x,y)\in\mathbb R^2}|A_u(P_{\le0}f)(x,y)|
 \le\|\chi_0\|_1\|P_{\le0}f\|_{\infty}\lesssim\|f\|_p.
\]
Since $[1,2]$ has length $1$, it follows that
\[
 \|M(P_{\le0}f)\|_{L_u^p([1,2])}\lesssim\|f\|_p,
\]
and we are done.

For high frequency:
\paragraph{Sobolev embedding.}
Suppose $\operatorname{supp}\widehat f\subset\{(\xi,\eta):|\xi|^2+|\eta|^2\sim2^{2k}\}$. Then for any $p\ge1$,
\[
 \sup_{(x,y)\in\mathbb R^2}|f(x,y)|\lesssim2^{2k/p}\|f\|_p.
\]
Indeed,
\[
 \|f\|_p\gtrsim|B_{2^{-k}}|^{1/p}\|f\|_\infty=2^{-2k/p}\|f\|_\infty.
\]
\begin{center}\includegraphics[width=.35\linewidth]{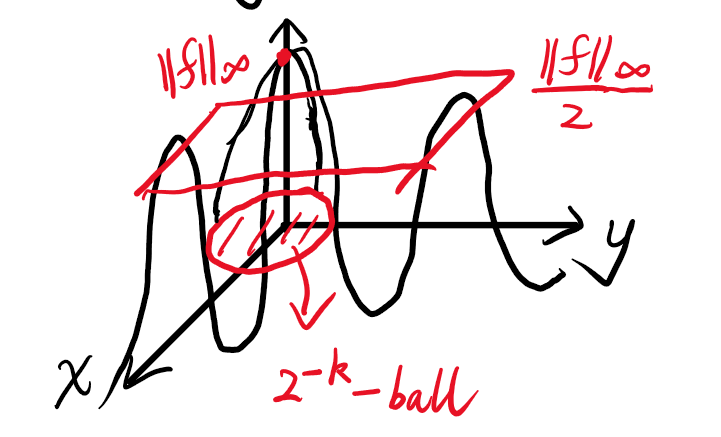}\end{center}
Consequently,
\[
 \sup_{(x,y)\in\mathbb R^2}|A_u(P_kf)(x,y)|\lesssim2^{2k/p}\|A_u(P_kf)\|_{L^p(\mathbb R^2)}.
\]
\paragraph{Local smoothing estimate.}
For any $p\ge6$,
\[
 \|A_u(P_kf)\|_{L^p(\mathbb R^2\times[1,2])}
 \lesssim_\varepsilon2^{-2k/p}2^{k\varepsilon/2}\|P_kf\|_p.
\]
Stationary phase principle:
\[
 A_u(P_kf)(x,y)\approx2^{-k/2}E_{\mathrm{cone}}(\widehat{P_kf})(x,y,u),
\]
where
\[
 \mathrm{cone}=\left\{\left(\xi,\eta,-\frac{\xi^2}{4\eta}\right):|\xi|\le10|\eta|,
 \ \sqrt{|\xi|^2+|\eta|^2}\sim2^k\right\}.
\]
Apply decoupling for the cone:
\begin{align*}
 \|A_u(P_kf)\|_p
 &\lesssim2^{-k/2}2^{\frac k2(\frac12-\frac3p+\varepsilon)}
 \left(\sum_{\substack{I\subset[-1,1]\\\ell(I)=2^{-k/2}}}\|A_u(P_{k,I}f)\|_p^2\right)^{1/2}\\
 &\approx2^{-k/4-3k/(2p)+k\varepsilon/2}
 \left(\sum_I\|P_{k,I}f\|_p^2\right)^{1/2}.
\end{align*}
Interpolation between $p=2$ and $p=\infty$ gives
\[
 \left(\sum_I\|P_{k,I}f\|_p^2\right)^{1/2}
 \lesssim2^{\frac k4(1-\frac2p)}\|P_kf\|_p.
\]
In summary,
\[
 \|A_u(P_kf)\|_p\lesssim2^{-2k/p}2^{k\varepsilon/2}\|P_kf\|_p,
\]
as desired. Finally, note that $\operatorname{supp}\widehat f\subset B(0,\delta^{-1})$, so $k=1,2,\ldots,O(\log\delta^{-1})$. Then
\[
 \sum_{k=1}^{O(\log\delta^{-1})}2^{k\varepsilon/2}\le C_\varepsilon\delta^{-\varepsilon}.
\]
\paragraph{Remark.} If $p<6$, using decoupling we only have
\[
 \|A_u(P_kf)\|_p\lesssim2^{-2k/p}2^{(\frac{3}{2p}-\frac14)k}2^{k\varepsilon/2}\|P_kf\|_p,
 \qquad\frac{3}{2p}-\frac14>0.
\]
This is not enough to complete the proof since
\[
 \sum_{k=1}^{O(\log\delta^{-1})}2^{\alpha_pk}2^{k\varepsilon/2}
 =C_\varepsilon\delta^{-\varepsilon}\delta^{-\alpha_p},
 \qquad\alpha_p=\frac{3}{2p}-\frac14>0,
\]
and $\delta^{-\alpha_p}\gg\delta^{-\varepsilon}$.

\subsection{Connection with incidence geometry}
For the circular incidence estimates behind the maximal theorem, see~\cite{Wolff1997}.
We consider the dual problem of $(**)$. Since $\operatorname{supp}\widehat f\subset B(0,\delta^{-1})$, the smallest scale on which we can distinguish $u$ is $\delta$. Choose
\[
 u_i\in[1+(i-1)\delta,1+i\delta],\qquad i=1,2,\ldots,\delta^{-1},
\]
and let $(x_i,y_i)$ be a point such that
\[
 |A_{u_i,\delta}f(x_i,y_i)|\ge\frac12\sup_{(x,y)\in\mathbb R^2}|A_{u_i,\delta}f(x,y)|.
\]
Then
\begin{align*}
 \|M^\delta f\|_{L_u^p([1,2])}
 &=\left(\sum_i\int_{1+(i-1)\delta}^{1+i\delta}|M^\delta f(u)|^p\,du\right)^{1/p}\\
 &\approx\left(\sum_i|A_{u_i,\delta}f(x_i,y_i)|^p\delta\right)^{1/p}\\
 &=\delta^{1/p}\sup\left\{\left|\sum_i a_iA_{u_i,\delta}f(x_i,y_i)\right|:\|(a_i)\|_{\ell^{p'}}\le1\right\}
 \quad\text{(duality)}.
\end{align*}
Here
\[
 A_{u_i,\delta}f(x_i,y_i)=\frac1\delta\int_{\mathbb R}\int_{-\delta}^{\delta}
 f(x_i+u_it,y_i+u_it^2+s)\chi_0(t)\,ds\,dt
 =\frac1\delta\langle f,1_{P_i}\rangle,
\]
where
\[
 P_i=\{(x_i+u_it,y_i+u_it^2+s)\in\mathbb R^2:t\in[-1,1],\ |s|\le\delta\}.
\]
Continuing,
\begin{align*}
 \|M^\delta f\|_{L_u^p([1,2])}
 &=\delta^{1/p-1}\sup\left\{\left|\sum_i a_i\langle f,1_{P_i}\rangle\right|:\|(a_i)\|_{\ell^{p'}}\le1\right\}\\
 &=\delta^{1/p-1}\sup\left\{\left|\left\langle f,\sum_i a_i1_{P_i}\right\rangle\right|:\|(a_i)\|_{\ell^{p'}}\le1\right\}\\
 &\le\delta^{-1/p'}\|f\|_p\sup\left\{\left\|\sum_i a_i1_{P_i}\right\|_{p'}:\|(a_i)\|_{\ell^{p'}}\le1\right\}.
\end{align*}
Therefore, Theorem 1 is equivalent to the following theorem.

\noindent\textbf{Theorem 2.} For $p\le3/2$,
\[
 \left\|\sum_{i=1}^{\delta^{-1}}a_i1_{P_i}\right\|_{L^p(\mathbb R^2)}
 \lesssim_\varepsilon\delta^{-\varepsilon}\delta^{1/p}\|(a_i)\|_{\ell^p}.
\]
Moreover, by dyadic pigeonholing, we have a more concise version of Theorem 2.

\noindent\textbf{Theorem 3.} Let $\mathbb P$ be a set of $\delta$-parabolas
\[
 P^\delta_{u,(x,y)}:=\{(x+ut,y+ut^2+s):t\in[-1,1],\ |s|\le\delta\},
\]
where $u\in[1,2]$ and $(x,y)\in\mathbb R^2$. Then for any finite subset $\mathbb P'\subset\mathbb P$ with $\delta$-separated ``directions $u$'',
\[
 \left\|\sum_{P^\delta\in\mathbb P'}1_{P^\delta}\right\|_{L^p(\mathbb R^2)}
 \lesssim_\varepsilon\delta^{-\varepsilon}\delta^{1/p}(\#\mathbb P')^{1/p}
\]
holds for all $p\le3/2$.

\paragraph{Proof of Theorem 2 $\Longleftrightarrow$ Theorem 3.}
``Theorem 2 $\Longrightarrow$ Theorem 3'' is trivial. Now we prove ``Theorem 3 $\Longrightarrow$ Theorem 2''. By normalizing $\{a_i\}$, suppose $\|(a_i)\|_{\ell^p}=1$. Let
\[
 I_0=\{i:|a_i|\le\delta^{100}\},\qquad
 I_k=\{i:2^{-k}\le|a_i|\le2^{-(k-1)}\},\quad k=1,2,\ldots,O(\log\delta^{-1}).
\]
Note that $\#I_0,\#I_k\le\delta^{-1}$, so
\[
 \left\|\sum_i a_i1_{P_i}\right\|_p
 \le\left\|\sum_{i\in I_0}a_i1_{P_i}\right\|_p
 +\sum_{k=1}^{O(\log\delta^{-1})}\left\|\sum_{i\in I_k}a_i1_{P_i}\right\|_p.
\]
For the first term,
\begin{align*}
 \left\|\sum_{i\in I_0}a_i1_{P_i}\right\|_p
 &\le\delta^{100}\left(\int\left|\sum_{i\in I_0}1_{P_i}\right|^p\right)^{1/p}\\
 &\lesssim_\varepsilon\delta^{100}\delta^{-\varepsilon}\delta^{1/p}(\#I_0)^{1/p}
 \le\delta^{99}\le\delta^{1/p}.
\end{align*}
For each $I_k$ in the second term,
\begin{align*}
 \left\|\sum_{i\in I_k}a_i1_{P_i}\right\|_p
 &\le2^{-(k-1)}\left\|\sum_{i\in I_k}1_{P_i}\right\|_p\\
 &\lesssim_\varepsilon2^{-(k-1)}\delta^{-\varepsilon}\delta^{1/p}(\#I_k)^{1/p}\\
 &\le2\delta^{-\varepsilon}\delta^{1/p}\left(\sum_{i\in I_k}|a_i|^p\right)^{1/p}
 \le2\delta^{-\varepsilon}\delta^{1/p}.
\end{align*}
Therefore,
\[
 \left\|\sum_i a_i1_{P_i}\right\|_p
 \lesssim_\varepsilon\delta^{1/p}+\delta^{-\varepsilon}\delta^{1/p}\log(\delta^{-1})
 \lesssim_\varepsilon\delta^{-\varepsilon}\delta^{1/p},
\]
as desired.\hfill$\square$

\paragraph{Remarks.}
\begin{enumerate}
\item Note that $\delta(\#\mathbb P')\sim\sum_{P^\delta\in\mathbb P'}|P^\delta|$. By Theorem 3 and H\"older's inequality,
\begin{align*}
 \sum_{P^\delta\in\mathbb P'}|P^\delta|
 &=\int\sum1_{P^\delta}
 \le\left|\bigcup P^\delta\right|^{1-1/p}\left\|\sum1_{P^\delta}\right\|_p\\
 &\lesssim_\varepsilon\delta^{-\varepsilon}
 \left(\sum|P^\delta|\right)^{1/p}\left|\bigcup P^\delta\right|^{1-1/p}.
\end{align*}
Thus
\[
 \sum_{P^\delta\in\mathbb P'}|P^\delta|
 \lesssim_\varepsilon\delta^{-\varepsilon}\left|\bigcup_{P^\delta\in\mathbb P'}P^\delta\right|,
\]
which implies the $\delta$-parabolas in $\mathbb P'$ are essentially pairwise disjoint.
\item The range $p\le3/2$ is sharp. Consider $\delta^{-1}$ $\delta$-parabolas that are all tangent at one point. Then
\[
 \sum1_{P^\delta}(x)=\#P^\delta=\delta^{-1}\qquad\text{for any }x\in R.
\]
\begin{center}\includegraphics[width=.3\linewidth]{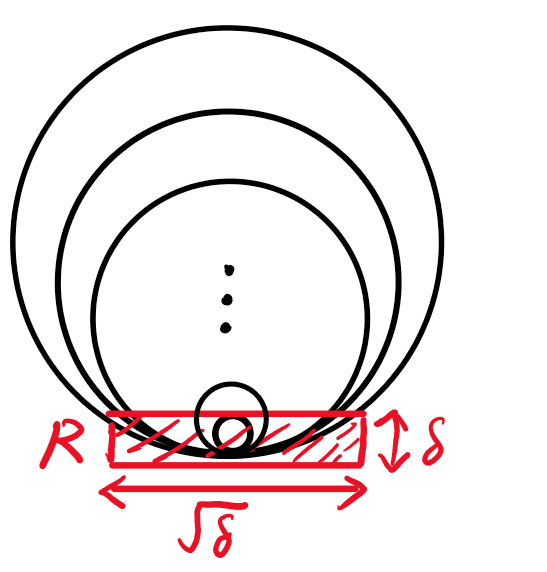}\end{center}
Hence,
\[
 1=\delta^{1/p}(\#P^\delta)^{1/p}\gtrsim_\varepsilon\delta^\varepsilon\left\|\sum1_{P^\delta}\right\|_p
 \ge\delta^\varepsilon\delta^{-1}|R|^{1/p}=\delta^{\varepsilon-1+3/(2p)}.
\]
Thus $(1-\varepsilon)p\le3/2$ for any $\varepsilon>0$, hence $p\le3/2$.
\end{enumerate}

\section{The Kakeya problem on \texorpdfstring{$\mathbb R^2$}{R2}}
The maximal-function estimate below is the planar C\'ordoba estimate~\cite{Cordoba1977}.

Recall a compact set $E\subset\mathbb R^2$ is a Kakeya set if for any $e\in S^1$, there is a line segment $\{x_e+te:t\in[0,1]\}\subset E$, and $\mathcal L^2(E)=0$.

\noindent\textbf{Question.} Do any Kakeya sets have Hausdorff dimension $2$?

\noindent\textbf{Theorem 1.} Let $\mathbb T$ be a set of $\delta$-tubes with $\delta$-separated directions. Then
\[
 \left\|\sum_{T\in\mathbb T}1_T\right\|_{L^2(\mathbb R^2)}^2
 \lesssim\log\delta^{-1}\left\|\sum_{T\in\mathbb T}1_T\right\|_1
 \quad(=\log\delta^{-1}\cdot\delta\cdot\#\mathbb T).
\]
\paragraph{Proof.} Direct computation shows
\[
 \left\|\sum_{T\in\mathbb T}1_T\right\|_2^2
 =\int\sum_{T_1,T_2\in\mathbb T}1_{T_1}1_{T_2}
 =\sum_{T_1,T_2\in\mathbb T}|T_1\cap T_2|.
\]
Since if $T_1\not\approx T_2$, $T_1\cap T_2\ne\varnothing$ and $\delta\le\angle(T_1,T_2)\lesssim1$, we can do dyadic pigeonholing as follows:
\[
 \sum_{T_1,T_2\in\mathbb T}|T_1\cap T_2|
 =\sum_{T_1\approx T_2}|T_1|+
 \sum_{j=1}^{O(\log\delta^{-1})}\sum_{\substack{T_1,T_2\\\angle(T_1,T_2)\sim2^{-j}}}|T_1\cap T_2|.
\]
If $\angle(T_1,T_2)\sim2^{-j}$, then $|T_1\cap T_2|\sim\delta^2/2^{-j}$.
\begin{center}\includegraphics[width=.25\linewidth]{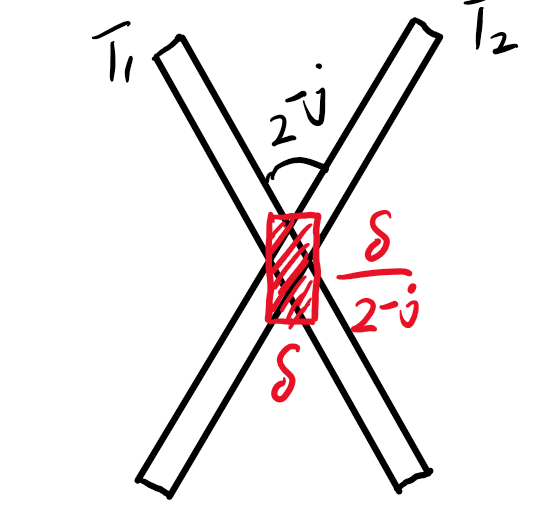}\end{center}
Another perspective: two extreme scenarios,
\[
 \angle(T_1,T_2)\sim1\quad\Longrightarrow\quad|T_1\cap T_2|\approx|\delta\text{-square}|=\delta^2,
\]
\[
 \angle(T_1,T_2)\sim\delta\quad\Longrightarrow\quad|T_1\cap T_2|\approx|\delta\text{-tube}|=\delta.
\]
The situation of intermediate values should correspond to interpolation of the two extreme scenarios.

Finally, fix $T_1$. By the $\delta$-separated property,
\[
 \#\{T_2\in\mathbb T:T_1\cap T_2\ne\varnothing,\ \angle(T_1,T_2)\sim2^{-j}\}
 \lesssim\frac{2^{-j}}\delta,
\]
and the implicit constant is independent of $T_1$.
\begin{center}\includegraphics[width=.25\linewidth]{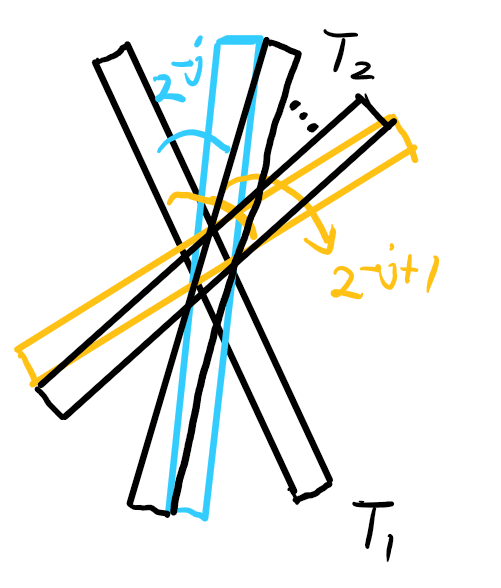}\end{center}
Hence,
\[
 \left\|\sum1_T\right\|_2^2
 \lesssim\delta\cdot\#\mathbb T+
 \sum_{j=1}^{O(\log\delta^{-1})}2^j\delta^2\cdot\frac{2^{-j}}\delta\cdot\#\mathbb T
\]
\[
 \lesssim\log\delta^{-1}\cdot\delta\cdot\#\mathbb T,
\]
as desired.\hfill$\square$

By a duality argument as we did in Section 7, we obtain the boundedness of the Kakeya maximal function: for $u\in[0,1]\simeq S^1$,
\[
 K_\delta f(u):=\sup_{(x,y)\in\mathbb R^2}\frac1\delta\int_{\mathbb R}\int_{-\delta}^{\delta}
 f(x+t,y+ut+s)\chi_0(t)\,ds\,dt.
\]
\begin{center}\includegraphics[width=.2\linewidth]{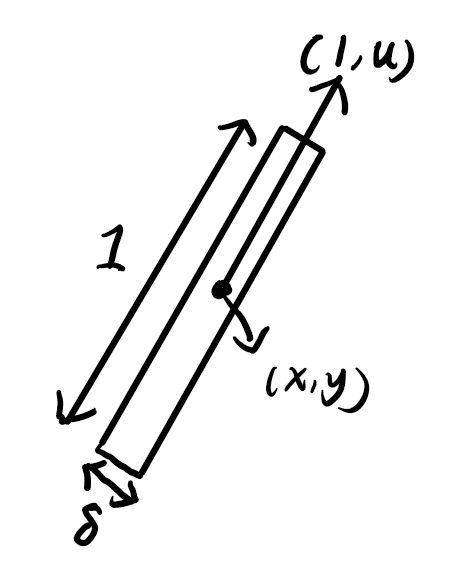}\end{center}
\noindent\textbf{Theorem 2.} For $0<\delta\le1$,
\[
 \|K_\delta f\|_{L_u^2(S^1)}\lesssim\sqrt{\log(\delta^{-1})}\|f\|_{L^2(\mathbb R^2)}.
\]
Test the characteristic function on a Kakeya set, $1_{E_\delta}$; we prove that:

\noindent\textbf{Theorem 3.} Any Kakeya set in the plane has dimension $2$. See~\cite{Davies1971,Mattila2015}.

\paragraph{Remark.} The loss $\sqrt{\log(\delta^{-1})}$ is sharp. To see this, consider the first example in Section 4. For $H$ steps:
\begin{center}\includegraphics[width=.95\linewidth]{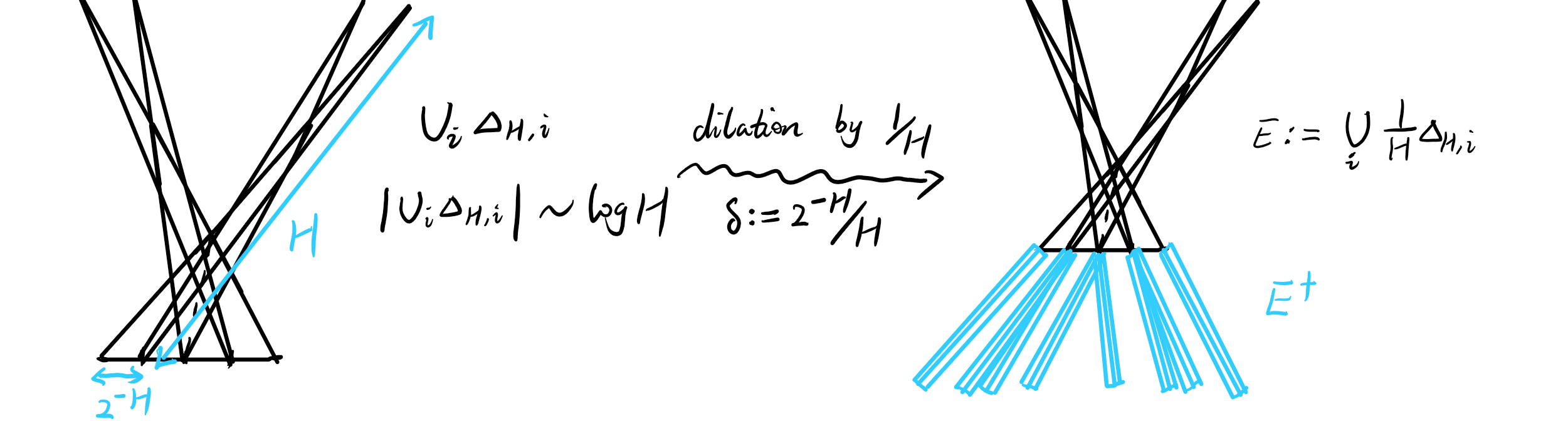}\end{center}
\[
 \left|\bigcup_i\Delta_{H,i}\right|\sim\log H.
\]
Dilate by $1/H$ and set
\[
 \delta:=\frac{2^{-H}}H,\qquad E:=\bigcup_i\frac1H\Delta_{H,i}.
\]
Extend the two long sides of each thin triangle $\Delta_{H,i}$ backwards by segments of length $1/2$, and then thicken these segments by $\delta$, that is, take their $\delta$-neighbourhoods. Denote the resulting union by $E^+$.

By construction, $|E|=\log H/H^2$. $E^+$ is the union of $\delta$-tubes with $\delta$-separated directions, and the number of directions is $\sim2^H$. Take $f=1_E$. Observe that for any of these directions $(1,u)$, there exists $(x,y)\in E^+$ such that the $\delta$-tube
\[
 T^\delta_{u,(x,y)}:=\{(x+t,y+ut+s):|t|\le1,\ |s|\le\delta\}
\]
satisfies
\[
 |T^\delta_{u,(x,y)}\cap E|\ge\frac12|T^\delta_{u,(x,y)}|.
\]
Thus $K^\delta f(u)\gtrsim1$.
Therefore,
\[
 \|K^\delta f\|_{L_u^2([1,2])}^2\gtrsim\delta\cdot\#\text{directions}=\frac1H.
\]
Every direction contributes a region with angle $\delta$ in $S^1$.
In addition,
\[
 \|f\|_2=|E|^{1/2}\sim\frac{(\log H)^{1/2}}H.
\]
It follows that
\[
 \frac{\|K^\delta f\|_2}{\|f\|_2}\gtrsim\left(\frac H{\log H}\right)^{1/2}
 \sim\left(\frac{\log\delta^{-1}}{\log\log\delta^{-1}}\right)^{1/2}.
\]
To obtain the sharp bound $(\log\delta^{-1})^{1/2}$, we need a more detailed construction for Kakeya sets, due to Schoenberg, '62 and Keich, '99. See~\cite{Schoenberg1962,Keich1999}.

\section{The Nikodym problem on \texorpdfstring{$\mathbb R^2$}{R2}}

A set $E\subset[0,1]^2$ is a Nikodym set if $\mathcal L^2(E)=0$ and for any $x\in[0,1]^2\setminus E$, there exists a line segment $l_x$ through $x$ such that $l_x\setminus\{x\}\subset E$.

The first construction of Nikodym sets is due to Nikodym, '27. See~\cite{Nikodym1927}.

\noindent\textbf{Question.} Does any Nikodym set in the plane have dimension $2$?

\subsection{Maximal operator}
For the planar maximal estimate, see~\cite{Cordoba1977,Keich1999}.
Suppose $0<\delta<1$, $x\in\mathbb R^2$, $e\in S^1$. Define the $\delta$-tube
\[
 T_e^\delta(x)=\mathcal N_\delta\bigl(\{x+te\in\mathbb R^2:t\in[0,1]\}\bigr).
\]
Define the Nikodym maximal function:
\[
 \mathcal N_\delta f(x)=\sup_{e\in S^1}\frac1{|T_e^\delta(x)|}\int_{T_e^\delta(x)}|f|.
\]
\begin{center}\includegraphics[width=.17\linewidth]{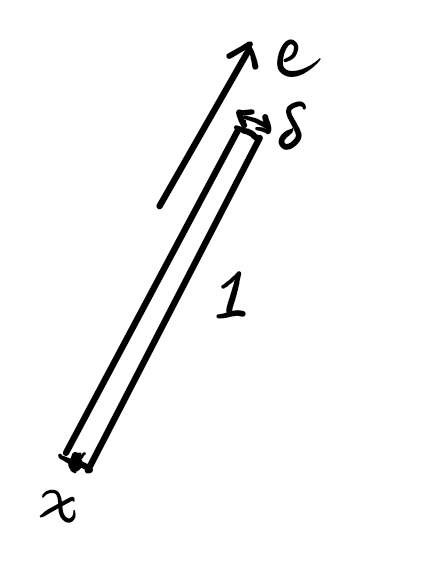}\end{center}
\noindent\textbf{Theorem.} For $p\ge2$,
\[
 \|\mathcal N_\delta f\|_{L^p(\mathbb R^2)}\lesssim_p(\log\delta^{-1})^{1/p}\|f\|_{L^p(\mathbb R^2)}.
\]
\subsection{Examples}
In the following two examples, we aim to show the loss $(\log\delta^{-1})^{1/p}$ is sharp.

\paragraph{1.} Consider the first example in Section 4. For $H$ steps,
\begin{center}\includegraphics[width=.95\linewidth]{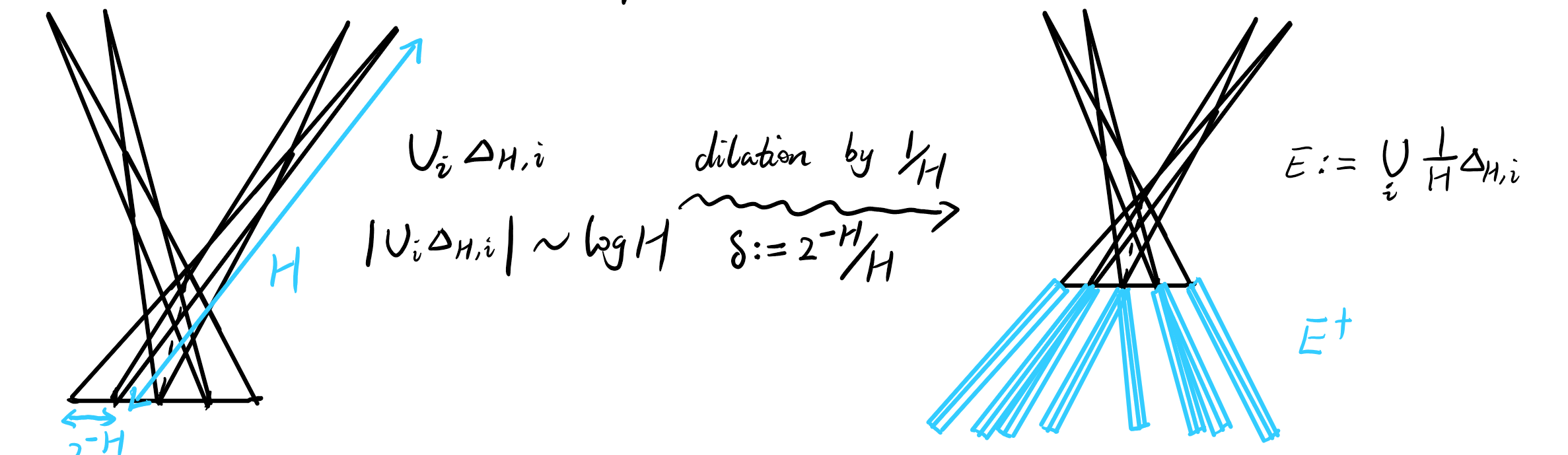}\end{center}
\[
 \left|\bigcup_i\Delta_{H,i}\right|\sim\log H.
\]
Dilate by $1/H$ and set
\[
 \delta:=2^{-H}/H,\qquad E:=\bigcup_i\frac1H\Delta_{H,i}.
\]
Extend the two long sides of each thin triangle $\Delta_{H,i}$ backwards by segments of length $1/2$,
and then thicken these segments by $\delta$, that is, take their $\delta$-neighbourhoods. Denote the resulting union by $E^+$.

By construction, $|E|=\log H/H^2$. $E^+$ is the union of $\delta$-tubes with $\delta$-separated directions, and the number of these tubes is $\sim2^H$. Take $f=1_E$. Observe that for any point $x\in E^+$, there exists a tube $T_e^\delta(x)$ such that
\[
 |T_e^\delta(x)\cap E|\approx\frac12|T_e^\delta(x)|,
\]
which implies that
\[
 \mathcal N_\delta f(x)\gtrsim1\quad\text{on }E^+.
\]
Thus
\[
 \|\mathcal N_\delta f\|_p\gtrsim|E^+|^{1/p}=(\delta\cdot\#\text{tubes})^{1/p}=\left(\frac1H\right)^{1/p}.
\]
Moreover,
\[
 \|f\|_p=|E|^{1/p}=\left(\frac{\log H}{H^2}\right)^{1/p},
\]
so
\[
 \frac{\|\mathcal N_\delta f\|_p}{\|f\|_p}\gtrsim\left(\frac H{\log H}\right)^{1/p}
 \sim\left(\frac{\log\delta^{-1}}{\log\log\delta^{-1}}\right)^{1/p}.
\]
\paragraph{2. (Schoenberg and Keich).} \cite{Schoenberg1962,Keich1999}
For every positive integer $k$, there exists a collection of rectangles
\[
 R_1,R_2,\ldots,R_{2^k}\subset\mathbb R^2
\]
of dimensions $1\times2^{-k}$ satisfying:
\begin{enumerate}[label=(\roman*)]
\item The two long sides of $R_i$ point in the direction $(1,i\cdot2^{-k})$ for every $i$.
\item $\displaystyle\mathcal L^2\left(\bigcup_{i=1}^{2^k}R_i\right)\le\frac{1000}k$.
\item If we denote by $R_i'$ the translation of $R_i$ in the direction $(1,i\cdot2^{-k})$ by distance $10$, then
\[
 R_i'\cap R_j'=\varnothing\qquad\text{whenever }i\ne j.
\]
\end{enumerate}
Compared with Example 1, let $\delta=2^{-k}$; then we obtain $\delta^{-1}$ $\delta$-tubes with $\delta$-separated directions. It is essentially a Kakeya set, also a Nikodym set. Hence we can expect it to produce the smallest loss.

\subsection{Fubini's nightmare}
For Katok\textquotesingle s example as presented by Milnor, see~\cite{Milnor1997}.
Anatole Katok has proved the following:

There exists a measurable set $E$ in the unit square $[0,1]^2$ with
$\mathcal L^2(E)=1$, together with a family of disjoint smooth real analytic curves $\Gamma_\beta$ which fill out this square, so that each curve $\Gamma_\beta$ intersects the set $E$ in at most a single point.
\begin{center}\includegraphics[width=.3\linewidth]{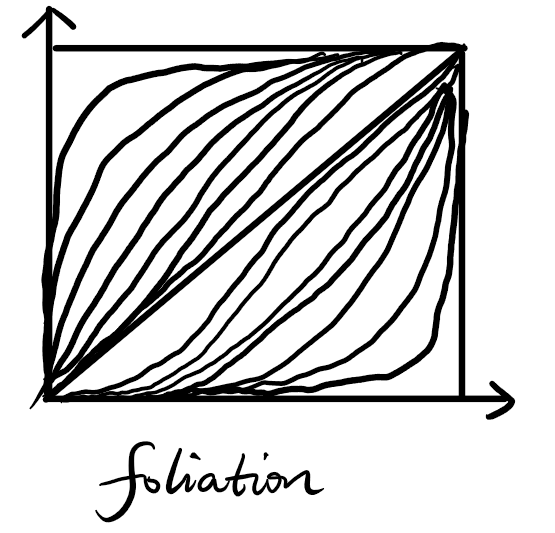}\end{center}
In other words, we can construct a set of full $2$-dimensional Lebesgue measure by selecting at most one point from each $\Gamma_\beta$. For this foliation $\{\Gamma_\beta\}$ and set $E$, we observe a Fubini's nightmare.

Fubini's nightmare is a seeming violation of Fubini's theorem. But there is no real contradiction to Fubini's theorem because despite smoothness of the fibers $\Gamma_\beta$, the foliation is not absolutely continuous, and neither are the conditional measures on fibers.

\subsection{Proof of the theorem}
Without loss of generality, we restrict the direction
\[
 e\in S^1\cap\left([0,1]\times\left[0,\frac\pi{1000}\right]\right).
\]
By Fubini's theorem, it suffices to prove
\[
 \int|\mathcal N_\delta f(x_1,x_2)|^2\,dx_2
 \lesssim(\log\delta^{-1})\|f\|_{L^2(I_{x_1}\times\mathbb R)}^2
\]
holds uniformly in $x_1$, where $I_{x_1}=[x_1,x_1+100]$.
\begin{center}\includegraphics[width=.25\linewidth]{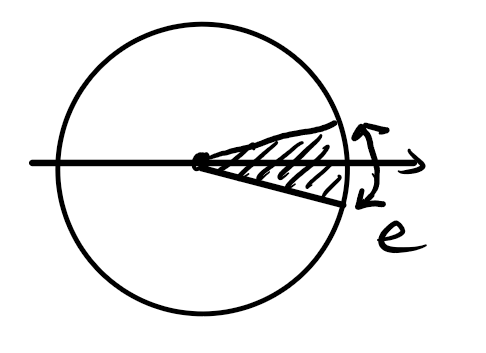}\end{center}
\begin{align*}
 \widehat{\mathcal N_\delta f}(\xi,\eta)
 &=\int_{\mathbb R^2}\frac1\delta\int_{\mathbb R}\int_{-\delta}^{\delta}
 f(x+te_1,y+te_2+s)\chi_0(t)\,ds\,dt\,e^{-i(x\xi+y\eta)}\,dx\,dy\\
 &=\widehat f(\xi,\eta)\int_{\mathbb R}e^{i(\xi e_1+\eta e_2)t}\chi_0(t)\,dt
 \underbrace{\frac1\delta\int_{-\delta}^{\delta}e^{is\eta}\,ds}_{\displaystyle=1_{B(0,\delta^{-1})}(\eta)\ \Longrightarrow\ |\eta|\le\delta^{-1}}.
\end{align*}
Since
\[
 \left|\frac\partial{\partial t}[(\xi e_1+\eta e_2)t]\right|
 =|\xi e_1+\eta e_2|\approx|\xi|\ge100
 \qquad\text{whenever }|\xi|\ge10|\eta|,
\]
Van der Corput's lemma shows that $\widehat f$ is essentially supported in $B(0,\delta^{-1})$.

Without loss of generality assume $x_1=0$. Note that $\widehat{\mathcal N_\delta f}$ is essentially supported in a ball of radius $\delta^{-1}$. We know $|\mathcal N_\delta f(0,x_2)|\approx\text{constant}$ when $x_2$ is localized at scale $\delta$. Let $I_k=[k\delta,(k+1)\delta]$ and let $R_k$ be the corresponding tube such that
\[
 |\mathcal N_\delta f(0,x_2)|\approx\frac1{|R_k|}\int_{R_k}|f|\qquad\forall x_2\in I_k.
\]
\begin{center}\includegraphics[width=.3\linewidth]{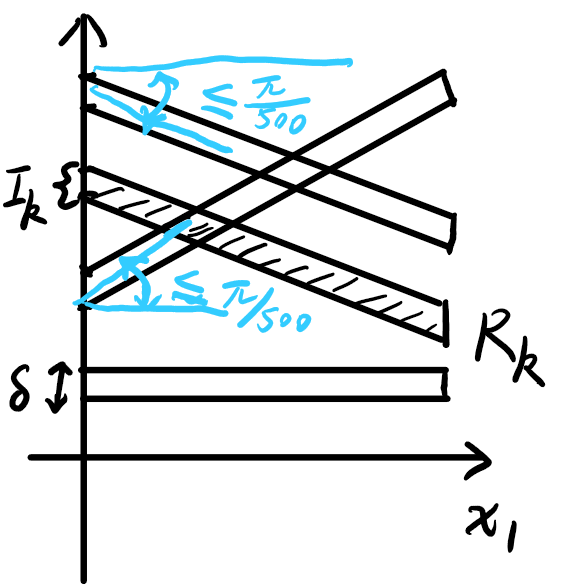}\end{center}
Using the dual argument as before,
\begin{align*}
 \|\mathcal N_\delta f(0,\cdot)\|_2
 &=\left(\sum_k\delta\left(\frac1{|R_k|}\int_{R_k}|f|\right)^2\right)^{1/2}\\
 &=\delta^{1/2}\left(\sum_k\left(\frac1{|R_k|}\langle f,1_{R_k}\rangle\right)^2\right)^{1/2}\\
 &=\delta^{1/2}\sup\left\{\left|\sum_k a_k\frac1{|R_k|}\langle f,1_{R_k}\rangle\right|:\|(a_k)\|_2\le1\right\}\\
 &=\delta^{1/2}\sup\left|\left\langle f,\sum_k a_k\frac{1_{R_k}}{|R_k|}\right\rangle\right|\\
 &\le\delta^{1/2}\|f\|_{L^2([0,100]\times\mathbb R)}
 \sup\left\|\sum_k a_k\frac{1_{R_k}}{|R_k|}\right\|_2.
\end{align*}
Note that $\bigcup R_k\subset[0,100]\times\mathbb R$. Then by dyadic pigeonholing, it suffices to prove
\[
 \left\|\sum_k|R_k|^{-1}1_{R_k}\right\|_{L^2(\mathbb R^2)}
 \lesssim(\log\delta^{-1})^{1/2}\delta^{-1/2}(\#k)^{1/2}.
\]
Direct computation:
\begin{align*}
 \left\|\sum_k|R_k|^{-1}1_{R_k}\right\|_2^2
 &=\int\sum_{k,l}|R_k|^{-1}|R_l|^{-1}1_{R_k\cap R_l}\\
 &=\delta^{-2}\sum_{k,l}|R_k\cap R_l|\\
 &=\delta^{-2}\left(\sum_{k=l}|R_k|
 +\sum_{\ell=-\log\delta^{-1}}^\infty\sum_{\operatorname{dist}(I_k,I_l)\sim2^\ell}|R_k\cap R_l|\right).
\end{align*}
Recall $e\in S^1\cap([0,1]\times[0,\pi/1000])$, so $\angle(R_k,R_l)\le\pi/500$. Therefore $R_k\cap R_l=\varnothing$ whenever $\operatorname{dist}(I_k,I_l)\ge100$, and
\[
 |R_k\cap R_l|\sim\frac{\delta^2}{\angle(R_k,R_l)}
 \lesssim\frac{\delta^2}{\operatorname{dist}(I_k,I_l)}.
\]
(The smallest angle between $R_k$ and $R_l$ is $\sim\operatorname{dist}(I_k,I_l)$.)
\begin{center}\includegraphics[width=.3\linewidth]{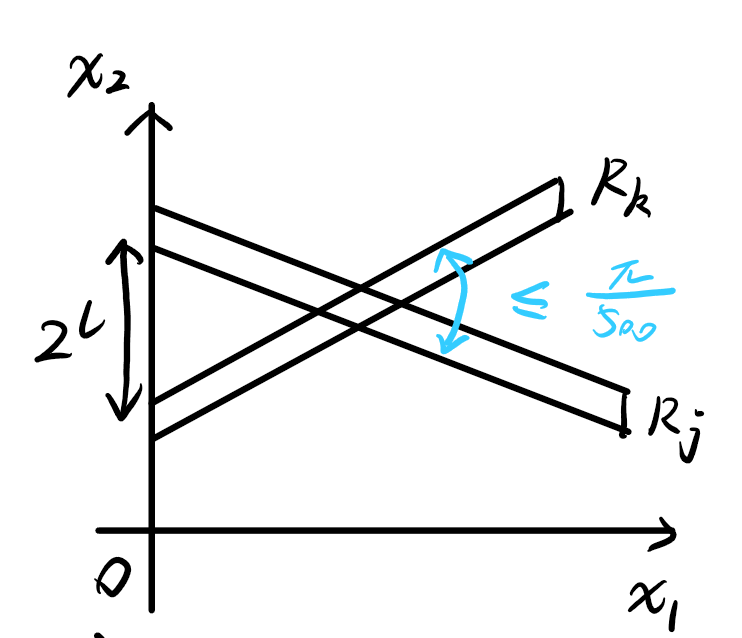}\quad\includegraphics[width=.3\linewidth]{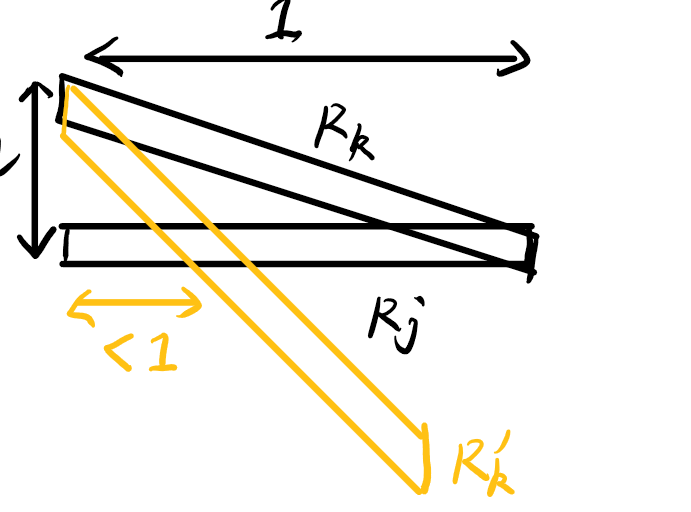}\end{center}
So
\begin{align*}
 \left\|\sum_k|R_k|^{-1}1_{R_k}\right\|_2^2
 &\lesssim\delta^{-2}\left(\delta(\#k)+\sum_{I_k}\sum_{\ell=-\log\delta^{-1}}^{100}
 \sum_{\substack{I_l\\\operatorname{dist}(I_k,I_l)\sim2^\ell}}2^{-\ell}\delta^2\right)\\
 &\lesssim\delta^{-2}\left(\delta(\#k)+(\#k)O(\log\delta^{-1})\frac{2^\ell}\delta\cdot2^{-\ell}\delta^2\right)\\
 &\lesssim\log\delta^{-1}\cdot\delta^{-1}(\#k),
\end{align*}
as desired. Finally, we have the trivial bound when $p=\infty$:
\[
 \|\mathcal N_\delta f\|_{L^\infty(\mathbb R^2)}\lesssim\|f\|_{L^\infty(\mathbb R^2)}.
\]
Interpolating between $2$ and $\infty$, we obtain
\[
 \|\mathcal N_\delta f\|_{L^p(\mathbb R^2)}\lesssim(\log\delta^{-1})^{1/p}\|f\|_{L^p(\mathbb R^2)}.
\]
\subsection{Comments from geometric measure theory}
For the projection-theoretic background to point--line duality, see~\cite{Mattila2015}.
We say the Nikodym set looks like it contains many more lines. By
standard point--line duality, we can view a nonvertical line as a point $(a,b)$ in $\mathbb R^2$, where $a$ is the slope and $b$ is the intersection with the $y$-axis. So we obtain a one-to-one correspondence:
\[
 \mathcal L\xrightarrow{D^{-1}}A,\qquad A\xrightarrow{D}\mathcal L.
\]
\begin{center}\includegraphics[width=.7\linewidth]{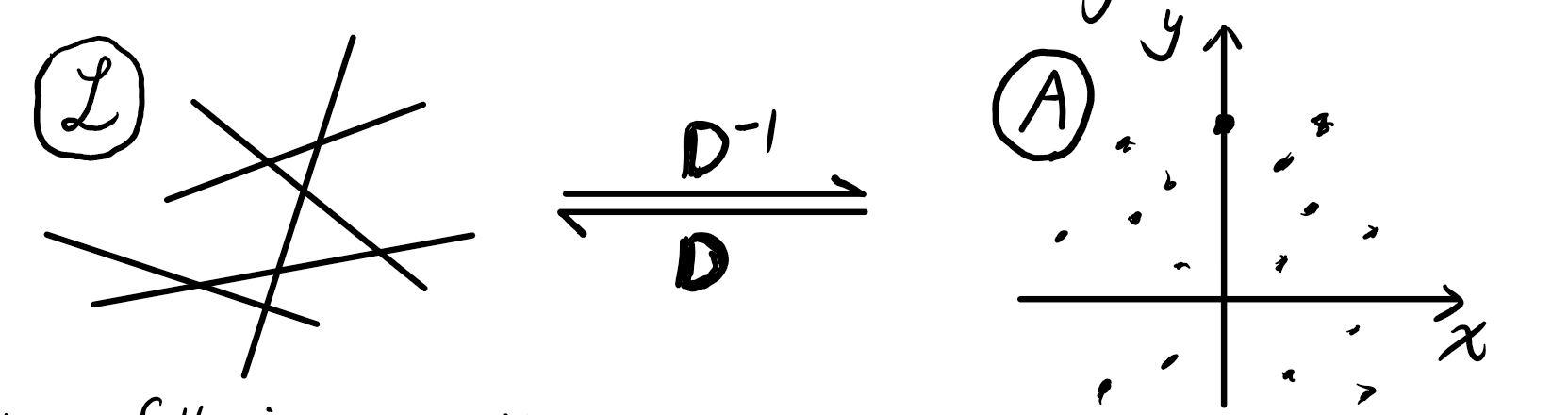}\end{center}
We have known the following result:
\[
 \boxed{\text{If }\dim_H(A)>1,\text{ then }m\left(\bigcup_{L\in\mathcal L}L\right)>0.}
\]
\paragraph{For Kakeya sets.}
We imagine the map $a\mapsto b$ is continuous. Therefore, $\dim A=1$.
\begin{center}\includegraphics[width=.3\linewidth]{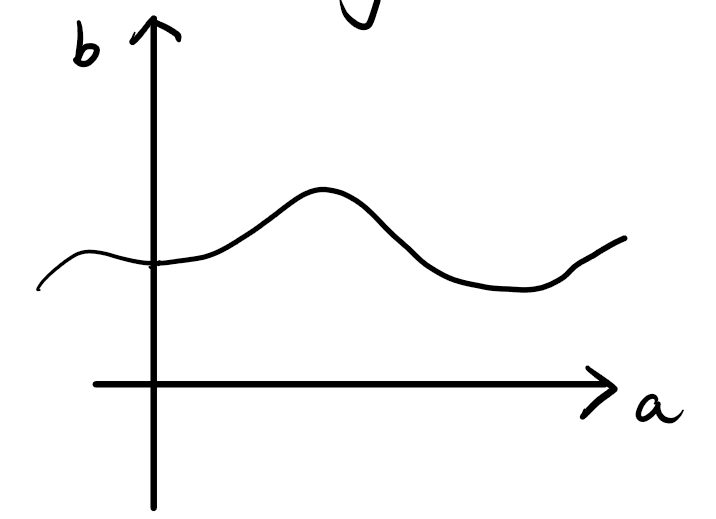}\end{center}
\paragraph{For Nikodym sets.}
It seems that $\dim A=2$ since for a.e. $x\in\mathbb R^2$, there is $e\in S^1$, giving $\mathbb R^2\ni x\mapsto e\in S^1$. But by the existence of Nikodym sets, we know $\dim A=1$. So Nikodym sets do not contain too many lines. And this fact is acceptable because there exists a bijection $[0,1]^2\to[0,1]$, so \emph{the existence of a bijection between two sets does not mean that their dimensions are equal}.

\subsection{Geometry of the cinematic curvature condition}
For precise formulations of cinematic curvature and the corresponding maximal estimates, see~\cite{MSS1992,Schlag1997}.
We know if the curves in the definition of Nikodym sets have cinematic curvature $\ne0$, then it must have positive Lebesgue measure, i.e., the corresponding Nikodym sets do not exist. Here we do not give a formal definition, but rather look at what geometric information it contains.

\paragraph{(Line) Nikodym sets: exist.}
In constructions of a Nikodym set, for a fixed line $l_x$ through $x$, when choosing other lines through $x'$ close enough
to $l_x$, our principle is to be incident to most points on $l_x$ to make the overlapping large enough, i.e., to make the lines sufficiently close to $l_x$. In particular, for $x$ and $l_x$, for any $x'$, there exists a line through $x'$ such that it coincides with $l_x$.
\begin{center}\includegraphics[width=.23\linewidth]{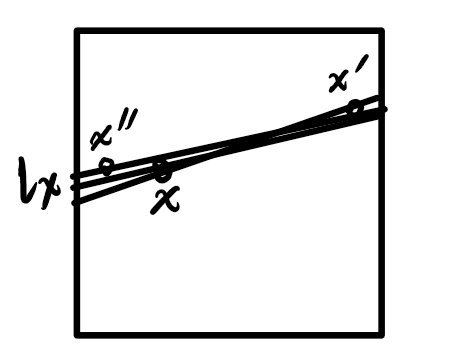}\quad\includegraphics[width=.23\linewidth]{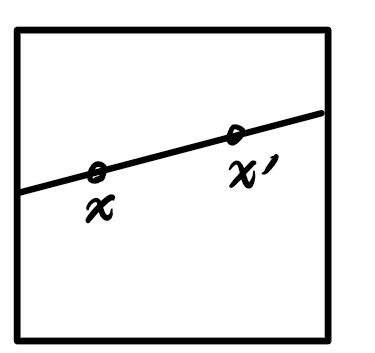}\end{center}
\paragraph{Circular/parabolic Nikodym sets: do not exist.}
\begin{center}\includegraphics[width=.65\linewidth]{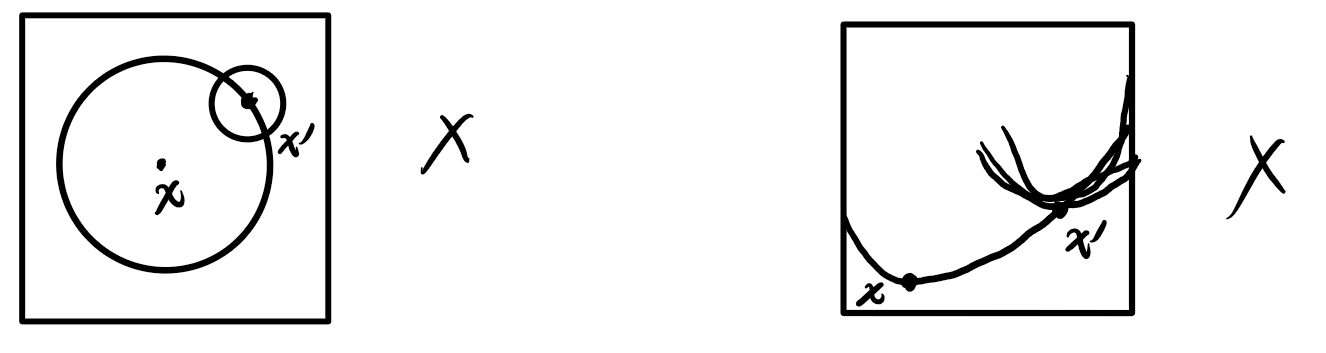}\end{center}
\paragraph{Another parabolic Nikodym set: exists.}
Consider: for every $x\notin E$, there exists a parabola containing $x$ of the form $t^2+bt+c$. It satisfies the geometric cinematic curvature condition.
On the other hand, it can be transformed to a (line) Nikodym set as follows:
\begin{center}\includegraphics[width=.23\linewidth]{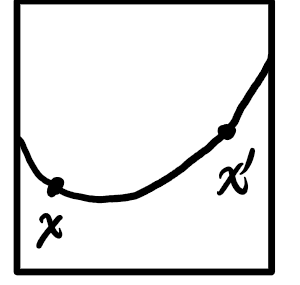}\end{center}
\[
 (X,Y)\mapsto(X,Y+X^2).
\]
\begin{align*}
 (x+t,y+ut)&\mapsto(x+t,y+ut+(x+t)^2)\\
 &=(x+t,y+x^2+(2x+u)t+t^2).
\end{align*}
\paragraph{Cubic Nikodym sets.}
Consider
\[
 (x+t,y+ut+t^3)\quad\text{(do not exist)}
 \qquad\text{vs.}\qquad
 (x+t,y+ut+vt^2+t^3)\quad\text{(exist)}.
\]
\[
 (X,Y)\mapsto(X,Y+X^3),
\]
\begin{align*}
 (x+t,y+ut)&\mapsto(x+t,y+ut+(x+t)^3)\\
 &=(x+t,y+x^3+(3x^2+u)t+\boxed{3xt^2}+t^3).
\end{align*}
It can be verified that the first form does not satisfy the geometric cinematic curvature condition.

\section{The Fourier restriction problems on \texorpdfstring{$\mathbb R^2$}{R2}}
For the classical planar restriction theorem and its background, see~\cite{Fefferman1970,Stein1993,Demeter2020}.
\subsection{Motivation}
For a function $f$, can $\widehat f$ be restricted to the sphere $S^{n-1}$?
\begin{itemize}
\item If $f\in L^1$, $\widehat f\in C_0$, and $\widehat f|_{S^{n-1}}$ can be defined pointwise.
\item If $f\in L^2$, $\widehat f\in L^2$ by the Plancherel theorem. Moreover, $f\mapsto\widehat f$ is a bijection. So $\widehat f|_{S^{n-1}}$ cannot be well-defined, since $\mathcal L^n(S^{n-1})=0$ and functions in $L^2(\mathbb R^n)$ are allowed to change values on any sets of measure zero.
\item If $f\in L^p$, $p>2$, $\widehat f$ is defined as a tempered distribution, and then $\widehat f|_{S^{n-1}}$ is meaningless.
\item If $f\in L^p$, $1<p<2$, we know $\widehat f$ can be defined as follows:
\[
f=f_1+f_2\in L^1+L^2\quad\Longrightarrow\quad
\widehat f:=\widehat f_1+\widehat f_2.
\]
\end{itemize}

\paragraph{Question.} What is the largest $p$ that makes $\widehat f|_{S^{n-1}}$ well-defined?
To study this question, it is sufficient to prove
\[
\|\widehat f\|_{L^q(S^{n-1})}\le C_{p,q}\|f\|_{L^p(\mathbb R^n)}
\]
for some $(p,q)$ and any Schwartz function $f\in\mathcal S(\mathbb R^n)$. If we have the above estimate, by the standard density argument we can extend the definition of $\widehat f|_{S^{n-1}}$ for $f\in\mathcal S(\mathbb R^n)$ to $f\in L^p(\mathbb R^n)$.
\subsection{Fourier restriction conjecture (\texorpdfstring{$n=2$}{n=2})}
\[
\|\widehat f\|_{L^q(S^1)}\le C_{p,q}\|f\|_{L^p(\mathbb R^2)}
\]
holds for $1\le p<\frac43$ and $1\le q\le\frac{p'}3$.
\paragraph{Duality: Fourier extension operator.}
\[
Eg(x,y)=E_{\mathbb P^1}g(x,y):=
\int_0^1g(\theta)e^{i(\theta x+\theta^2y)}\,d\theta.
\]

Define $Rf=\widehat f|_{\mathbb P^1}$, identified with the function $\theta\mapsto\widehat f(\theta,\theta^2)$ on $[0,1]$. For $f\in\mathcal S(\mathbb R^2)$ and $g\in C_c^\infty((0,1))$, it is easy to check that
\[
\langle Rf,g\rangle=\langle f,Eg\rangle
.
\]
\paragraph{Conjecture.}
\[
\|Eg\|_{L^q(\mathbb R^2)}\le C_{p,q}\|g\|_{L^p([0,1])}
\]
holds for $q>4$ and $\frac1p+\frac3q\le1$.
\begin{center}\includegraphics[width=.38\linewidth]{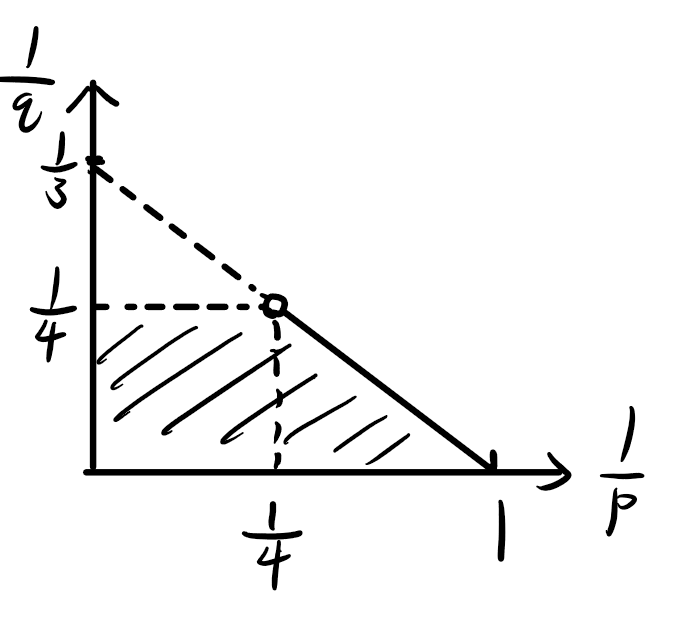}\end{center}
Moreover, for the endpoint $q=4$, $p=4$, we expect that there is an $\varepsilon$-loss locally:
\[
\|Eg\|_{L^4(B_R)}\le C_\varepsilon R^\varepsilon\|g\|_{L^4([0,1])}
\]
for any $R\ge1$ and $\varepsilon>0$.

\subsection{Examples}
\paragraph{1.} $f=1$. Then $\|f\|_p=1$. When $(x,y)$ is away from the origin, by the stationary phase principle,
\[
|Ef(x,y)|=\left|\int_0^1e^{i(x\theta+y\theta^2)}\,d\theta\right|
\gtrsim |y|^{-1/2}\sim |(x,y)|^{-1/2}.
\]
Thus $Ef\in L^q$ only if $q>4$.
Next we show why we need $\frac1p+\frac3q\le1$.
\[
\int_0^1e^{i(x\theta+y\theta^2)}\,d\theta
=\sum_{\ell(I)=R^{-1/2}}\int_Ie^{i(x\theta+y\theta^2)}\,d\theta.
\]
\begin{center}\includegraphics[width=.85\linewidth]{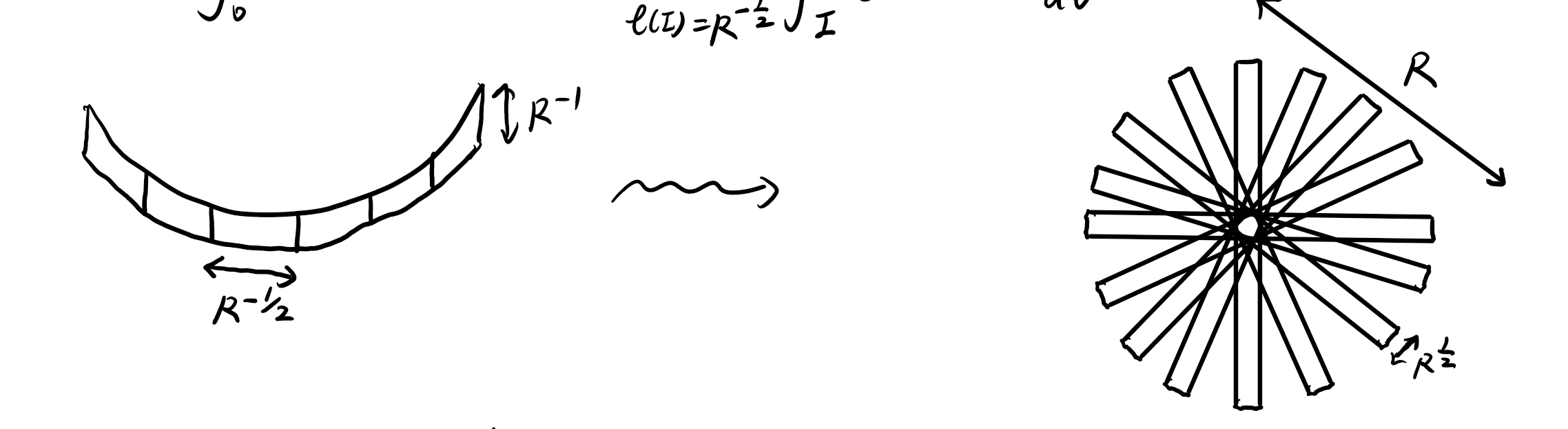}\end{center}
\paragraph{Special $I$:} $I=[0,R^{-1/2}]$.

$|y|\le R$, $|\theta|\le R^{-1/2}$ imply $|y\theta^2|\le1$. To obtain $|x\theta|\le1$, we require $|x|\le R^{1/2}$. So
\[
\int_Ie^{i(x\theta+y\theta^2)}\,d\theta\sim |I|=R^{-1/2}
\quad\text{if }|x|\le R^{1/2},\ |y|\le R.
\]
More precisely,
\begin{align*}
\int_0^{R^{-1/2}}e^{i(x\theta+y\theta^2)}\,d\theta
&=R^{-1/2}\int_0^1e^{i(xR^{-1/2}\theta+yR^{-1}\theta^2)}\,d\theta\\
&\sim\begin{cases}R^{-1/2},&|x|\lesssim R^{1/2},\\0,&|x|\gtrsim R^{1/2},\end{cases}
\end{align*}
by van der Corput. (The change of variables is $\theta\mapsto R^{1/2}\theta$.)

This gives $Ef_I\sim R^{-1/2}1_{T_I}$, where $T_I$ is the tube $[0,R^{1/2}]\times[0,R]$ dual to $I$.
\paragraph{General $I$:} $I=[c_I,c_I+R^{-1/2}]$.
Using $\theta\mapsto R^{1/2}(\theta-c_I)$,
\begin{align*}
\int_Ie^{i(x\theta+y\theta^2)}\,d\theta
&=R^{-1/2}\int_0^1e^{i\{x(R^{-1/2}\theta+c_I)+y(R^{-1/2}\theta+c_I)^2\}}\,d\theta\\
&=R^{-1/2}e^{i(xc_I+yc_I^2)}\int_0^1
 e^{i[R^{-1/2}\theta(x+2c_Iy)+R^{-1}\theta^2 y]}\,d\theta.
\end{align*}
Here $x+2c_Iy=(x,y)\cdot(1,2c_I)$, and the last phase term is at most $1$. We need
\[
|(x,y)\cdot(1,2c_I)|\lesssim R^{1/2}.
\]
Note that $(1,2c_I)\perp(2c_I,-1)$, so
\[
\angle((x,y),(2c_I,-1))\lesssim R^{-1/2}.
\]
Since $|y|\le R$, we obtain the width $\sim R^{1/2}$.
\begin{center}\includegraphics[width=.65\linewidth]{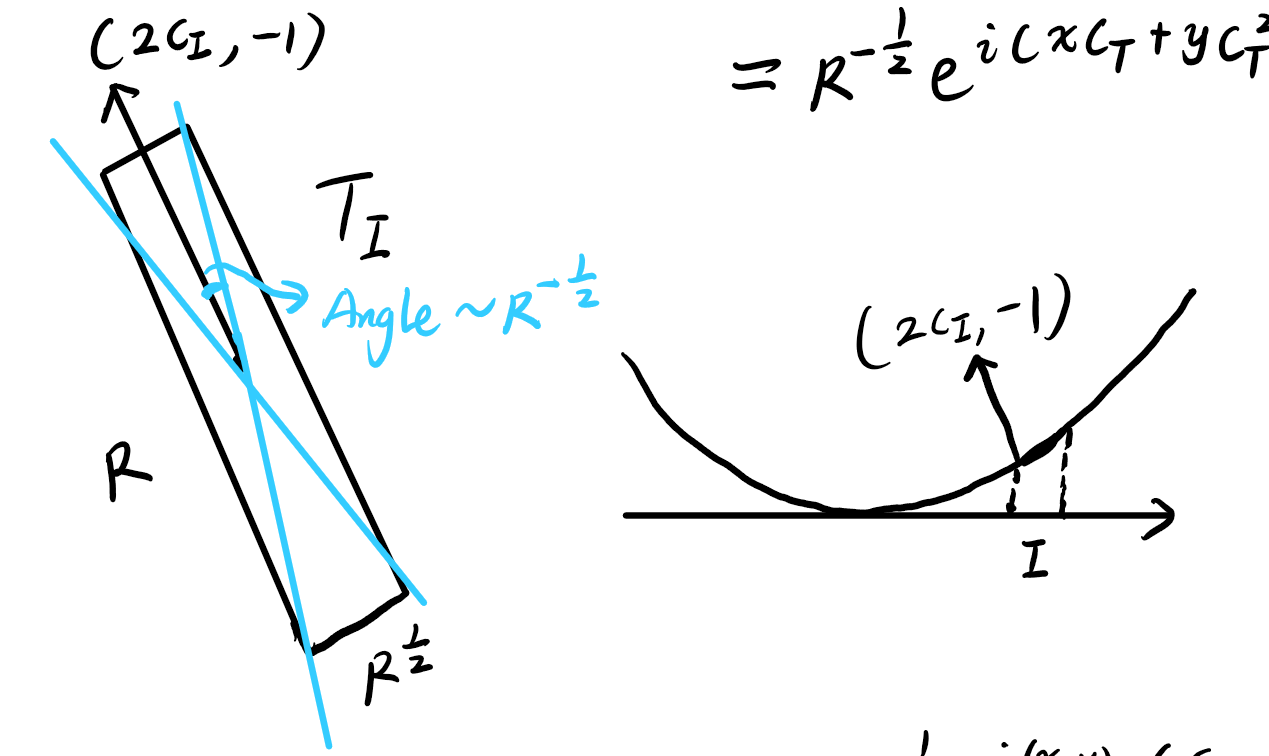}\end{center}
This gives
\[
Ef_I\sim R^{-1/2}e^{i(x,y)\cdot(c_I,|c_I|^2)}1_{T_I}.
\]
The decomposition $Ef=\sum_I Ef_I\sim\sum_Ia_{T_I}1_{T_I}$ is actually the wave packet decomposition, and each $Ef_I=a_{T_I}1_{T_I}$ is a single wave packet.
\paragraph{2. Knapp example.}

Take $f=1_{[0,R^{-1/2}]}$. Then $Ef\sim R^{-1/2}1_{[0,R^{1/2}]\times[0,R]}$ is a single wave packet. Thus
\[
\|Ef\|_q\sim R^{-1/2}R^{3/(2q)},\qquad
\|f\|_p=R^{-1/(2p)}.
\]
Hence
\[
\frac{\|Ef\|_q}{\|f\|_p}\sim R^{\frac3{2q}-\frac12+\frac1{2p}}
\quad\Longrightarrow\quad
\frac3{2q}-\frac12+\frac1{2p}\le0
\quad\Longrightarrow\quad\frac1p+\frac3q\le1.
\]
\subsection{Stein--Tomas inequality}
For the classical Stein--Tomas theorem, see~\cite{Tomas1975,Stein1993}; the local proof below uses~\cite{BD2015}.
\paragraph{Theorem 1.} For the pair $(2,6)$, we have
\[
\|Ef\|_{L^6(B_R)}\lesssim_\varepsilon R^\varepsilon\|f\|_2.
\]
\paragraph{Remark.} The classical method to prove this inequality is the $TT^*$ method. Here we use decoupling to give a very simple proof.
\paragraph{Proof.} Apply decoupling for the parabola:
\[
\|Ef\|_{L^6(B_R)}\lesssim_\varepsilon R^\varepsilon
\left(\sum_{\ell(I)=R^{-1/2}}\|Ef_I\|_{L^6(B_R)}^2\right)^{1/2}.
\]

Next, for each $I$, we aim to prove
\[
\|Ef_I\|_{L^6(B_R)}\lesssim\|f_I\|_2.
\]
WLOG, set $I=[0,R^{-1/2}]$. We expand $f_I$ as a Fourier series:
\[
f_I(\theta)=\sum_{n\in\mathbb Z}a_ne^{inR^{1/2}\theta}.
\]
Then
\begin{align*}
Ef_I(x,y)&=\sum_{n\in\mathbb Z}a_n\int_0^{R^{-1/2}}
 e^{inR^{1/2}\theta}e^{i(x\theta+y\theta^2)}\,d\theta\\
&\sim\sum_{n\in\mathbb Z}a_nR^{-1/2}1_{T_n},
\end{align*}
where $T_n=\{|x+nR^{1/2}|\le R^{1/2},\ |y|\le R\}$.
By Parseval's inequality,
\[
\|f_I\|_2^2\simeq\sum_{n\in\mathbb Z}|a_n|^2R^{-1/2}.
\]
\begin{center}
\begin{tikzpicture}
\node[anchor=south west,inner sep=0] (image) at (0,0)
  {\includegraphics[width=.62\linewidth]{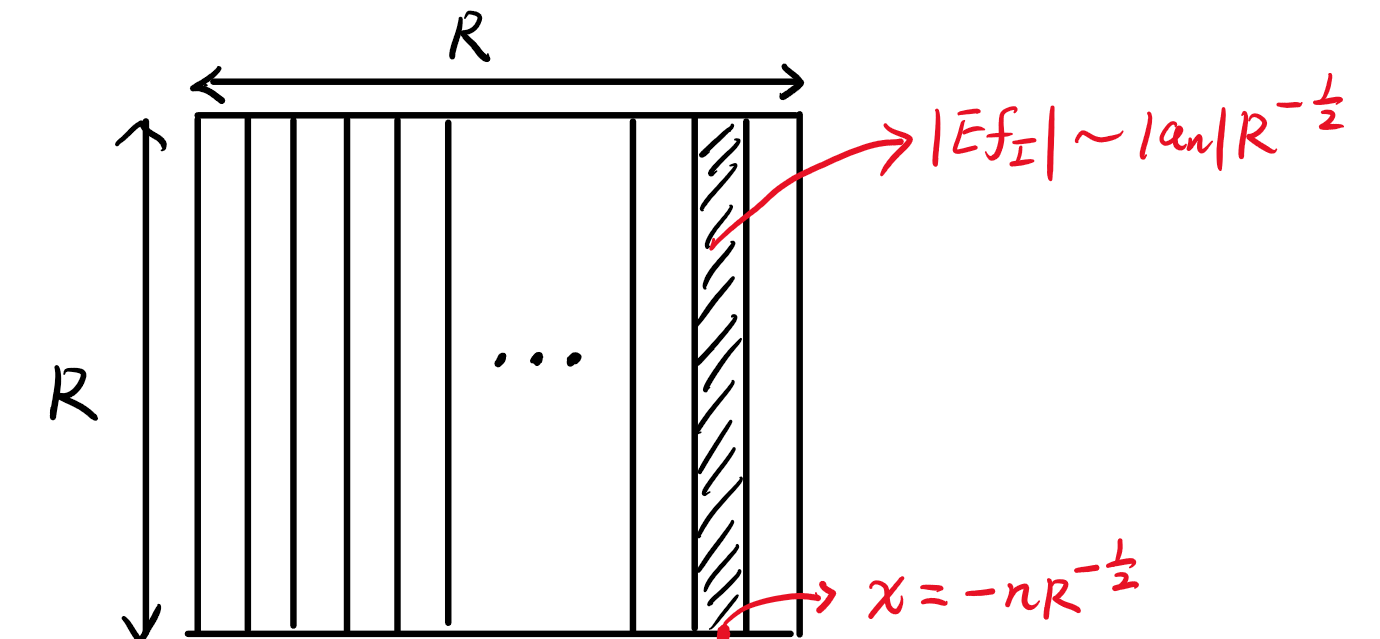}};
\begin{scope}[x={(image.south east)},y={(image.north west)}]
\fill[white] (.61,0) rectangle (1,.17);
\node[anchor=west,text=red!85!black] at (.62,.075) {$x=-nR^{1/2}$};
\end{scope}
\end{tikzpicture}
\end{center}
The annotation on the selected tube is $|Ef_I|\sim|a_n|R^{-1/2}$; its bottom coordinate is $x=-nR^{1/2}$.
Consequently,
\begin{align*}
\|Ef_I\|_{L^6(B_R)}
&\simeq\left(\sum|a_n|^6R^{-3}\cdot R^{3/2}\right)^{1/6}\\
&\sim R^{-1/4}\left(\sum|a_n|^6\right)^{1/6}
\le R^{-1/4}\left(\sum|a_n|^2\right)^{1/2}.
\end{align*}

The last expression is $\sim\|f_I\|_2$ (the width of each tube in the preceding diagram is $R^{1/2}$). Finally,
\[
\|Ef\|_{L^6(B_R)}\lesssim_\varepsilon R^\varepsilon
\left(\sum_I\|f_I\|_2^2\right)^{1/2}\sim R^\varepsilon\|f\|_2.\qquad\square
\]
\subsection{Refined decoupling}
In the above proof of the Stein--Tomas inequality,
\[
\left(\sum|a_n|^6\right)^{1/6}\le\left(\sum|a_n|^2\right)^{1/2}
\]
has caused a large waste. The only way to obtain equality is $a_n=0$ for all $n\in\mathbb N$ except one. Equivalently, $Ef_I$ concentrates on exactly one tube.
\paragraph{Assumption.}
\[
Ef=\sum_I Ef_I=\sum_I\sum_{T_I\parallel I^*}Ef_{T_I}
\]
satisfies:
\begin{enumerate}
\item For all $I$, $|Ef_I|$ is comparable $\Longleftrightarrow |a_n|$ is comparable.
\item There exists $M$ such that $\#\{T_I\}\sim M$ for all $I$.
\end{enumerate}
Then for each $I$,
\[
\frac{(\sum|a_n|^6)^{1/6}}{(\sum|a_n|^2)^{1/2}}
\sim\frac{M^{1/6}}{M^{1/2}}=M^{-1/3}.
\]

This implies
\[
\|Ef_I\|_{L^6(B_R)}\lesssim M^{-1/3}\|f_I\|_2.
\]
Therefore, under this assumption, we obtain a refined estimate:
\[
\|Ef\|_{L^6(B_R)}\lesssim_\varepsilon R^\varepsilon M^{-1/3}\|f\|_2.
\]
\subsection{Applications of the Stein--Tomas inequality}
\paragraph{1. Strichartz estimate for the Schr\"odinger equation.} \cite{Strichartz1977}
Consider
\[
\begin{cases}
i\partial_tu(x,t)-\partial_x^2u(x,t)=0&\text{in }\mathbb R\times(0,\infty),\\
u(x,0)=u_0(x)&\text{on }\mathbb R.
\end{cases}
\]
It has solution
\[
u(x,t)=(e^{it\Delta}u_0)(x)=E(\widehat{u_0})(x,t).
\]
We have

\[
\|u\|_{L^6(\mathbb R\times(0,\infty))}\lesssim\|\widehat{u_0}\|_2=\|u_0\|_2,
\]
and also have a refined Strichartz estimate if $u_0$ satisfies an extra condition.
\paragraph{2. Wave equation.}
\[
\begin{cases}
\partial_t^2u-\Delta_xu=0&\text{in }\mathbb R^2\times(0,\infty),\\
u(\cdot,0)=u_0,\quad\partial_tu(\cdot,0)=0&\text{on }\mathbb R^2.
\end{cases}
\]
For the following local-in-time estimates, assume
\[
\operatorname{supp}\widehat{u_0}
\subset\{\xi\in\mathbb R^2:N/2\leq|\xi|\leq2N\},
\qquad N\geq1,
\]
and restrict time to $[1,2]$. Thus
$u(t)=\cos(t|D|)u_0$, where $|D|=(-\Delta_x)^{1/2}$.
\begin{itemize}
\item Strichartz estimate: \cite{Strichartz1977}
\[
\|u\|_{L^6(\mathbb R^2\times[1,2])}
\lesssim N^{1/2}\|u_0\|_2.
\]
\item Local smoothing estimate (Guth--Wang--Zhang, high--low method): \cite{GWZ2020}
for every $\varepsilon>0$,
\[
\|u\|_{L^4(\mathbb R^2\times[1,2])}
\lesssim_\varepsilon N^\varepsilon\|u_0\|_4.
\]
\item Mixed norm Strichartz estimate: \cite{Wolff1997} for every $\varepsilon>0$,
\[
\big\|\|u(x,t)\|_{L_x^\infty(\mathbb R^2)}\big\|_{L_t^3([1,2])}
\lesssim_\varepsilon N^{1/2+\varepsilon}\|u_0\|_3.
\]
This is related to Wolff's circular maximal estimate. If $f$ has
the same frequency localization, the corresponding circle-average
estimate is
\[
\left\|\left\|\int_{S^1}f(x+ty)\,d\sigma(y)
\right\|_{L_x^\infty(\mathbb R^2)}\right\|_{L_t^3([1,2])}
\lesssim_\varepsilon N^\varepsilon\|f\|_3.
\]
The connection uses stationary phase: at frequency $N$ and
$t\in[1,2]$, circular averaging carries an additional amplitude
factor $N^{-1/2}$ relative to the half-wave propagators.
In the thickened-circle formulation, the thickness is
$\delta\sim N^{-1}$ and the corresponding loss is
$\delta^{-\varepsilon}$.
\end{itemize}
\subsection{\texorpdfstring{$L^4$}{L4} estimate}
For the bilinear restriction method in the plane, see~\cite{Fefferman1970,Demeter2020}.
Now we prove the $L^4$ estimate by the bilinear method.
\paragraph{Theorem 2.} For the pair $(4,4)$, we have
\[
\|Ef\|_{L^4(B_R)}\lesssim(\log R)^{1/4}\|f\|_4.
\]
\paragraph{Proof (work at scale $R^{-1}$).}
\begin{align*}
\|Ef\|_{L^4(B_R)}^4
&=\|Ef\cdot\overline{Ef}\|_{L^2(B_R)}^2\\
&=\left\|\iint_{[0,1]^2}f(\xi)\overline{f(\eta)}
 e^{i[x_1(\xi-\eta)+x_2(\xi^2-\eta^2)]}\,d\xi\,d\eta\right\|_{L^2(B_R)}^2.
\end{align*}
Let $F(\xi,\eta)=f(\xi)\overline{f(\eta)}$, $y_1=\xi-\eta$, $y_2=\xi^2-\eta^2$, and $[0,1]^2\mapsto\Omega$. The last expression is
\[
\left\|\iint_\Omega F(\xi(y),\eta(y))e^{i(x_1y_1+x_2y_2)}
\frac{1}{|\xi(y)-\eta(y)|}\,dy\right\|_{L^2(B_R)}^2.
\]
Choose a smooth function $\phi$ with Fourier transform supported on $B(0,1)$ and such that $\phi\ge1$ on $B(0,1)$. Let $\phi_R=\phi(R^{-1}\cdot)$.

The preceding expression is bounded by
\begin{align*}
&\left\|\left(\iint_\Omega F(\xi(y),\eta(y))e^{i(x_1y_1+x_2y_2)}
 \frac{1}{|\xi(y)-\eta(y)|}\,dy\right)\phi_R\right\|_{L^2(\mathbb R^2)}^2\\
&\qquad=\left\|\left(1_\Omega F(\xi(y),\eta(y))\frac{1}{|\xi(y)-\eta(y)|}\right)*\widehat\phi_R\right\|_{L^2(\mathbb R^2)}^2
\end{align*}
by Plancherel.
Apply the Whitney decomposition to $[0,1]\setminus\{\xi=\eta\}$ at scale $R^{-1}$:
\[
|\xi-\eta|\sim2^{-k}\rightsquigarrow\Omega_k,\quad k=1,2,\ldots,O(\log R),
\qquad |\xi-\eta|\le R^{-1}\rightsquigarrow\Omega_0.
\]
So $[0,1]^2\mapsto\Omega=\Omega_0\cup\bigcup_{k=1}^{O(\log R)}\Omega_k$.
\begin{center}\includegraphics[width=.25\linewidth]{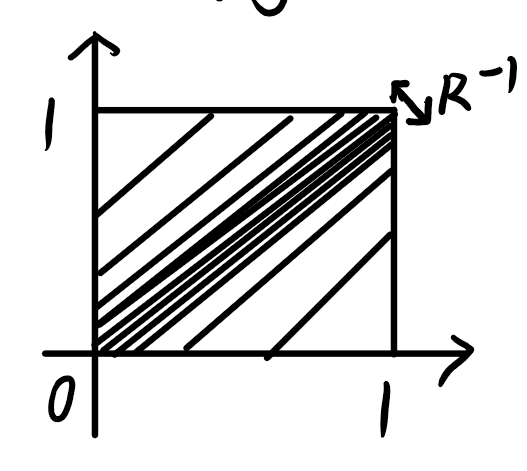}\end{center}
\begin{align*}
\|Ef\|_{L^4(B_R)}^4
&=\left\|\left(1_{\Omega_0}F(\xi(y),\eta(y))\frac1{|\xi(y)-\eta(y)|}\right)*\widehat\phi_R\right\|_{L^2(\mathbb R^2)}^2\\
&\quad+\sum_{k=1}^{O(\log R)}\left\|\left(1_{\Omega_k}F(\xi(y),\eta(y))\frac1{|\xi(y)-\eta(y)|}\right)*\widehat\phi_R\right\|_{L^2(\mathbb R^2)}^2.
\end{align*}
\paragraph{Low frequency $\Omega_0$.}

By Plancherel,
\begin{align*}
&\left\|\left(1_{\Omega_0}F(\xi(y),\eta(y))\frac1{|\xi(y)-\eta(y)|}\right)*\widehat\phi_R\right\|_{L^2(\mathbb R^2)}\\
&\quad\sim\left\|\iint_{\Omega_0}F(\xi(y),\eta(y))\frac1{|\xi(y)-\eta(y)|}e^{ix\cdot y}\,dy\right\|_{L^2(B_R)}\\
&\quad=\left\|\iint_{|\xi-\eta|\le R^{-1}}F(\xi,\eta)
 e^{i[x_1(\xi-\eta)+x_2(\xi^2-\eta^2)]}\,d\xi\,d\eta\right\|_{L^2(B_R)}.
\end{align*}
(In this case, $|\xi-\eta|\le R^{-1}$, $|\xi^2-\eta^2|=|\xi+\eta||\xi-\eta|\le R^{-1}$, and $|x|\le R$, so there is actually nothing about oscillations.)
The last expression is
\begin{align*}
&\le\left\|\iint_{|\xi-\eta|\le R^{-1}}|F(\xi,\eta)|\,d\xi\,d\eta\right\|_{L^2(B_R)}\\
&=|B_R|^{1/2}\iint_{|\xi-\eta|\le R^{-1}}|f(\xi)||f(\eta)|\,d\xi\,d\eta\\
&\le R^{1/2}\left(\iint_{|\xi-\eta|\le R^{-1}}|f(\xi)|^2|f(\eta)|^2\,d\xi\,d\eta\right)^{1/2}\\
&\sim R^{1/2}\left(\sum_j\int_{I_j}|f|^2\int_{J_j}|f|^2\right)^{1/2},
\end{align*}
where $\{|\xi-\eta|\le R^{-1}\}=\bigcup_{j=1}^R I_j\times J_j$, with side lengths $R^{-1}$.
\begin{center}\includegraphics[width=.24\linewidth]{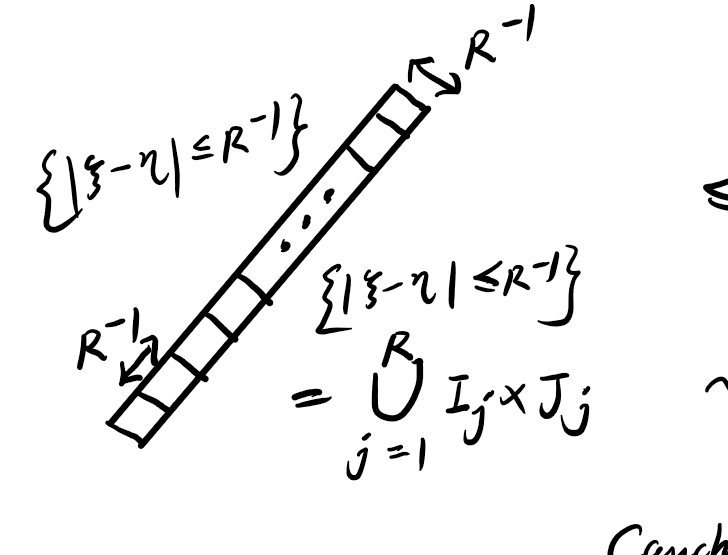}\end{center}

By Cauchy--Schwarz, this is at most
\begin{align*}
&R^{1/2}|I_j|^{1/4}|J_j|^{1/4}
 \left(\sum_j\left(\int_{I_j}|f|^4\right)^{1/2}\left(\int_{J_j}|f|^4\right)^{1/2}\right)^{1/2}\\
&\quad\le\left(\sum_j\int_{I_j}|f|^4\right)^{1/4}
 \left(\sum_j\int_{I_j}|f|^4\right)^{1/4}=\|f\|_4^2.
\end{align*}
\paragraph{High frequency $\Omega_k$.}
By Young's inequality,
\begin{align*}
&\left\|\left(1_{\Omega_k}F(\xi(y),\eta(y))\frac1{|\xi(y)-\eta(y)|}\right)*\widehat\phi_R\right\|_{L^2(\mathbb R^2)}\\
&\quad\le\|\widehat\phi_R\|_1\left(\int_{\Omega_k}|F(\xi(y),\eta(y))|^2
 \frac1{|\xi(y)-\eta(y)|^2}\,dy\right)^{1/2}\\
&\quad\sim\left(\iint_{|\xi-\eta|\sim2^{-k}}|F(\xi,\eta)|^2\frac1{|\xi-\eta|}\,d\xi\,d\eta\right)^{1/2}\\
&\quad\sim2^{k/2}\left(\sum_j\int_{I_j}|f|^2\int_{J_j}|f|^2\right)^{1/2}\\
&\quad\le2^{k/2}|I_j|^{1/4}|J_j|^{1/4}
 \left(\sum_j\int_{I_j}|f|^4\right)^{1/4}
 \left(\sum_j\int_{J_j}|f|^4\right)^{1/4}.
\end{align*}
Here $\{|\xi-\eta|\sim2^{-k}\}=\bigcup_jI_j\times J_j$, with side lengths $2^{-k}$; the last step is Cauchy--Schwarz.
\begin{center}\includegraphics[width=.24\linewidth]{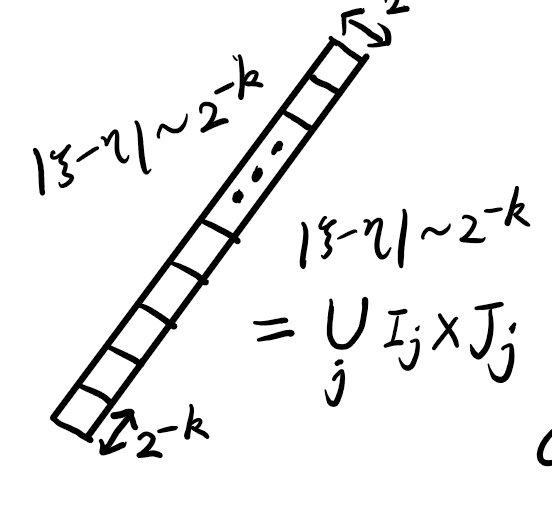}\end{center}

The last expression is $\le\|f\|_4^2$. In summary, we have
\[
\|Ef\|_{L^4(B_R)}\lesssim(\log R)^{1/4}\|f\|_4.\qquad\square
\]
\paragraph{Another proof of Theorem 2 (work at scale $R^{-1/2}$).}
By a similar argument as before, we have the following bilinear estimate:
\[
\big\||E_If\,E_Jf|^{1/2}\big\|_{L^4(B_R)}^4
\le\operatorname{dist}(I,J)^{-1}\|f_I\|_2^2\|f_J\|_2^2,
\]
for any intervals $I,J$ with $\operatorname{dist}(I,J)>0$.
Apply the Whitney decomposition to $[0,1]^2\setminus\{\xi=\eta\}$ at scale $R^{-1/2}$: draw dyadic cubes; if they touch $\{\xi=\eta\}$, do a dyadic partition again. Let
\[
\mathcal W_l=\{\text{dyadic cubes of length }\sim2^{-l}\},
\qquad l=1,2,\ldots,\frac{\log R}{2\log2}.
\]
\begin{center}\includegraphics[width=.32\linewidth]{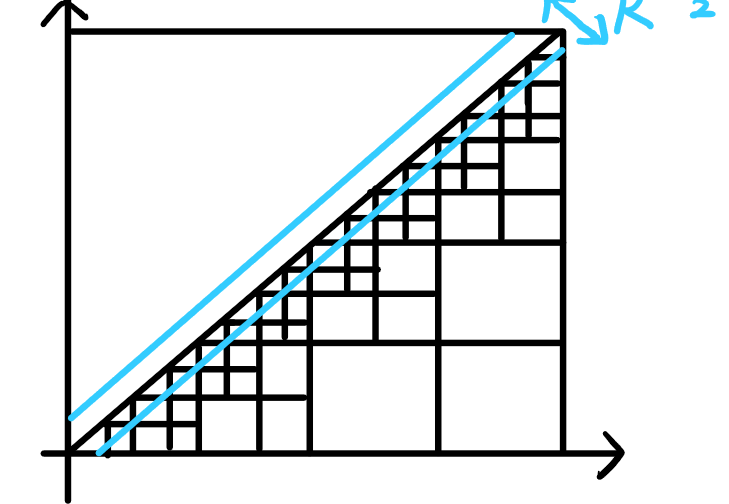}\end{center}

Let $\mathcal W_0$ be the $R^{-1/2}$-neighborhood of the diagonal line $\{\xi=\eta\}$.
\paragraph{Observations.}
\begin{enumerate}
\item $[0,1]^2=\mathcal W_0\cup\bigcup_{l=1}^{O(\log R)}\bigcup_{I_1\times I_2\in\mathcal W_l}I_1\times I_2$.
\item Suppose $Q=I_1\times I_2$ is a dyadic cube. Then $\operatorname{dist}(I_1,I_2)\sim\ell(Q)$.
\end{enumerate}
So we have
\begin{align*}
\|Ef\|_{L^4(B_R)}^2
&=\|Ef\,\overline{Ef}\|_{L^2(B_R)}\\
&=\left(\iint_{B_R}\left|\iint_{[0,1]^2}f(\xi)\overline{f(\eta)}
 e^{ix\cdot(\xi-\eta,\xi^2-\eta^2)}\,d\xi\,d\eta\right|^2dx\right)^{1/2}\\
&=\left(\iint_{B_R}\left|\sum_l\sum_{I_1\times I_2\in\mathcal W_l}
 \iint_{I_1\times I_2}\cdots\,d\xi\,d\eta\right|^2dx\right)^{1/2}\\
&\le\sum_l\left(\iint_{B_R}\left|\sum_{I_1\times I_2\in\mathcal W_l}
 \iint_{I_1\times I_2}\cdots\,d\xi\,d\eta\right|^2dx\right)^{1/2}
\end{align*}
by Minkowski. Since $\#l\sim\log R$, by pigeonholing it is sufficient to show, for any $l$, the following estimate.

\begin{align*}
&\left\|\sum_{I_1\times I_2\in\mathcal W_l}\iint_{I_1\times I_2}
 f(\xi)\overline{f(\eta)}e^{ix\cdot(\xi-\eta,\xi^2-\eta^2)}\,d\xi\,d\eta\right\|_{L^2(B_R)}\\
&\qquad=\left\|\sum_{I_1\times I_2\in\mathcal W_l}Ef_{I_1}\overline{Ef_{I_2}}\right\|_{L^2(B_R)}
\lesssim\|f\|_4^2.
\end{align*}
(The phase is $x_1\xi+x_2\xi^2-(x_1\eta+x_2\eta^2)$.)
When $l\ge1$, note that
\[
\operatorname{supp}\widehat{(Ef_{I_1}\overline{Ef_{I_2}})}\subset\theta_{I_1}+\theta_{I_2}.
\]
If $I_1\times I_2,I'_1\times I'_2\in\mathcal W_l$, then
\[
\operatorname{dist}(\theta_{I_1}+\theta_{I_2},\theta_{I'_1}+\theta_{I'_2})
\ge2^{-l}\ge R^{-1/2}\gg R^{-1}.
\]
\begin{center}\includegraphics[width=.35\linewidth]{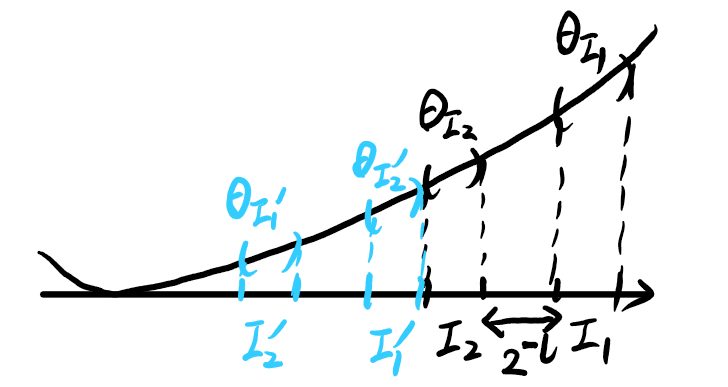}\end{center}
Therefore $\sum Ef_{I_1}\overline{Ef_{I_2}}$ is almost orthogonal, i.e.,
\begin{align*}
\left\|\sum_{I_1\times I_2\in\mathcal W_l}Ef_{I_1}\overline{Ef_{I_2}}\right\|_{L^2(B_R)}
&\lesssim\left(\sum_{I_1\times I_2\in\mathcal W_l}\|Ef_{I_1}\overline{Ef_{I_2}}\|_{L^2(B_R)}^2\right)^{1/2}\\
&\le\left(\sum_{I_1\times I_2\in\mathcal W_l}\operatorname{dist}(I_1,I_2)^{-1}
 \|f_{I_1}\|_2^2\|f_{I_2}\|_2^2\right)^{1/2}.
\end{align*}

This is $\le\|f\|_4^2$. When $l=0$, i.e., $|\xi-\eta|\le R^{-1/2}$, we have
\begin{align*}
\left\|\sum_{\mathcal W_0}Ef_{I_1}\overline{Ef_{I_2}}\right\|_{L^2(B_R)}
&=\left\|\sum_{\ell(I)=R^{-1/2}}|Ef_I|^2\right\|_{L^2(B_R)}\\
&\sim\left(\sum_I\|Ef_I\|_{L^4(B_R)}^4\right)^{1/2}.
\end{align*}
Here $\mathcal W_0=\bigcup_{\ell(I)=R^{-1/2}}I\times I$.
\begin{center}\includegraphics[width=.25\linewidth]{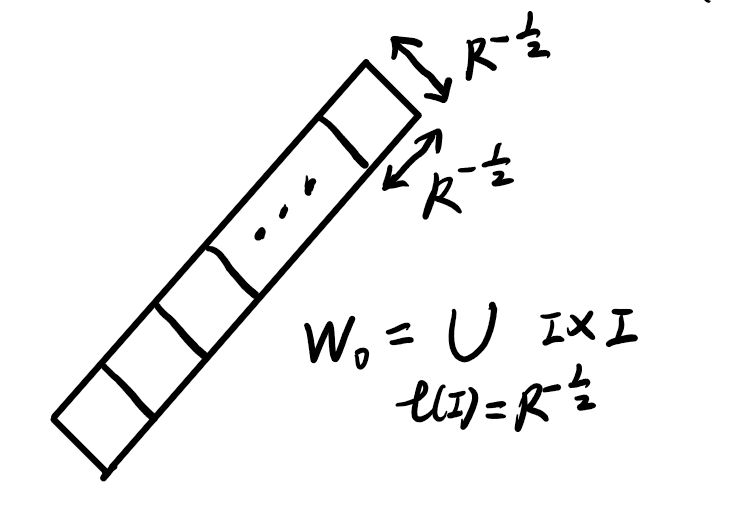}\end{center}
For each $I$, consider the wave packet decomposition for $Ef_I$:
\[
f_I\text{ is supported on }I
\quad\Longrightarrow\quad
\widehat{f_I}=\sum_{\ell(J_k)=R^{1/2}}b_k1_{J_k},
\]
and
\[
Ef_I=\sum_kb_kE(\check1_{J_k})=\sum_kb_k1_{T_k},
\]
where $T_k$ is a tube of dimensions $R^{1/2}\times R$.
\begin{center}\includegraphics[width=.38\linewidth]{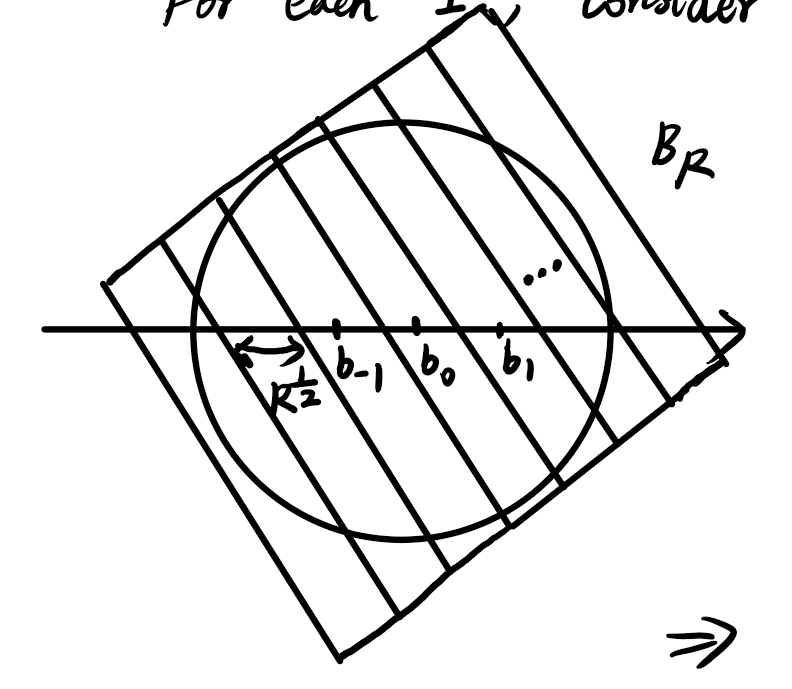}\end{center}
Thus
\[
\|Ef_I\|_{L^4(B_R)}\sim\left(\sum_k|b_k|^4R^{3/2}\right)^{1/4}.
\]

This equals $R^{3/8}(\sum_k|b_k|^4)^{1/4}$. Also,
\[
\|f_I\|_2=\|\widehat{f_I}\|_2\sim\left(\sum_k|b_k|^2R^{1/2}\right)^{1/2}
=R^{1/4}\left(\sum_k|b_k|^2\right)^{1/2}.
\]
Hence we get
\[
\|Ef_I\|_{L^4(B_R)}\lesssim R^{1/8}\|f_I\|_2\lesssim\|f_I\|_4,
\]
and
\[
\left\|\sum_{\mathcal W_0}|Ef_I|^2\right\|_{L^2(B_R)}
\lesssim\left(\sum_I\|f_I\|_4^4\right)^{1/2}=\|f\|_4^2.\qquad\square
\]
\paragraph{Compared to scale $R^{-1}$.}
If we work at a smaller scale $R^{-1}$, then $Ef_I\approx\text{constant}$ on $B_R$, so we just expect the following $L^1$ estimate:
\begin{align*}
\|Ef_I\|_{L^4(B_R)}
&\le|B_R|^{1/4}\|Ef_I\|_\infty
\le R^{1/2}\|f_I\|_1\\
&\le R^{1/2}(R^{-1})^{1-1/4}\|f_I\|_4
=R^{-1/4}\|f_I\|_4.
\end{align*}
We obtain a gain $R^{-1/4}$!

At scale $R^{-1/2}$, we use the wave packet decomposition to apply $L^2$ orthogonality. However, this only separates wave packets with wave speed $R^{-1/2}$ travelling in the same direction; it cannot distinguish smaller wave speeds. In fact, the smallest distinguishable wave speed should be $R^{-1}$, so working at scale $R^{-1}$ gives an additional gain.

From the computational point of view, at scale $R^{-1/2}$ we still have the oscillatory term $e^{i(x\cdot(\xi-\eta,\xi^2-\eta^2))}$, and handling it using the wave packet decomposition introduces some losses.

\section{The Bochner--Riesz problem on \texorpdfstring{$\mathbb R^2$}{R2}}

\subsection{Motivation: convergence of Fourier series}
Suppose $f:[0,1]\to\mathbb C$, and
\[
a_n=\widehat f(n)=\int_0^1f(t)e^{-2\pi int}\,dt.
\]
Do we have $f(x)\sim\sum_{n\in\mathbb Z}a_ne^{2\pi inx}$?
\paragraph{Convergence in $L^p$ norm.}
Let
\[
S_Nf(x):=\sum_{n=-N}^Na_ne^{2\pi inx}.
\]
Does $S_Nf\to f$ in $L^p$? The uniform boundedness principle shows that it suffices to prove
\[
\|S_Nf\|_p\lesssim_p\|f\|_p.
\]
\paragraph{Transference principle.}
By routine scaling and limiting arguments, the boundedness of $S_N$ is equivalent to the boundedness of its continuous version:
\[
S_Rf(x)=\int_{-R}^R\widehat f(\xi)e^{2\pi i\xi x}\,d\xi
=(1_{[-R,R]}\widehat f)^\vee(x).
\]

It is well known that the Hilbert transform is bounded, from which one can deduce that $S_R$ is also bounded.
\paragraph{Remark.} The transference principle is not always valid, especially when the operator contains nonlinear terms. For example,
\[
\overline S_Nf(x)=\sum_{|n|\le N}f(x-n^3)\frac1n.
\]
In addition, it can be used in ergodic systems.
\paragraph{Convergence almost everywhere pointwise.}
Carleson (65', $p=2$) and Hunt (67', $p>1$) prove that~\cite{Carleson1966,Hunt1968}
\[
\lim_{N\to\infty}S_Nf(x)=f(x)\quad\text{a.e. if }f\in L^p(\mathbb T).
\]
Their original proofs are very hard. In higher dimensions, even in $\mathbb R^2$, this is still open.
\subsection{Bochner--Riesz operator on \texorpdfstring{$\mathbb R^2$}{R2}}
Suppose

\[
S_Nf(x)=\sum_{\substack{|n|\le N\\n\in\mathbb Z^2}}\widehat f(n)e^{2\pi in\cdot x}.
\]
By the transference principle and scaling,
\[
\|S_N\|_{p\to p}<\infty\quad\Longleftrightarrow\quad
\|\mathcal S_1\|_{p\to p}<\infty,
\]
where $\mathcal S_1f(x):=(1_{B_1}\widehat f)^\vee(x)$ is the operator associated to the ball multiplier. Fefferman shows that, when $n\ge2$, $\|\mathcal S_1\|_{p\to p}<\infty$ if and only if $p=2$! See~\cite{Fefferman1971}.
On the other hand, if we replace $1_{B_1}$ by a smooth bump function $\varphi$, it is easy to check that $(\varphi\widehat f)^\vee$ is bounded for all $p>1$. Therefore, what happens if we add some smoothness to the multiplier?
Let
\[
m^\alpha(\xi)=(1-|\xi|^2)_+^\alpha.
\]
Note that $m^\alpha\to1_{B_1}$ as $\alpha\to0$, and the larger $\alpha$, the smoother $m^\alpha$.
\paragraph{Bochner--Riesz conjecture.}
If $n\ge2$, $p\ge1$, and $\alpha>\max\{n|\frac12-\frac1p|-\frac12,0\}$, we have

\[
\|m^\alpha(D)f\|_{L^p(\mathbb R^n)}\lesssim_{n,p,\alpha}\|f\|_{L^p(\mathbb R^n)},
\qquad m^\alpha(D)f:=(m^\alpha\widehat f)^\vee.
\]
\paragraph{Endpoints.}
\begin{align*}
p=2,&\quad\alpha>0&&\text{by Plancherel's theorem};\\
p=\infty,&\quad\alpha>\frac{n-1}{2},&&\|\check m^\alpha\|_1<\infty\ \text{iff }\alpha>\frac{n-1}{2};\\
p=\frac{2n}{n-1},&\quad\alpha>0,&&\text{and, by duality,}\\
p=\frac{2n}{n+1},&\quad\alpha>0.
\end{align*}
\begin{center}\includegraphics[width=.36\linewidth]{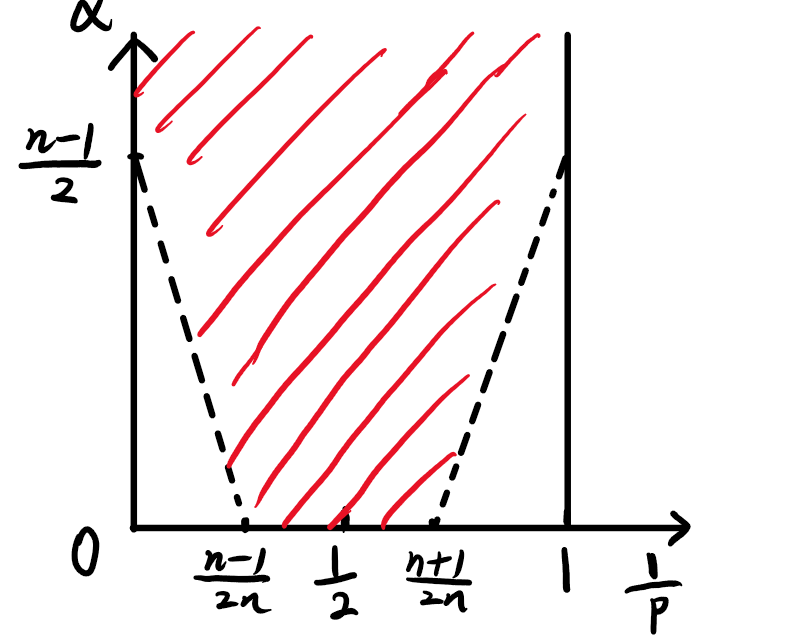}\end{center}
The planar Bochner--Riesz theorem is due to Carleson and Sj\"olin~\cite{CarlesonSjolin1972}.
\paragraph{Theorem 1.} When $n=2$ and $\alpha>\max\{2|\frac12-\frac1p|-\frac12,0\}$,
\[
\|m^\alpha(D)f\|_{L^p(\mathbb R^2)}\lesssim_{\alpha,p}\|f\|_{L^p(\mathbb R^2)}.
\]
In particular, for $p=4$, $\alpha>0$, we have
\[
\|m^\alpha(D)f\|_{L^4(\mathbb R^2)}\lesssim_\alpha\|f\|_{L^4(\mathbb R^2)}.
\]
\subsection{Reduction of the theorem}
Since all the mass concentrates on the neighborhood of the circle $S^1$, we decompose the unit ball $B_1$ into annuli in terms of the distance to $S^1$:

\[
1_{\{|\xi|\le1\}}(\xi)=\psi_0(|\xi|)+\sum_{\substack{k\ge1\\k\in\mathbb N}}\psi_k(1-|\xi|^2),
\]
where $\psi_0$ is supported on $\{|\xi|\le\frac34\}$, $\psi_k=\varphi_0(2^k\cdot)$, and $\varphi_0$ is supported on $\{\frac14<|\xi|<\frac34\}$. Thus
\begin{align*}
m^\alpha(\xi)
&=m^\alpha(\xi)1_{\{|\xi|\le1\}}(\xi)\\
&=m^\alpha(\xi)\psi_0(|\xi|)+\sum_km^\alpha(\xi)\psi_k(1-|\xi|^2)
=:m_0(\xi)+\sum_km_k(\xi).
\end{align*}
For $m_0$, $\|m_0(D)f\|_{p\to p}<\infty$, since $m_0$ is smooth.
For $m_k$, note that $m_k\sim2^{-\alpha k}\psi_k$; more precisely, $m_k=2^{-\alpha k}\widetilde\psi_k$ for a slightly different bump function $\widetilde\psi_k$. To simplify the notation, we still denote $\widetilde\psi_k$ as $\psi_k$. Hence it suffices to prove:
\paragraph{Theorem 2.} For $\alpha>\max\{2|\frac12-\frac1p|-\frac12,0\}$,
\[
\left\|\int_{\mathbb R^2}\widehat f(\xi)\psi_k(1-|\xi|^2)e^{ix\cdot\xi}\,d\xi\right\|_{L^p(\mathbb R^2)}
\lesssim2^{\alpha k}\|f\|_{L^p(\mathbb R^2)}.
\]
\paragraph{When $p=1$.}
Since $\psi_k(1-|\xi|^2)$ is supported on $|\xi|\sim2^{-k}$, we can decompose it into disjoint caps.

The caps $\theta$ have dimensions $2^{-k/2}\times2^{-k}$, and
\[
\psi_k(1-|\xi|^2)=\sum_\theta\varphi_\theta(\xi).
\]
\begin{center}\includegraphics[width=.26\linewidth]{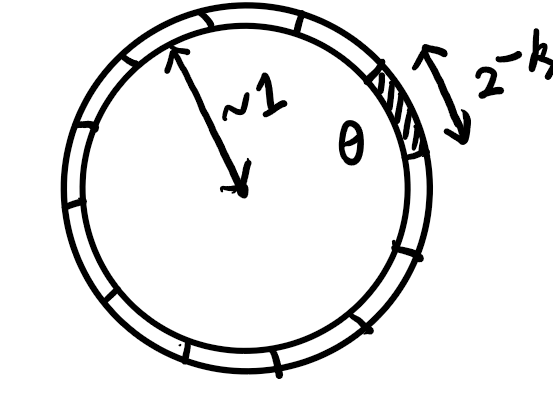}\end{center}
Then
\begin{align*}
\left\|\int_{\mathbb R^2}\widehat f(\xi)\psi_k(1-|\xi|^2)e^{ix\cdot\xi}\,d\xi\right\|_{L^1(\mathbb R^2)}
&=\left\|\sum_\theta(\check\varphi_\theta*f)\right\|_{L^1(\mathbb R^2)}\\
&\le\sum_\theta\|\check\varphi_\theta*f\|_{L^1(\mathbb R^2)}\\
&\le\sum_\theta\|\check\varphi_\theta\|_{L^1}\|f\|_{L^1}\\
&\lesssim2^{k/2}\|f\|_1\le2^{k\alpha}\|f\|_1,
\end{align*}
since $\check\varphi_\theta$ is $L^1$-normalized, for any $\alpha>\frac12$, as desired.
\paragraph{When $p=6$, $\alpha>\frac16$, using decoupling.}
Let $f_\theta=\check\varphi_\theta*f$. Then
\[
\left\|\sum_\theta f_\theta\right\|_{L^6(\mathbb R^2)}
\lesssim_\varepsilon2^{k\varepsilon}\left(\sum_\theta\|f_\theta\|_{L^6(\mathbb R^2)}^2\right)^{1/2}.
\]

Interpolate between
\[
\left(\sum_\theta\|f_\theta\|_2^2\right)^{1/2}\lesssim\|f\|_2
\quad\text{(Plancherel)},\qquad
\left(\sum_\theta\|f_\theta\|_\infty^2\right)^{1/2}\lesssim2^{k/4}\|f\|_\infty.
\]
We obtain
\[
\left(\sum_\theta\|f_\theta\|_6^2\right)^{1/2}\lesssim2^{k/6}\|f\|_6,
\]
and hence
\[
\left\|\int\widehat f(\xi)\psi_k(1-|\xi|^2)e^{ix\cdot\xi}\,d\xi\right\|
\lesssim_\varepsilon2^{(1/6+\varepsilon)k}\|f\|_6.
\]
Choose $\varepsilon=\alpha-\frac16>0$. We get the desired result.
\paragraph{Proof of Theorem 2 using restriction estimates (Carleson--Sj\"olin).} \cite{CarlesonSjolin1972}
WLOG, we prove the Bochner--Riesz operator for the parabola instead of the circle:
\[
\left\|\int_{\mathbb R^2}\widehat f(\xi)\chi_0(\xi_1)\psi_k(\xi_2-\xi_1^2)e^{ix\cdot\xi}\,d\xi\right\|_{L^4(\mathbb R^2)}
\lesssim_\varepsilon2^{\varepsilon k}\|f\|_{L^4(\mathbb R^2)}.
\]
Note that

\begin{align*}
&\int\widehat f(\xi)\chi_0(\xi_1)\psi_k(\xi_2-\xi_1^2)e^{ix\cdot\xi}\,d\xi\\
&\quad=\int_{\mathbb R^2}\int_{\mathbb R^2}f(y)\chi_0(\xi_1)\psi_k(\xi_2-\xi_1^2)e^{-iy\cdot\xi}e^{ix\cdot\xi}\,dy\,d\xi\\
&\quad=\int_{\mathbb R^2}f(y)\left(\int_{\mathbb R^2}\chi_0(\xi_1)\psi_k(\xi_2-\xi_1^2)e^{i(x-y)\cdot\xi}\,d\xi\right)dy.
\end{align*}
Denote
\[
K_k(x)=\int_{\mathbb R^2}\chi_0(\xi_1)\psi_k(\xi_2-\xi_1^2)e^{ix\cdot\xi}\,d\xi.
\]
The expression then equals $\int_{\mathbb R^2}f(y)K_k(x-y)\,dy=(K_k*f)(x)$.
Let $\xi_2-\xi_1^2\mapsto\xi_2$, $\xi_1\mapsto\xi_1$. We have $\frac{\partial(\xi_1,\xi_2)}{\partial(\eta_1,\eta_2)}=1$, and
\begin{align*}
K_k(x)&=\int_{\mathbb R^2}\chi_0(\xi_1)\psi_k(\xi_2)e^{ix\cdot(\xi_1,\xi_2+\xi_1^2)}\,d\xi\\
&=\int_{\mathbb R}\psi_k(\xi_2)e^{ix_2\xi_2}\,d\xi_2
 \cdot\int_{\mathbb R}\chi_0(\xi_1)e^{i[x\cdot(\xi_1,\xi_1^2)]}\,d\xi_1\\
&=\check\psi_k(x_2)E_{\mathbb P^1}\chi_0(x).
\end{align*}

By stationary phase, this is
\begin{align*}
K_k(x)&\sim\check\psi_k(x_2)|x_2|^{-1/2}e^{-ix_1^2/(4x_2)}
 \qquad\text{(we can assume }|x_1|\le10|x_2|\text{)}\\
&\sim2^{-k}1_{\{|x_2|\sim2^k\}}(x_2)|x_2|^{-1/2}e^{-ix_1^2/(4x_2)}\\
&\sim2^{-3k/2}e^{-ix_1^2/(4x_2)}1_{\{|x_2|\sim2^k\}}.
\end{align*}
Thus
\[
f*K_k(x)\approx\iint f(y)2^{-3k/2}e^{-i\frac{(x_1-y_1)^2}{4(x_2-y_2)}}
1_{\{|x_2|\sim2^k\}}(x_2-y_2)\,dy_1\,dy_2.
\]
By Minkowski's inequality, we have
\[
\|f*K_k\|_{L^4}\lesssim2^{-3k/2}\int_{\mathbb R}
\left\|\int_{\mathbb R}f(y_1,y_2)e^{-i\frac{(x_1-y_1)^2}{4(x_2-y_2)}}
1_{\{|x_2|\sim2^k\}}(x_2-y_2)\,dy_1\right\|_{L^4}\,dy_2.
\]
It suffices to prove
\[
\left\|\int_{\mathbb R}f(y_1,y_2)e^{-i\frac{(x_1-y_1)^2}{4(x_2-y_2)}}
1_{\{|x_2|\sim2^k\}}(x_2-y_2)\,dy_1\right\|_{L^4(\mathbb R^2)}
\lesssim_\varepsilon2^{\varepsilon k}2^{3k/4}\|f(\cdot,y_2)\|_{L^4(\mathbb R)}.
\]
Let $x_1\mapsto x_1$, $4(x_2-y_2)\mapsto x_2$. The integrand is
\[
\int_{\mathbb R}f(y_1,y_2)e^{-i(\frac{x_1^2}{x_2}-\frac{2x_1y_1}{x_2}+\frac{y_1^2}{x_2})}
1_{\{|x_2|\sim2^k\}}(x_2)\,dy_1.
\]
The term $x_1^2/x_2$ is crossed out in the handwritten phase, with the annotation ``since $|x_1|\le10|x_2|$''.

Let $2x_1/x_2=\widetilde x_1$, $-1/x_2=\widetilde x_2$. Then the Jacobian is $\sim|\widetilde x_2|^{-3}\sim2^{3k}$, and
\begin{align*}
&\left\|\int_{\mathbb R}f(y_1,y_2)e^{-i\frac{(x_1-y_1)^2}{4(x_2-y_2)}}
1_{\{|x_2|\sim2^k\}}(x_2-y_2)\,dy_1\right\|_{L^4}\\
&\qquad\approx2^{3k/4}\left\|\int_{\mathbb R}f(y_1,y_2)e^{i(\widetilde x_1y_1+\widetilde x_2y_1^2)}\,dy_1\right\|_{L^4(B_{2^k})}\\
&\qquad=2^{3k/4}\|E_{\mathbb P^1}(f(\cdot,y_2))\|_{L^4(B_{2^k})}\\
&\qquad\lesssim2^{\varepsilon k}2^{3k/4}\|f(\cdot,y_2)\|_{L^4(\mathbb R)}
\qquad\text{(restriction)}.
\end{align*}
Finally,
\[
\|f*K_k\|_{L^4(\mathbb R^2)}\lesssim2^{\varepsilon k}2^{-3k/4}
\int_{\mathbb R}\|f(\cdot,y_2)\|_4\,dy_2.
\]
Note that $|y_2|\sim2^k$. We get
\[
\|f*K_k\|_{L^4(\mathbb R^2)}
\le2^{\varepsilon k}\left(\int\|f(\cdot,y_2)\|_4^4\,dy_2\right)^{1/4}
=2^{\varepsilon k}\|f\|_{L^4(\mathbb R^2)}.\qquad\square
\]
\section{Proof of decoupling inequality for the parabola}
For the bilinear and asymmetric approaches developed in this section, see~\cite{Li2021,GLYZ2021}; the original decoupling theorem is~\cite{BD2015}.
\subsection*{Notation}
\[
\mathbb P^1=\{(\xi,\xi^2):\xi\in[0,1]\}.
\]
$\mathcal P(\delta)$ is the $\delta^2$-neighborhood of $\mathbb P^1$:
\[
\mathcal P(\delta)=\{(\xi,\xi^2+t):\xi\in[0,1],\ 0\le t\le\delta^2\}.
\]
We write
\[
\mathcal P(\delta)=\bigcup_{\ell(I)=\delta}\mathcal U_I,
\]
where $\mathcal U_I$ is the parallelogram of dimension $\delta\times\delta^2$.
\begin{center}\includegraphics[width=.40\linewidth]{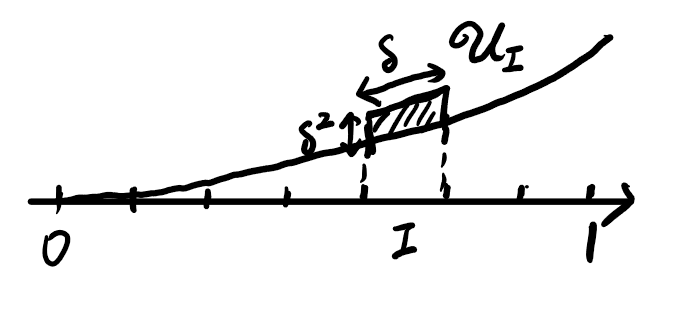}\end{center}
$D(\delta)$, the decoupling constant, is the smallest constant such that
\[
\left\|\sum_{\ell(I)=\delta}f_I\right\|_{L^6(\mathbb R^2)}
\le D(\delta)\left(\sum_I\|f_I\|_{L^6(\mathbb R^2)}^2\right)^{1/2}.
\]
\paragraph{Theorem 1.} $D(\delta)\le C_\varepsilon\delta^{-\varepsilon}$.

\paragraph{Main idea to prove Theorem 1.}
Bilinear estimates, Plancherel theorem and induction on scales.

\subsection*{Bilinear decoupling}
Take two intervals $J,J'\subset[0,1]$, $\ell(J)=\ell(J')=\frac14$ and
$\operatorname{dist}(J,J')\gtrsim1$.
$\mathcal P(J,\delta)$ denotes a partition of $J$ into disjoint $\delta$-intervals.
$B(\delta)$, the bilinear decoupling constant, is the smallest constant such that
\begin{align*}
&\left\|\left|\sum_{I\in\mathcal P(J,\delta)}f_I\cdot
\sum_{I'\in\mathcal P(J',\delta)}f_{I'}\right|^{1/2}\right\|_{L^6(\mathbb R^2)}\\
&\qquad\le B(\delta)
\left(\sum_{I\in\mathcal P(J,\delta)}\|f_I\|_{L^6(\mathbb R^2)}^2\right)^{1/4}
\left(\sum_{I'\in\mathcal P(J',\delta)}\|f_{I'}\|_{L^6(\mathbb R^2)}^2\right)^{1/4}.
\end{align*}
\subsection*{Relation between $D(\delta)$ and $B(\delta)$}
By H\"older's inequality, $B(\delta)\le D(\delta)$.
Conversely, $D(\delta)$ is also dominated by $B(\delta)$ up to a logarithmic loss:
\[
D(\delta)\lesssim(\log\delta^{-1})B(\delta).
\]
More precisely, we have the following proposition.
\paragraph{Proposition 2.} Take $\delta=2^{-N}$. Then
\[
D(\delta)\lesssim\left(1+\sum_{n=2}^{N}B(2^{-N+n-2})^2\right)^{1/2}.
\tag*{}
\]
\begin{proof}
Apply Whitney decomposition of $[0,1]^2$ at scale $\delta$.
Let $W_n$ consist of squares with side length $2^{-n}$.
\begin{center}\includegraphics[width=.36\linewidth]{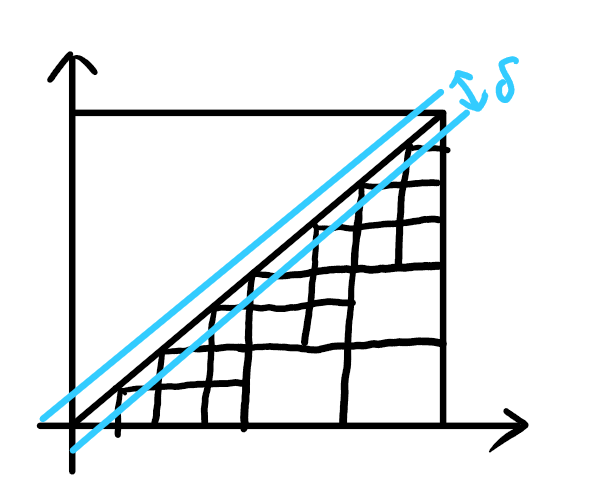}\end{center}
Then we have
\begin{align*}
\left\|\sum_{I\in\mathcal P(\delta)}f_I\right\|_{L^6(\mathbb R^2)}
&=\left\|\left(\sum_{I\in\mathcal P(\delta)}f_I\right)
\left(\sum_{I'\in\mathcal P(\delta)}\overline{f_{I'}}\right)\right\|_{L^3(\mathbb R^2)}^{1/2}\\
&=\left\|\sum_{n=2}^N\sum_{J\times J'\in W_n}f_Jf_{J'}\right\|_{L^3(\mathbb R^2)}^{1/2}\\
&\le\left(\sum_{n=2}^N\sum_{J\times J'\in W_n}\|f_Jf_{J'}\|_{L^3(\mathbb R^2)}\right)^{1/2},
\end{align*}
where $f_J=\sum_{I\in\mathcal P(J,\delta)}f_I$, and assume $f$ is real-valued.
By parabolic rescaling and the definition of $B(\delta)$, we obtain
\begin{align*}
\left\|\sum_{I\in\mathcal P(\delta)}f_I\right\|_6
\le\Bigg[&\sum_{J\times J'\in W_N}\|f_Jf_{J'}\|_3\\
&+\sum_{n=2}^{N-1}\sum_{J\times J'\in W_n}
 B\left(\frac{2^{-N}}{2^{-n}}\right)^2
 \left(\sum_{I\in\mathcal P(J,\delta)}\|f_I\|_6^2\right)^{1/2}
 \left(\sum_{I'\in\mathcal P(J',\delta)}\|f_{I'}\|_6^2\right)^{1/2}\Bigg]^{1/2}.
\end{align*}
Using Cauchy--Schwarz inequality, and noting that for a fixed $J$,
\[
\#\{J':J\times J'\in W_n\}\le3\qquad\text{for all }n,
\]
we get
\begin{align*}
\left\|\sum_{I\in\mathcal P(\delta)}f_I\right\|_6
&\lesssim\left(\sum_{I\in\mathcal P(\delta)}\|f_I\|_6^2
 +\sum_{n=2}^{N-1}B(2^{-N+n})^2
 \sum_{I\in\mathcal P(\delta)}\|f_I\|_6^2\right)^{1/2}\\
&=\left(1+\sum_{n=2}^{N-1}B(2^{-N+n})^2\right)^{1/2}
 \left(\sum_{I\in\mathcal P(\delta)}\|f_I\|_6^2\right)^{1/2},
\end{align*}
as desired.
\end{proof}
\paragraph{Remark.} Through precise calculation, Zane Li proves that~\cite{Li2021}
\[
D(\delta)\lesssim2^{\frac{\log\delta^{-1}}{\log\log\delta^{-1}}}.
\]
Since
\[
\log\delta^{-1}\le2^{\log\delta^{-1}/\log\log\delta^{-1}}
\le C_\varepsilon\delta^{-\varepsilon},
\]
the Bourgain--Demeter method cannot improve $D(\delta)$ to $\log\delta^{-1}$.

\subsection*{Induction on scales}
We will do induction on $a,b$.
\begin{center}\includegraphics[width=.52\linewidth]{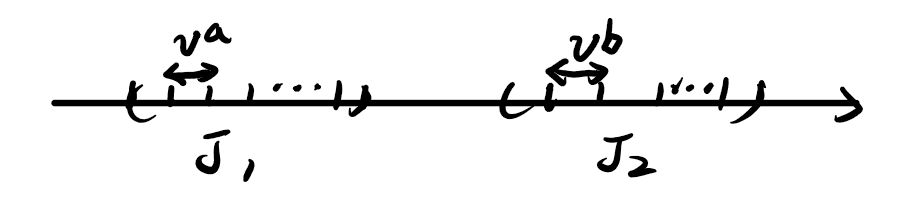}\end{center}
To do induction, the key step is the following lemma.
Fix $J_1,J_2\subset[0,1]$, $\ell(J_1)=\ell(J_2)=\nu^a$ and
$\operatorname{dist}(J_1,J_2)\gtrsim1$.
Let $B_a(\delta)$ be the smallest constant such that
\[
\big\||f_{J_1}f_{J_2}|^{1/2}\big\|_{L^6(\mathbb R^2)}
\le B_a(\delta)
\left(\sum_{I_1\in\mathcal P(J_1,\delta)}\|f_{I_1}\|_6^2\right)^{1/4}
\left(\sum_{I_2\in\mathcal P(J_2,\delta)}\|f_{I_2}\|_6^2\right)^{1/4}.
\]
\paragraph{Lemma 3.}
\[
B_a(\delta)\lesssim B_{2a}(\delta)^{1/2}D\left(\frac\delta{\nu^a}\right)^{1/2}.
\]
\begin{proof}
Here we give an intuitive but not rigorous proof, and in fact we will show
\[
B_a(\delta)\lesssim B_{2a}(\delta)^{1/3}D\left(\frac\delta{\nu^a}\right)^{2/3}.
\]
By Cauchy--Schwarz inequality,
\begin{align*}
\big\||f_{J_1}f_{J_2}|^{1/2}\big\|_{L^6(\mathbb R^2)}
&=\left(\int|f_{J_1}|^3|f_{J_2}|^3\right)^{1/6}\\
&\le\left(\int|f_{J_1}|^2|f_{J_2}|^4\right)^{1/12}
\left(\int|f_{J_1}|^4|f_{J_2}|^2\right)^{1/12}.
\end{align*}
By symmetry, it suffices to estimate the first term.
Since $\operatorname{supp}(\widehat f_{J_2})\subset\mathcal U_{J_2}$,
we know $|f_{J_2}|$ is essentially constant on each tube $T_{J_2}$,
denoted by $C_{T_{J_2}}$, by the local constant property.
\begin{center}\includegraphics[width=.70\linewidth]{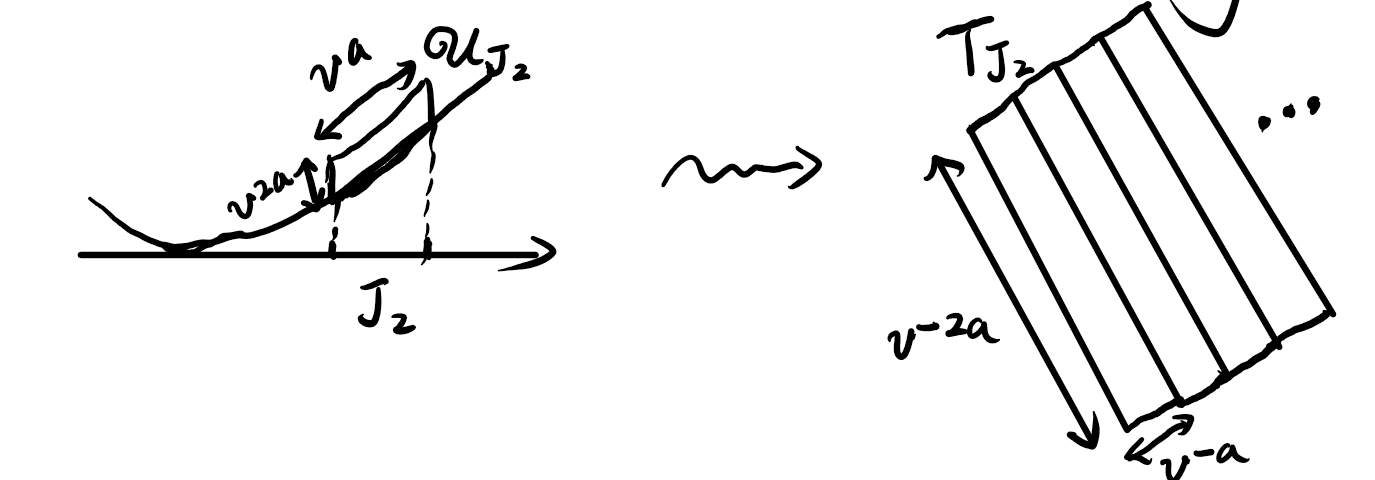}\end{center}
So we have
\[
\int|f_{J_1}|^2|f_{J_2}|^4
=\sum_{T_{J_2}}\int_{T_{J_2}}|f_{J_1}|^2|f_{J_2}|^4
\approx\sum_{T_{J_2}}C_{T_{J_2}}^4\int_{T_{J_2}}|f_{J_1}|^2.
\]
We claim that, for any line $L$ whose direction is normal to $\mathcal U_{J_2}$,
\[
\int_{L\cap T_{J_2}}|f_{J_1}|^2
\le\sum_{J_1'\in\mathcal P(J_1,\nu^{2a})}
\int_{L\cap T_{J_2}}|f_{J_1'}|^2.
\]
\begin{center}\includegraphics[width=.37\linewidth]{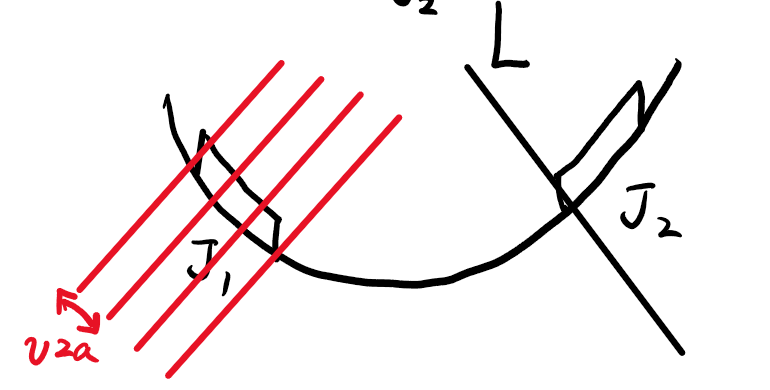}\end{center}
In fact, the Fourier transform of $f_{J_1'}$ along $L$ is supported on the
projection of $\operatorname{supp}(\widehat f_{J_1'})$ onto $L$.
Write $\mathcal F_2f(x,\eta)$ for the Fourier transform of $f$ in the second variable only.
More precisely,
\begin{align*}
\int_{\mathbb R}f_{J_1'}(x,y)e^{-ix\xi}\,dx
&=\int(\mathcal F_2f_{J_1'})(x,\eta)e^{iy\eta}e^{-ix\xi}\,d\eta\,dx\\
&=\int\widehat f_{J_1'}(\xi,\eta)e^{iy\eta}\,d\eta.
\end{align*}
Thus
\[
\operatorname{supp}\big(\widehat{f_{J_1'}|_L}\big)
\subset\pi_L\big(\operatorname{supp}(\widehat f_{J_1'})\big)
\subset\pi_L(\mathcal U_{J_1'});
\]
see also Section 5.
Recall $\operatorname{dist}(J_1,J_2)\gtrsim1$; it follows that
$\theta_{J_1},\theta_{J_2}$ are transverse. Hence
$\pi_L(\mathcal U_{J_1'})$ is contained in an interval of length $\sim\nu^{2a}$.
And since $L\cap T_{J_2}$ has length $\sim\nu^{-2a}$, the decomposition
$\sum_{J_1'}f_{J_1'}$ is $L^2$-orthogonal, i.e.,
\[
\left\|\sum_{J_1'}f_{J_1'}\right\|_{L^2(L\cap T_{J_2})}^2
\lesssim\sum_{J_1'}\|f_{J_1'}\|_{L^2(L\cap T_{J_2})}^2,
\]
as desired.
Therefore, by Fubini's theorem, we have
\begin{align*}
\int|f_{J_1}|^2|f_{J_2}|^4
&\lesssim\sum_{T_{J_2}}C_{T_{J_2}}^4\sum_{J_1'}\int_{T_{J_2}}|f_{J_1'}|^2\\
&\sim\sum_{J_1'}\int|f_{J_1'}|^2|f_{J_2}|^4\\
&\le\sum_{J_1'}\left(\int|f_{J_1'}|^4|f_{J_2}|^2\right)^{1/2}
\left(\int|f_{J_2}|^6\right)^{1/2}.
\end{align*}
Here the integrand in the second line is written as
$|f_{J_1'}|^2|f_{J_2}|\cdot|f_{J_2}|^3$.
By the same argument as before, but instead of using scale $\nu^{4a}$ to
partition $J_2$, we use $J_1'\subset B_{\nu^{2a}}$ and squares of side length
$\nu^{-2a}$ tiling $\mathbb R^2$ directly.
\begin{center}\includegraphics[width=.40\linewidth]{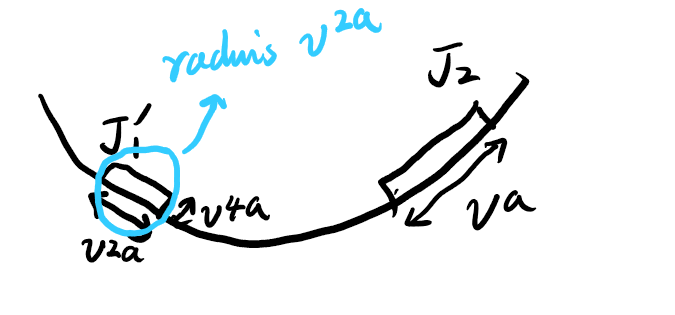}\end{center}
We obtain
\[
\int|f_{J_1}|^2|f_{J_2}|^4
\lesssim\sum_{J_1'}
\left(\sum_{J_2'\in\mathcal P(J_2,\nu^{2a})}
\int|f_{J_1'}|^4|f_{J_2'}|^2\right)^{1/2}
\left(\int|f_{J_2}|^6\right)^{1/2}.
\]
Write $|f_{J_1'}|^4|f_{J_2'}|^2=|f_{J_1'}|^2|f_{J_2'}|^2\cdot|f_{J_1'}|^2$.
Then the preceding expression is bounded by
\begin{align*}
&\sum_{J_1'}\left(\sum_{J_2'}
\left(\int|f_{J_1'}|^3|f_{J_2'}|^3\right)^{2/3}
\left(\int|f_{J_1'}|^6\right)^{1/3}\right)^{1/2}
 D\left(\frac\delta{\nu^a}\right)^3
 \left(\sum_{I_2\in\mathcal P(J_2,\delta)}\|f_{I_2}\|_6^2\right)^{3/2}\\
&\le\sum_{J_1'\subset J_1}\left(\sum_{J_2'\subset J_2}
 \left[B_{2a}(\delta)
 \left(\sum_{I_1\subset J_1'}\|f_{I_1}\|_6^2\right)^{1/4}
 \left(\sum_{I_2\subset J_2'}\|f_{I_2}\|_6^2\right)^{1/4}\right]^{(2/3)\times6}\right)^{1/2}\\
&\qquad\times D\left(\frac\delta{\nu^{2a}}\right)
 D\left(\frac\delta{\nu^a}\right)^3
 \left(\sum_{I_2\subset J_2}\|f_{I_2}\|_6^2\right)^{3/2}
 \left(\sum_{I_1\subset J_1'}\|f_{I_1}\|_6^2\right)^{1/2}\\
&\le B_{2a}(\delta)^2D\left(\frac\delta{\nu^a}\right)^4
 \left(\sum_{J_1'\subset J_1}\sum_{I_1\subset J_1'}\|f_{I_1}\|_6^2\right)^{1/2}
 \left(\sum_{J_2'\subset J_2}\sum_{I_2\subset J_2'}\|f_{I_2}\|_6^2\right)^{1/2}\\
&\qquad\times
 \left(\sum_{J_1'\subset J_1}\sum_{I_1\subset J_1'}\|f_{I_1}\|_6^2\right)^{1/2}
 \left(\sum_{I_2\subset J_2}\|f_{I_2}\|_6^2\right)^{3/2}\\
&=B_{2a}(\delta)^2D\left(\frac\delta{\nu^a}\right)^4
 \left(\sum_{I_1\subset J_1}\|f_{I_1}\|_6^2\right)
 \left(\sum_{I_2\subset J_2}\|f_{I_2}\|_6^2\right)^2.
\end{align*}
In the summation over $J_2'$, sum directly; in the summation over $J_1'$, use
Cauchy--Schwarz.
In summary, we have
\[
\big\||f_{J_1}f_{J_2}|^{1/2}\big\|_6
\lesssim B_{2a}(\delta)^{1/3}D\left(\frac\delta{\nu^a}\right)^{2/3}
\left(\sum_{I_1\subset J_1}\|f_{I_1}\|_6^2\right)^{1/4}
\left(\sum_{I_2\subset J_2}\|f_{I_2}\|_6^2\right)^{1/4}.
\]
This implies
\[
B_a(\delta)\lesssim B_{2a}(\delta)^{1/3}D\left(\frac\delta{\nu^a}\right)^{2/3}.
\]
\end{proof}
\subsection*{Asymmetric bilinear constant}
For the asymmetric constants and smoothing weights, see~\cite{GLYZ2021}.
Let $a,b\in[0,1]$. Let $B_{a,b}(\delta)$ be the smallest constant such that,
for any intervals $I\in\mathcal P(\delta^a)$, $I'\in\mathcal P(\delta^b)$ and
$\operatorname{dist}(I,I')\ge\frac14$,
\begin{align*}
\int (|F_I|^2*\phi_I)(|F_{I'}|^4*\phi_{I'})
\le B_{a,b}(\delta)^6
\left(\sum_{J\in\mathcal P(I,\delta)}\|F_J\|_6^2\right)^{2/2}
\left(\sum_{J'\in\mathcal P(I',\delta)}\|F_{J'}\|_6^2\right)^{4/2},
\end{align*}
where $\widehat\phi_I=1$ on $\mathcal U_I$, $\phi$ is positive and
$\|\phi\|_{L^1}=1$.
\paragraph{Lemma 4.} $|F_I|^p\lesssim_p|F_I|^p*\phi_I$.
\begin{proof}
Since $\operatorname{supp}(\widehat F_I)\subset\mathcal U_I$,
$\widehat F_I=\widehat F_I\cdot\widehat\phi_I$, hence $F_I=F_I*\phi_I$.
Therefore,
\begin{align*}
|F_I(x)|^p
&=\left|\int F_I(x-y)\phi_I(y)\,dy\right|^p\\
&\le\left(\int|F_I(x-y)|^p\phi_I(y)\,dy\right)
 \left(\int\phi_I(y)\,dy\right)^{p-1}.
\end{align*}
Thus $|F_I(x)|^p\lesssim_p(|F_I|^p*\phi_I)(x)$.
\end{proof}
Applying Plancherel theorem, we have the following.
\paragraph{Lemma 5.} Let $\sigma\in2^{-\mathbb N}$, $\ell(I)=\sigma$.
Suppose $f_I:\mathbb R\to\mathbb C$, $\operatorname{supp}(\widehat f_I)\subset I$.
Then for any interval $B\subset\mathbb R$ of length $2/\sigma$, we have
\[
\frac1{\ell(B)}\int_B\left|\sum_{I\in\mathcal P(\sigma)}f_I\right|^2
\lesssim_M\sum_{I\in\mathcal P(\sigma)}\int|f_I|^2\Phi_B,
\]
where
\[
\Phi_B(x)=\frac1{\ell(B)}(1+\sigma\operatorname{dist}(x,B))^{-M}.
\]
\paragraph{Lemma 6.} If $0\le a\le2b\le1$, $I\in\mathcal P(\delta^a)$,
$I'\in\mathcal P(\delta^b)$, $\operatorname{dist}(I,I')\ge\frac14$, then
\[
\int_{\mathbb R^2}(|F_I|^2*\phi_I)(|F_{I'}|^4*\phi_{I'})
\lesssim\sum_{I_-\in\mathcal P(I,\delta^{2b})}
\int(|F_{I_-}|^2*\phi_{I_-})(|F_{I'}|^4*\phi_{I'}).
\]
\begin{proof}
Without loss of generality, assume $I'=[-\delta^b/2,\delta^b/2]$.
Let $L$ be the $y$-axis, $L_z=L+z$.
\begin{center}\includegraphics[width=.28\linewidth]{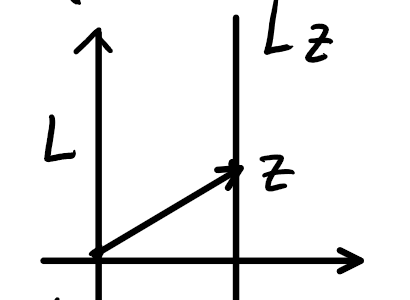}\end{center}
By Fubini's theorem,
\begin{align*}
&\int_{\mathbb R^2}(|F_I|^2*\phi_I)(|F_{I'}|^4*\phi_{I'})\\
&\qquad=\frac1{\delta^{-2b}}\int_{\mathbb R^2}
\int_{B_{L_z}(z,\delta^{-2b})}(|F_I|^2*\phi_I)(|F_{I'}|^4*\phi_{I'}).
\end{align*}
Note the local constant property:
\[
\sup_{x\in B_{L_z}(z,\delta^{-2b})}(|F_{I'}|^4*\phi_{I'})(x)
\lesssim(|F_{I'}|^4*\phi_{I'})(z).
\]
The preceding expression is therefore bounded by
\[
\frac1{\delta^{-2b}}\int_{\mathbb R^2}(|F_{I'}|^4*\phi_{I'})(z)
\left(\int_{B_{L_z}(z,\delta^{-2b})}|F_I|^2*\phi_I\right)\,dz.
\]
Suppose $z=(z_1,z_2)$. We have
\begin{align*}
\int_{B_{L_z}(z,\delta^{-2b})}|F_I|^2*\phi_I
&=\int_{|z_2'-z_2|\le\delta^{-2b}}(|F_I|^2*\phi_I)(z_1,z_2')\,dz_2'\\
&=\int_{|z_2'-z_2|\le\delta^{-2b}}\int_{\mathbb R^2}
 |F_I(z_1-w_1,z_2'-w_2)|^2\phi_I(w_1,w_2)\,dw_1\,dw_2\,dz_2'.
\end{align*}
Make the substitutions
\[
z_2'-w_2\mapsto z_2',\qquad z_1-w_1\mapsto w_1,\qquad z_2-w_2\mapsto w_2.
\]
Then the expression becomes
\begin{align*}
&\int_{\mathbb R^2}\left(\int_{|z_2'-w_2|\le\delta^{-2b}}
 |F_I(w_1,z_2')|^2\,dz_2'\right)
 \phi_I(z_1-w_1,z_2-w_2)\,dw_1\,dw_2\\
&\lesssim\int_{\mathbb R^2}\left(\delta^{-2b}
 \sum_{I_-\in\mathcal P(I,\delta^{2b})}
 \|F_{I_-}\|_{L^2(\Phi_{B_{L_w}(w,\delta^{-2b})})}^2\right)
 \phi_I(z_1-w_1,z_2-w_2)\,dw_1\,dw_2,
\end{align*}
by Lemma 5.
Substituting this back, we obtain
\begin{align*}
&\sum_{I_-\in\mathcal P(I,\delta^{2b})}\int_{\mathbb R^2}
 (|F_{I'}|^4*\phi_{I'})(z)
 \left(\int_{\mathbb R^2}\|F_{I_-}\|_{L^2(\Phi_{B_{L_w}(w,\delta^{-2b})})}^2
 \phi_I(z-w)\,dw\right)\,dz\\
&=\sum_{I_-\in\mathcal P(I,\delta^{2b})}\int_{\mathbb R^2}
 (|F_{I'}|^4*\phi_{I'})(z)
 \left(|F_{I_-}|^2*_L\Phi_{B_L(0,\delta^{-2b})}*\phi_I\right)(z)\,dz.
\end{align*}
Here $*_L$ means the one-dimensional convolution, i.e.,
\[
\big(|F_{I_-}|^2*_L\Phi_{B_L(0,\delta^{-2b})}\big)(w)
=\int_{\mathbb R}|F_{I_-}(w_1,t)|^2\Phi_{B_L(0,\delta^{-2b})}(w_2-t)\,dt.
\]
Recall $\langle f,g*h\rangle=\langle f*g,h\rangle$. We have
\begin{align*}
&\sum_{I_-\in\mathcal P(I,\delta^{2b})}\int_{\mathbb R^2}
 \big(|F_{I'}|^4*\phi_{I'}*_L\Phi_{B_L(0,\delta^{-2b})}\big)
 \big(|F_{I_-}|^2*\phi_I\big)\\
&\lesssim\sum_{I_-\in\mathcal P(I,\delta^{2b})}\int
 (|F_{I'}|^4*\phi_{I'})(|F_{I_-}|^2*\phi_{I_-}).
\end{align*}
For the second factor, Lemma 4 with $p=2$ gives
\[
|F_{I_-}|^2*\phi_I\lesssim|F_{I_-}|^2*\phi_{I_-}*\phi_I.
\]
\end{proof}
To conclude the above proof, we need the following lemma.
\paragraph{Lemma 7.}
\[
\phi_{I'}*_L\Phi_{B_L(0,\delta^{-2b})}\lesssim\phi_{I'},
\qquad \phi_{I_-}*\phi_I\lesssim\phi_{I_-}.
\]
\begin{proof}
The proof is illustrated by the essential supports below.
\begin{center}\includegraphics[width=.87\linewidth]{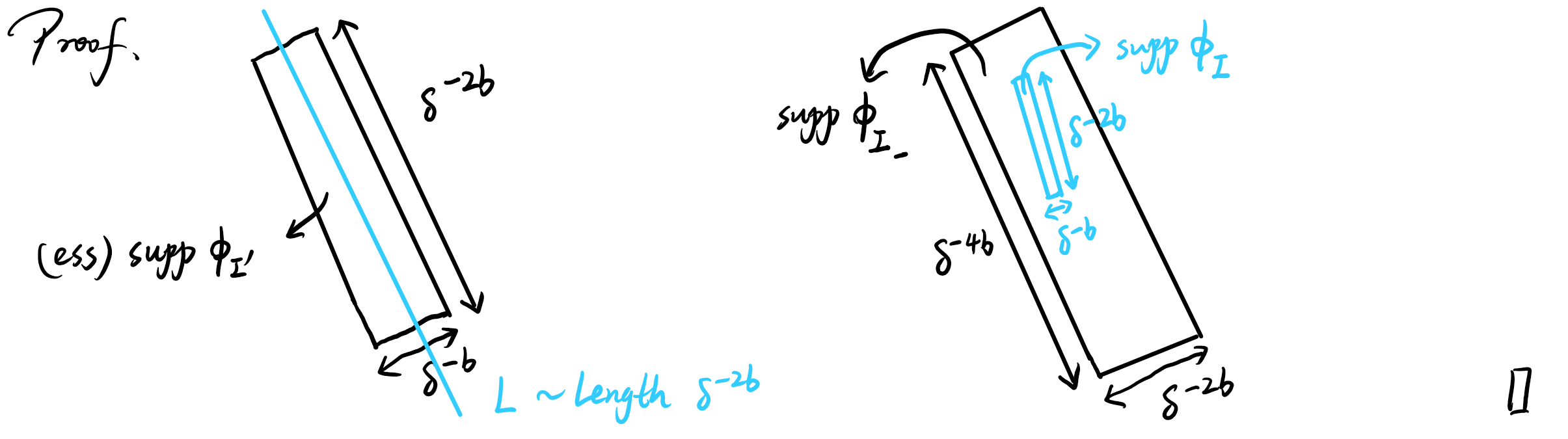}\end{center}
The line segment $L$ has length $\sim\delta^{-2b}$.
The essential support of $\phi_{I'}$ has dimensions $\delta^{-b}\times\delta^{-2b}$;
that of $\phi_{I_-}$ has dimensions $\delta^{-2b}\times\delta^{-4b}$.
\end{proof}
\paragraph{Corollary 8.} $B_{a,b}(\delta)\le B_{2b,b}(\delta)$ for $a\le2b\le1$.

\subsection*{H\"older's inequality}
\paragraph{Proposition 9.} Assume $b\le\frac14$. We have
\[
B_{2b,b}(\delta)\lesssim D\left(\frac\delta{\delta^b}\right)^{1/2}B_{4b,2b}(\delta)^{1/2}.
\]
\begin{proof}
Since
\begin{align*}
(|F_I|^2*\phi_I)(x)
&=\int|F_I(x-y)|^2\phi_I(y)\,dy\\
&\le\left(\int|F_I(x-y)|^4\phi_I(y)\,dy\right)^{1/2}
 \left(\int\phi_I(y)\,dy\right)^{1/2}\\
&\lesssim\big((|F_I|^4*\phi_I)(x)\big)^{1/2},
\end{align*}
and
\begin{align*}
(|F_{I'}|^4*\phi_{I'})(x)
&=\int|F_{I'}(x-y)|^4\phi_{I'}(y)\,dy\\
&\le\left(\int|F_{I'}(x-y)|^2\phi_{I'}(y)\,dy\right)^{1/2}
\left(\int|F_{I'}(x-y)|^6\phi_{I'}(y)\,dy\right)^{1/2}\\
&=\big((|F_{I'}|^2*\phi_{I'})(x)\big)^{1/2}\big((|F_{I'}|^6*\phi_{I'})(x)\big)^{1/2},
\end{align*}
we have
\begin{align*}
\int(|F_I|^2*\phi_I)(|F_{I'}|^4*\phi_{I'})
&\le\int(|F_I|^4*\phi_I)^{1/2}(|F_{I'}|^2*\phi_{I'})^{1/2}
 (|F_{I'}|^6*\phi_{I'})^{1/2}\\
&\le\left(\int|F_{I'}|^6*\phi_{I'}\right)^{1/2}
 \left(\int(|F_{I'}|^2*\phi_{I'})(|F_I|^4*\phi_I)\right)^{1/2}.
\end{align*}
The first factor is $\lesssim(\|F_{I'}\|_6^6)^{1/2}$;
in the second factor, $\ell(I')=\delta^b$ and $\ell(I)=\delta^{2b}$.
Consequently,
\begin{align*}
B_{2b,b}(\delta)
&\lesssim D\left(\frac\delta{\delta^b}\right)^{1/2}B_{b,2b}(\delta)^{1/2}\\
&\lesssim D\left(\frac\delta{\delta^b}\right)^{1/2}B_{4b,2b}(\delta)^{1/2},
\end{align*}
by Corollary 8.
\end{proof}
\subsection*{Bootstrapping argument}
Let $\eta\ge0$ be the infimum over all $\varepsilon>0$ such that
\[
D(\delta)\le C_\varepsilon\delta^{-\varepsilon}
\qquad\text{for any }\delta\in(0,1).
\]
Let $A(b)$ be the infimum over all $A>0$ such that
\[
B_{2b,b}(\delta)\le C_A\delta^{-A}
\qquad\text{for any }\delta\in(0,1).
\]
Now we know
\[
A(b)\le\frac12\eta(1-b)+\frac12 A(2b),
\]
so
\[
\frac{2\eta-2A(b)}b\ge2\frac{\eta-A(2b)}{2b}+\eta.
\]
Define
\[
A=\liminf_{b\to0}\frac{\eta-A(b)}b.
\]
If $|A|<\infty$, taking $b\to0$ we obtain
\[
2A\ge2A+\eta\quad\Longrightarrow\quad\eta=0.
\]
This completes our proof of Theorem 1. So it is sufficient to prove $|A|<\infty$.
\begin{proof}
We start with a crude upper bound.
Since $\#\mathcal P(\delta^a)\sim\delta^{-a}$,
$\#\mathcal P(\delta^b)\sim\delta^{-b}$, we have
\[
B(\delta)\lesssim\delta^{-a/6}\delta^{-b/2}B_{a,b}(\delta).
\tag{$*$}
\]
This implies $\eta\le5b/6+A(b)$ if we choose $a=2b$.
It follows that $A\le5/6$.
To prove $(*)$, consider
\begin{align*}
\int|F_I|^3|F_{I'}|^3
&\le\left(\int|F_I|^2|F_{I'}|^4\right)^{1/2}
 \left(\int|F_I|^4|F_{I'}|^2\right)^{1/2}\\
&=\left(\int\left|\sum_{J\in\mathcal P(I,\delta^a)}F_J\right|^2
 \left|\sum_{J'\in\mathcal P(I',\delta^b)}F_{J'}\right|^4\right)^{1/2}\\
&\qquad\times\left(\int\left|\sum_{J'\in\mathcal P(I',\delta^a)}F_{J'}\right|^2
 \left|\sum_{J\in\mathcal P(I,\delta^b)}F_J\right|^4\right)^{1/2}.
\end{align*}
For the first term, we have
\begin{align*}
\text{First term}
&\lesssim\delta^{-a}\delta^{-3b}
 \sum_{J\in\mathcal P(I,\delta^a)}\sum_{J'\in\mathcal P(I',\delta^b)}
 \int|F_J|^2|F_{J'}|^4\\
&\lesssim\delta^{-a}\delta^{-3b}
 \sum_{J\in\mathcal P(I,\delta^a)}\sum_{J'\in\mathcal P(I',\delta^b)}
 \int(|F_J|^2*\phi_J)(|F_{J'}|^4*\phi_{J'}).
\end{align*}
Continuing,
\begin{align*}
\text{First term}
&\le\delta^{-a}\delta^{-3b}B_{a,b}(\delta)^6\sum_J\sum_{J'}
 \left(\sum_{J_-\in\mathcal P(J,\delta)}\|F_{J_-}\|_6^2\right)^{2/2}
 \left(\sum_{J_-'\in\mathcal P(J',\delta)}\|F_{J_-'}\|_6^2\right)^{4/2}\\
&\le\delta^{-a}\delta^{-3b}B_{a,b}(\delta)^6
 \left(\sum_J\sum_{J_-}\|F_{J_-}\|_6^2\right)
 \left(\sum_{J'}\sum_{J_-'}\|F_{J_-'}\|_6^2\right)^2\\
&=\delta^{-a}\delta^{-3b}B_{a,b}(\delta)^6
 \left(\sum_{I_-\in\mathcal P(I,\delta)}\|F_{I_-}\|_6^2\right)
 \left(\sum_{I_-'\in\mathcal P(I',\delta)}\|F_{I_-'}\|_6^2\right)^2.
\end{align*}
By symmetry, we know
\[
\text{Second term}\lesssim\delta^{-a}\delta^{-3b}B_{a,b}(\delta)^6
\left(\sum_{I_-'\in\mathcal P(I',\delta)}\|F_{I_-'}\|_6^2\right)
\left(\sum_{I_-\in\mathcal P(I,\delta)}\|F_{I_-}\|_6^2\right)^2.
\]
In summary, we get
\[
\int|F_I|^3|F_{I'}|^3
\lesssim\delta^{-a}\delta^{-3b}B_{a,b}(\delta)^6
\left(\sum_{I_-}\|F_{I_-}\|_6^2\right)^{3/2}
\left(\sum_{I_-'}\|F_{I_-'}\|_6^2\right)^{3/2},
\]
which implies
\[
B(\delta)\le\delta^{-a/6}\delta^{-b/2}B_{a,b}(\delta).
\]
To obtain a lower bound for $A$, we choose $\ell(I)=\delta^{2b}$,
$\ell(I')=\delta^b$, $\operatorname{dist}(I,I')\gtrsim1$, and consider
\begin{align*}
\int(|F_I|^2*\phi_I)(|F_{I'}|^4*\phi_{I'})
&\lesssim\left(\int(|F_I|^2*\phi_I)^3\right)^{1/3}
 \left(\int(|F_{I'}|^4*\phi_{I'})^{3/2}\right)^{2/3}\\
&\lesssim\|F_I\|_6^2\|F_{I'}\|_6^4,
\end{align*}
by Young's inequality. Hence
\[
B_{2b,b}(\delta)\le D(\delta^{1-2b})^{1/3}D(\delta^{1-b})^{2/3},
\]
so
\[
A(b)\le\frac13\eta(1-2b)+\frac23\eta(1-b).
\]
Thus
\[
\frac{\eta-A(b)}b\ge\frac43\eta\ge0,
\qquad\text{and hence}\qquad A\ge0.
\]
\end{proof}

\clearpage
\phantomsection

\end{document}